\documentclass[a4paper, draft, reqno]{amsart}
\usepackage{epsf}
\usepackage{amsmath}
\usepackage{amssymb}
\usepackage{latexsym}

\usepackage[cmtip,matrix,arrow]{xy}
\UseComputerModernTips
\CompileMatrices

\newcommand{\zed}{{\mathbb Z}}
\newcommand{\RR}{{\mathbb R}}
\newcommand{\NN}{{\mathbb N}}
\newcommand{\pf}{\noindent\textbf{Proof:} }
\DeclareMathOperator{\Hom}{Hom}
\DeclareMathOperator{\Ext}{Ext}

\DeclareMathOperator{\soc}{soc}
\DeclareMathOperator{\im}{im}
\DeclareMathOperator{\hd}{hd}
\DeclareMathOperator{\ch}{ch}

\newcommand{\frob}{^{\mbox{\rm\tiny F}}}

\newcommand{\otherwise}{\mbox{\rm otherwise}}
\newcommand{\wif}{\mbox{\rm if }}
\newcommand{\wand}{\mbox{\rm and }}

\renewcommand{\mod}[0]{\mbox{\rm mod }}

\newcommand{\GL}{\mathrm{GL}}
\newcommand{\SL}{\mathrm{SL}}
\newcommand{\St}{\mathrm{St}}
\newcommand{\notequiv}{\not\equiv}

\DeclareMathOperator{\ind}{ind}

\DeclareMathOperator{\pr}{pr}
\DeclareMathOperator{\Stab}{Stab}

\newcommand{\alphac}{\alpha\check{\ }\,}

\begin{document}
\theoremstyle{plain}
\newtheorem{theorem}{Theorem}[section]
\newtheorem{fakeprop}{Proposition}[section]
\newtheorem{proposition}[theorem]{Proposition}
\newtheorem{lemma}[theorem]{Lemma}
\newtheorem{corollary}[theorem]{Corollary}
\newtheorem{conjecture}[theorem]{Conjecture}
\theoremstyle{definition}
\newtheorem{remark}[theorem]{Remark}
\newtheorem{ass}[theorem]{Assumption}
\newtheorem{definition}[theorem]{Definition}
\newtheorem{example}[theorem]{Example}


\setlength{\parskip}{1ex}


\title{Homomorphisms between Weyl modules for $\SL_3(k)$}

\author{Anton Cox
\and Alison Parker} 
\thanks{The first author was partially supported by Nuffield grant scheme NUF-NAL 02}  
\address{Centre for Mathematical Science \\
City University  \\
Northampton Square, London, EC1V 0HB\\
 England.} 
\email{A.G.Cox@city.ac.uk}
\address{
School of Mathematics and Statistics F07\\
University of Sydney \\
NSW 2006\\
Australia.}
\email{alisonp@maths.usyd.edu.au}

\begin{abstract}
We classify all homomorphisms between
Weyl modules for $\SL_3(k)$ when $k$ is an algebraically closed field
of characteristic at least three, and show that the $\Hom$-spaces are all
at most one-dimensional. As a corollary we obtain all homomorphisms
between Specht modules for the symmetric group when the labelling
partitions have at most three parts and the prime is at least
three. 
We
conclude by showing how a result of Fayers and Lyle on Hom-spaces for
Specht modules is related to earlier work of Donkin for algebraic
groups.
\end{abstract}
\maketitle



\section{Introduction}

Let $G$ be a reductive algebraic group over an algebraically closed
field of characteristic $p>0$. An important class of modules for such
a group are the Weyl modules $\Delta(\lambda)$, labelled by dominant
weights; these can be constructed (relatively) explicitly, and their
heads provide a full set of simple modules for $G$. (Equivalently one
can study the duals of these modules, denoted $\nabla(\lambda)$ which
have the advantage of being induced from one-dimensional modules for a
Borel subgroup).  In determining the structure of such modules, or
indeed their cohomology, the calculation of Hom-spaces between them is
a useful tool.

Relatively little is known in general about such Hom-spaces. In type
$A$, when $\lambda$ and $\mu$ are related by a (suitable) single
reflection, explicit non-zero homomorphisms from $\Delta(\lambda)$ to
$\Delta(\mu)$ were constructed (with some restrictions) by Carter and
Lusztig \cite{carlus}, and (more generally) by Carter and Payne
\cite{cpmaps}. The corresponding cases in other types were considered
by Franklin \cite{franklin}. While it is clear that there should be a
hierarchy of families of homomorphisms corresponding to different
powers of $p$, the only case where the above results provide a
complete classification is when $G$ is $\SL_2$, where it is relatively
easy to determine all Hom-spaces exactly \cite{coxerd}.

The only other general results in this area, by Andersen
\cite{andweyl} and Koppinen \cite{kop1}, concern homomorphisms
between modules labelled by weights which are \lq close
together\rq. Typically such results show that certain Hom-spaces are
non-zero, or in some cases one-dimensional.  For weights which are far
apart and not related by a single reflection almost nothing is known.

In this paper we will determine all homomorphisms between Weyl modules
for $\SL_3$ when $p\geq 3$, and provide a recursive procedure for
determining the composition factors arising in the image (or kernel)
of such maps in most cases. From these results we will also classify all
homomorphisms between Specht modules for the symmetric groups
corresponding to three part partitions, when $p\geq 3$.

After a section of preliminaries, we review the $\SL_3$ data concerning
$p$-filtrations that we will need from \cite{par1}. This describes
certain filtrations of induced modules which will allow us to proceed
by induction, together with the set of $p$-good homomorphisms which
will be fundamental in our later constructions. We also recall a
theorem of Carter and Payne \cite{cpmaps} on the existence of certain
homomorphisms.  These will be the two key sets of data which we need
to determine all possible homomorphisms.

With the notation developed up to that point in place, in Section
\ref{strategy} we can give the strategy to be followed in the
remainder of the paper, and in particular the translation functor
arguments that form the basis of our argument. This section also
contains a more precise description of the main results that will be
obtained.

The rest of the paper takes the form of a single inductive argument on
the weight labelling our induced module. In Section \ref{onedimsec} we
show (Theorem \ref{onedim}) that $\Hom(\nabla(\lambda),\nabla(\mu))$
is at most one-dimensional. Unfortunately the proof relies on the fact
that for one very specific configuration of weights $\tau$ and $\nu$
there are no homomorphisms from $\nabla(\tau)$ to $\nabla(\nu)$
(Assumption \ref{assume}). Thus we cannot complete this proof until we
have classified all possible homomorphisms from induced modules
labelled by smaller weights, and hence it forms the first step in our
inductive argument.

We then come to the heart of the paper, Sections
\ref{cpcomp}--\ref{exception}. In Section \ref{cpcomp} we 
analyse the homomorphisms constructed by Carter and Payne in
the $\SL_3$ case. We begin by giving (up to a controllable
ambiguity in certain non-generic cases) a recursive procedure for
describing such homomorphisms (Proposition \ref{mapsare}), and use
this to determine which composites of such maps are non-zero
(Corollary \ref{composite}).

 We begin in Section \ref{homstart} to determine precisely which
Hom-spaces are non-zero, by starting with those obtained inductively
from composites of Carter-Payne maps. This uses translation functor
arguments, together with the explicit construction of homomorphisms
from the previous section.  For weights close to the boundary of the
dominant region it is more convenient to use a different argument
(Proposition \ref{symprop}), based on Doty's description of the
structure of symmetric powers \cite{doty} as reinterpreted in
\cite{cox4}.

As well as the Carter-Payne maps, there is another class of maps
obtained by \lq twisting', which are easy to describe. Unfortunately
not every map comes from Carter-Payne composites or twisting, and
hence in Section \ref{exception} we must construct the remaining maps
by hand (by gluing together maps defined on appropriate parts of the
$p$-filtrations of the two modules), to complete the
classification. With this we see that Assumption \ref{assume} does
indeed hold, and the main results in the paper now follow by
induction.

In Section \ref{examples} we illustrate our main results with some
examples for the case $p=3$.

Using a theorem of Carter and Lusztig
there is a close relationship between homomorphisms between induced
modules for 
$\SL_3(k)$ and between Specht modules for symmetric groups
$k\Sigma_d$.  In Section \ref{rdsect} we show how our
results also give a classification of homomorphisms between Specht
modules labelled by partitions with at most three non-zero parts when
$p\ge3$. 
Finally, we recall a result of Donkin relating Ext-groups between
induced modules for a reductive group and for its Levi factors, and
show how this can be used to derive a tensor product theorem of Fayers
and Lyle \cite{fayly} if $p$ is not two.


\section{Preliminaries}\label{prelim}
%
%
In this section we shall review the basic results required in this
paper.  Except where otherwise indicated, this material can be found
in found in \cite[II, Chapters 1--6]{jantzen}. Although we shall state
these results for an arbitrary reductive
algebraic group $G$, defined over an algebraically closed field $k$ of
characteristic $p > 0$, for most of this paper we will only be
interested in the case when $G = \SL_3(k)$.

We fix a maximal torus $T \subset G$ and hence a weight lattice
$X(T)$. There is an associated root system $R$, in which we choose a
set of positive roots $R^+$. The corresponding set of simple roots
will be denoted $S$. If $G$ is semisimple and simply--connected there
is a basis $\{\varpi_{\alpha}\mid \alpha \in S\}$ of fundamental weights
for $X(T)$.  Let $s_{\alpha}$ be the reflection on $E= X(T)
\otimes_{\zed} \RR$ given by $s_{\alpha} \lambda = \lambda-\langle
\lambda, \alphac \rangle \alpha$ where $\alphac$
is
the coroot associated to $\alpha $ in $X(T)^*$ and $\langle-,-
\rangle$ is the usual bilinear form.  The Weyl group $W$ is the group
generated by such reflections, while the affine Weyl group $W_p$ is
the semidirect product of $W$ with the group $p\zed R$ (acting by
translation on $E$). 

Setting $\rho = \frac{1}{2} \sum_{\alpha \in R^+} \alpha \in E$, we
define the dot action $w. \lambda = w( \lambda + \rho)-\rho $ of $W_p$
on $X(T)$.   Associated to this action is
a system of facets; these are all sets of the form
\begin{equation*}
\begin{split}
F= \{ \lambda \in E \mid \ 
&\langle \lambda+\rho, \alphac
\rangle = n_\alpha p\quad \forall\, \alpha \in R^+_0(F),
\\
&(n_\alpha -1)p < \langle \lambda +\rho, \alphac \rangle < n_\alpha p
\quad \forall\, \alpha \in R_1^+(F)\} 
\end{split}
\end{equation*}
for some integers $n_{\alpha}$ and a disjoint decomposition
$R^+=R_0^+(F) \cup R_1^+(F)$. A facet $F$ is called an alcove if
$R_0^+(F)= \emptyset$ and a wall if $|R_0^+(F)|=1$.  We
shall always assume that $p\geq h$ (the maximum of the Coxeter numbers
of the connected components of $G$), which
ensures that every alcove contains a weight.  The closure
(respectively upper/lower closure) of a facet $F$ 
is the set obtained from $F$ by
replacing both (respectively the right/lefthand) strict inequalities
occurring in the defining equations for $F$
 by non-strict
inequalities.  The closure of any alcove is a fundamental domain for
the (dot) action of the affine Weyl group, and we call the alcove containing
the origin the fundamental alcove. Occasionally we will need to
consider $p^e$-facets; these are defined as above but with all
occurrences of $p$ replaced by $p^e$ (and hence are associated to the
action of $W_{p^e}$).
 
We will also need to consider the set of dominant weights
$$X^+=\{ \lambda \in  X(T)\mid 0< \langle \lambda + \rho, \alphac
 \rangle\ \mbox{for all}\ \alpha \in S\}
$$
and its subset of $p$-restricted weights $X_1(T)$ defined by imposing
the additional constraint that $\langle \lambda + \rho, \alphac
\rangle \leq p$ for all $\alpha\in S$.
Any weight $\lambda$ can be uniquely written in the form $\lambda =
\lambda' + p \lambda''$ with $\lambda' \in X_1(T)$, and \emph{any
decomposition of a weight in this way will be assumed to be of this
form}. 


Given a Borel $T \subset B \subset G$ we can define the modules
$\nabla(\lambda) = \ind_ B^G k_\lambda$ where $k_{\lambda}$ is the
one-dimensional $B$-module of weight $\lambda$.  By choosing $B$
appropriately, we may arrange that $\nabla(\lambda)$ is non--zero
precisely when $\lambda$ is dominant.  There is a contravariant
duality on $G$, and the Weyl module
$\Delta(\lambda)$ is the contravariant dual of $\nabla (\lambda)
$. Both $\nabla (\lambda)$ and $\Delta(\lambda)$
has the same character, given by Weyl's character formula. Clearly we
have $\Hom(\Delta(\mu), \Delta(\lambda))\cong\Hom(\nabla(\lambda),
\nabla(\mu))$, and hence for the remainder of this paper we will
consider only induced rather than Weyl modules.  (Note that throughout
we will abuse notation and write $\Hom$ for $\Hom_G$.)

For each $\lambda \in X^+$, the module $\nabla(\lambda)$ has simple
socle $L(\lambda)$, and all simple modules arise in this manner.  All
other composition factors $L(\mu)$ of $\nabla(\lambda)$ satisfy $\mu <
\lambda$ in the dominance order determined by $R^+$, and the strong
linkage principle implies that $\mu \in W_p. \lambda $.  Thus
$\Hom(\nabla(\lambda), \nabla(\mu))$ is non--zero only if $\mu \leq
\lambda$ and $\mu \in W_p. \lambda$.  Note that for a given weight
$\lambda$, the highest weight of any composition factor of
$\nabla(\lambda)$ can be specified by determining the facet in which
this highest weight lies.  By Steinberg's tensor product Theorem, for
any dominant weight $\lambda$ we have an isomorphism $L(\lambda) \cong
L(\lambda') \otimes L(\lambda'')\frob $, where $F$ is the Frobenius
morphism.  In the special case $\lambda = (p-1)\rho $ we have that
$\nabla(\lambda) \cong L(\lambda)$, and we call this the Steinberg
module $\St$. All simple modules are contravariantly self-dual.

Given weights $\lambda, \mu$ in the closure of some alcove $F$,
there is a unique dominant weight $\nu$ in $W(\mu-\lambda)$.  We
define the translation functor $T_{\lambda}^{\mu}$ from $\lambda$ to
$\mu$ on a module $V$ by 
$T_{\lambda}^{\mu} V= \pr_{\mu} (L(\nu) 
\otimes \pr_{\lambda} V)$, where $\pr_{\tau} V$ 
is the largest submodule of $V$ all of whose composition
factors have highest weights in $W_p.\tau $. We summarise the
properties of translation functors that we shall require in the
following proposition (proofs can be found in 
\cite[II, Chapter 7]{jantzen}).

\begin{proposition} \label{tiltfacts} Let $\lambda$ and $\mu$ be in the
closure of some facet $F$. Then
\begin{enumerate}
\item[(i)]{$T_{\lambda}^{\mu}$ is an exact functor, and is adjoint to
$T_{\mu}^{\lambda}$.}  
\item[(ii)]{If $\lambda $ and $\mu $ lie in the same facet then 
$T_{\lambda}^{\mu} $ is an equivalence of categories from the category of
modules $V$ with $\pr_{\lambda} V=V$ to the category of modules $V$
with $\pr_{\mu} V=V$.} 
\item[(iii)]{If $\mu$ is in the closure of the
facet containing $\lambda$ then $T_{\lambda}^{\mu} \nabla(\lambda)
\cong \nabla(\mu)$.} 
\item[(iv)]{More generally, $T_{ \lambda}^{\mu}
\nabla (\lambda)$ has a filtration $0=T_0 \subset \cdots \subset T_{t-1}
\subset T_t= T_{\lambda}^{\mu} \nabla(\lambda)$ such that the
factors $T_i/T_{i-1}$ are all of the form $\nabla(w_i. \mu)$ where
$w_i \in \Stab_{W_p} (\lambda)$ with $w_i.\mu \in X^+$, each
such factor occurring exactly once.  Further, this filtration may be
taken such that $w_i. \mu < w_j. \mu $ implies that $i < j$.}
\end{enumerate}
\end{proposition} 

We say that a module $V$ has a good filtration if it has a filtration
with all factors of the form $\nabla(\lambda_i) $ with $\lambda_i \in
X^+$, and denote by $(V: \nabla(\mu))$ the number of factors of the
form $\nabla(\mu)$ in such a filtration.  (This is independent of
our choice of good filtration.) We call any module $V$ such that both
$V$ and its dual have a good filtration a tilting module.  By
\cite{dontilt} there is for each dominant weight $\lambda$ a unique
indecomposable tilting module $T(\lambda)$ with that highest weight,
and every indecomposable tilting module arises in this way.

Similarly we say that a module $V$ has a good $p$-filtration if it has
a filtration all of whose factors of the form
$\nabla_p(\lambda)=\nabla (\lambda_i'' )\frob \otimes L (\lambda_i')$
(recall our standing convention for such a decomposition of
$\lambda$).  In \cite[II, Proposition 9.11]{jantzen} Jantzen gives a
criterion for $\nabla(\lambda)$ to have a good $p$-filtration, while
for $p\geq 2h-2$ Andersen \cite[3.7 Corollary]{andpfilt} has shown
that such good $p$-filtrations always exist. For our purposes it is
enough to note \cite[3.13]{jancrelle} that $\nabla(\lambda)$ with
$\lambda \in X^+$, always has a
good $p$-filtration in the case of $\SL_3(k)$. As these are the only
type of $p$-filtrations which we will consider, we will usually omit
the qualifier \lq good\rq.

For any dominant weight $\lambda$ we have that 
$\nabla(p\lambda+(p-1)\rho) \cong \nabla(\lambda)\frob\otimes\St$. 
Indeed the functor $V \longmapsto V\frob \otimes \St$ 
induces \cite[II, 10.5(1)]{jantzen} an
equivalence of categories from the category of all modules to the
category of those modules all of whose composition factors are of the
form $L(\lambda)$ with $\lambda' = (p-1)\rho$. More generally, the
functor $V \longmapsto V\frob \otimes L(\lambda)$ where $\lambda\in
X_1(T)$ induces an isomorphism of submodule lattices \cite[Lemma
1.2]{erdhenk2}.  

We denote the (scheme theoretic) kernel of the Frobenius morphism by
$G_1$.  The simple modules for this subgroup are precisely the
restrictions of the simple $G$-modules with $p$-restricted highest
weights.  As $\Hom_G(U,V) \cong (V\otimes U^*)^G \cong ((V \otimes
U^* )^{G_1} )^{G/G_{1}}$, and $G \cong G/G_1 $ via $F$, we deduce
that
\begin{equation}\label{homreduce}
\Hom(\nabla_p(\lambda),\nabla_p(\mu))\cong\left\{
\begin{array}{ll}
\Hom(\nabla(\lambda''),\nabla(\mu''))& \wif \lambda'=\mu'\\ 
0 & \otherwise\end{array}\right.
\end{equation}
for all dominant weights $\lambda$ and $\mu$.

Any $G$-module $W$ which is trivial as a $G_1$-module is of the form
$V\frob$ for some $G$-module $V$, and we define $W^{(-1)}=V$. The (usual)
dual of $\nabla(\lambda)$ is a Weyl module, which we denote by
$\Delta(\lambda^*)$.  Finally, if
$U, V, X, Y$ are finite dimensional $G$-modules then we have
\begin{equation}\label{infext}
\Ext_{G_1}^i(U\otimes X\frob, V\otimes Y\frob)\cong \Ext_{G_1}^i(U, V)\otimes
(X^*)\frob\otimes Y\frob
\end{equation}
(as a $G$-module) for all $i\geq 0$.




\section{$\SL_3$ data}
\label{localdata}
%
%
In this section we will describe explicitly the good $p$-filtration
for $\nabla(\lambda)$ in the case of $\SL_3(k)$. The detailed
structure of these filtrations was given in \cite[Theorem 4.12]{par1},
and we summarise those results here.  Further information about these
filtrations also appears in \cite[Chapter 2]{kh}. We will use this to
give an explicit description of \emph{$p$-good homomorphisms} (defined
below) 
between induced
modules, and recall a general result about homomorphisms due to Carter
and Payne which will be needed in what follows.

\begin{definition}
A \emph{$p$-good homomorphism} is a homomorphism $\phi$ between two modules
which both have $p$-good filtrations, such that the 
image and kernel of $\phi$ both have good $p$-filtrations.
\end{definition}

For $\SL_3(k)$ we shall denote the two simple roots by $\alpha_1$ and
$\alpha_2$. A weight
$\lambda \in X(T)$ can be identified with a pair of integers $(a,b)$
via $\lambda = a \varpi_{\alpha_1} +b \varpi_{\alpha_2}$.  With this
convention we can identify $X^+ $ with the subset $\NN^2 $ of $\zed^2
$. The set of $p$-restricted weights $X_1(T)$ contains two alcoves;
the fundamental alcove $C_0$ is the alcove containing the origin.
Any alcove which is a translate of the fundamental alcove will
be called a {\it down} alcove, while translates of the remaining
$p$-restricted alcove will be called {\it up} alcoves.  We extend our
earlier notation and write $\lambda = p(a'',b'')+(a',b')$ where
$(a',b') \in X_1(T)$.  In order to describe the good $p$-filtration of
$\nabla(\lambda)$, it is enough (by the strong linkage principle) to
indicate the facets in which the highest weights of the factors in the
$p$-filtration occur.  By the translation principle (Proposition
\ref{tiltfacts}(ii)) the configuration of facets obtained depends only
on the facet containing $\lambda$ and not on the weight itself.  Thus
we will represent the factors $\nabla_p(\lambda_i) $ by numbers
corresponding to an appropriate labelling of the facets in each case.

In the diagrams representing good $p$-filtrations, the factors of a
filtration will be represented by labels with lines connecting them.
Two factors will be joined if and only if there is a nontrivial
extension (in $\nabla(\lambda)$) between them.  The diagrams are
oriented so that sections that embed in $\nabla(\lambda) $ occur at
the bottom of the diagram.  If $\lambda_i''$ is not dominant for some
factor $\nabla_p(\lambda_i)$ occurring in one of our diagrams then we
interpret this factor as the zero module (i.e. we ignore that part of
the diagram).

\begin{figure}[ht]
\centerline{\epsffile{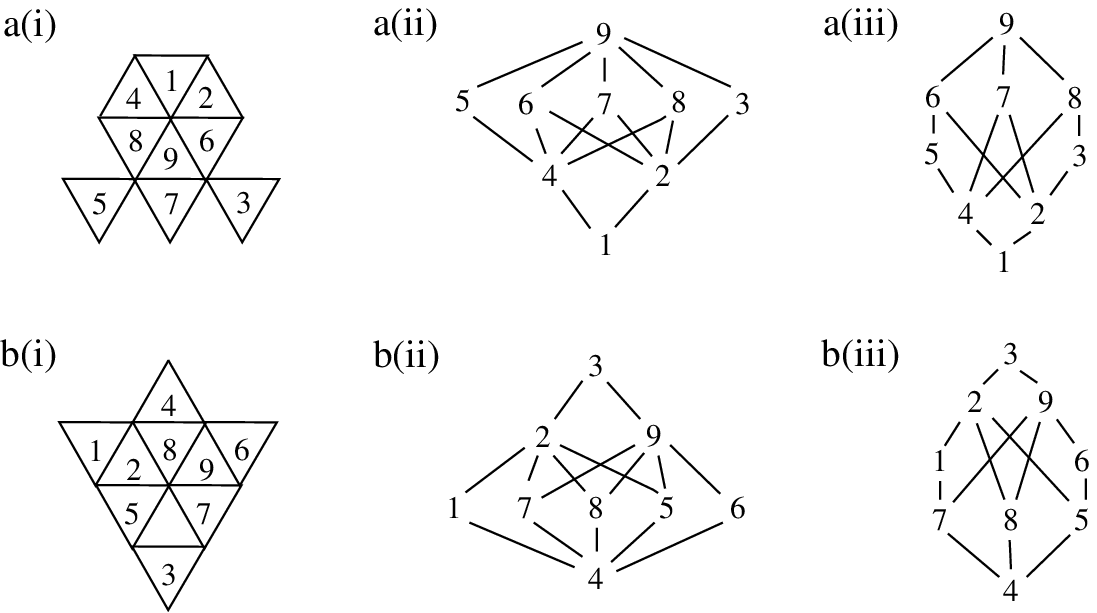}}
\caption{\label{alcoves}}
\end{figure}

We now list the various cases that can arise, according to
\cite[Theorem 4.12]{par1}.  If $a''$ and $b'' \geq 1$ 
and $(a',b')$ lies in the fundamental alcove then there are nine
factors in the $p$-filtration as indicated in Figure
\ref{alcoves}a(i), and the structure of this filtration (when $a''
$ and $b'' \notequiv 0 \pmod p$) is given in Figure
\ref{alcoves}a(ii).  If $a'' \equiv 0 \pmod p$ (respectively 
$b'' \equiv 0 \pmod p$) then there is an additional extension of
$\nabla_p(\lambda_5)$ by $\nabla_p(\lambda_6)$ (respectively of
$\nabla_p(\lambda_3)$ by $\nabla_p(\lambda_8)$).  If both $a'' $ and
$b'' \equiv 0 \pmod p$ then the structure is given in Figure
\ref{alcoves}a(iii); the other cases are similar.

If $(a',b')$ lies in the interior of the other $p$-restricted alcove
then there are also nine factors in the $p$-filtration, as indicated
in Figure \ref{alcoves}b(i), and the structure of the filtration again
depends on $a''$ and $b''$. If both $a''$ and $b'' \notequiv -1 \pmod
p$ then it is given in Figure \ref{alcoves}b(ii) while if $a'' \equiv
-1 \pmod p$ (respectively if $b'' \equiv -1 \pmod p$) then there
is an additional extension of $\nabla_p(\lambda_5)$ by 
$\nabla_p(\lambda_6)$ 
(respectively of $\nabla_p(\lambda_7)$ by $\nabla_p(\lambda_1)$). 
If both $a''$ and $b'' \equiv -1 \pmod p$ then the
structure is illustrated in Figure \ref{alcoves}b(iii); the other
cases are similar.

\begin{figure}[ht]
\centerline{\epsffile{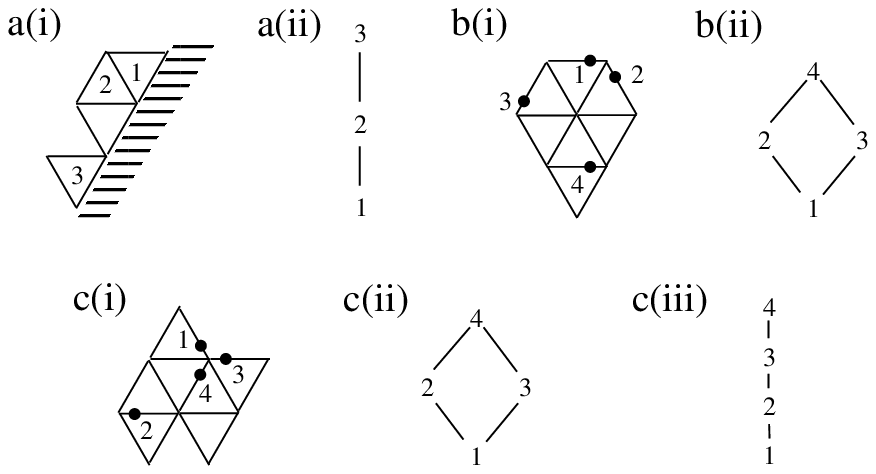}}
\caption{\label{walls}}
\end{figure}

There are two remaining alcove cases: where $\lambda$ lies in a down
alcove and either $a'' =0$ or $b''=0$. As these cases are symmetric we
only consider the latter.  So suppose that $\lambda =
p(a'',0)+(a',b')$ with $(a',b') \in C_0$.  Then $\nabla(\lambda)$ has
three factors in the $p$-filtration, as indicated in Figure
\ref{walls}a(i) (where the shaded region represents the boundary of
the dominant region), and the structure of the filtration is given in
Figure \ref{walls}a(ii). When we need to distinguish between the two
down alcove cases we shall refer to those as in Figure
\ref{alcoves}a(i) as {\it internal} down alcoves, and those as in
Figure \ref{walls}a(i) (and the symmetric version thereof) as {\it
just dominant} down alcoves. We shall also refer to the up
alcoves immediately above just dominant down alcoves as just dominant,
and similarly the walls between two just dominant alcoves. All other
walls will be referred to as interior walls.
To distinguish further the just dominant down alcove cases we
shall refer to that illustrated in Figure \ref{walls}a(i) as the
{\it right-hand} case, and to the symmetric version (which we number
in the symmetric fashion) as the {\it
left-hand} case.  

Next we consider the various wall cases.  If $a'+b'=p-2$ then
$\lambda$ lies on one of the horizontal walls in our diagrams.  There
are now four factors in our filtration, as indicated in Figure
\ref{walls}b(i), and the structure of the filtration is given in
Figure \ref{walls}b(ii). The remaining two wall cases are symmetric,
so we only consider the case where $\lambda = p(a'',b'')+(p-1,b')$.
We refer to this case as the {\it left-hand diagonal wall.} Once again
there are four factors in the filtration as indicated in Figure
\ref{walls}c(i). If $a'' \notequiv -1 \pmod p$ then the structure is
as in Figure \ref{walls}c(ii), otherwise it is as in Figure
\ref{walls}c(iii). For right-hand walls we will use the numbering
of the four factors given by  reflecting the diagram (so that the top and
bottom factors in the filtration are also labelled by $4$ and $1$
respectively). The only remaining case is when $\lambda$ lies on a
vertex, i.e. $\lambda' = (p-1)\rho$.  But in this case there is only a
single factor in the $p$-filtration, as noted in the previous section.



In what follows we will frequently make use of 
\begin{lemma} Suppose that $G=\SL_3(k)$.  For all $\lambda \in X^+$
the module $\nabla_p(\lambda)$ has a simple head.
\end{lemma} 
\pf Jantzen has shown \cite[6.9]{jancrelle} that $\nabla(\mu)$ has simple
head for all dominant weights $\mu $, provided that $p > 3$, and in
\cite[Proposition 4.11]{par1} the restriction on $p$ has been
removed. Now the result follows from \cite[Lemma 1.2]{erdhenk2}.
\qed

Our goal is to determine $\Hom(\nabla(\lambda), \nabla(\mu))$ for
all pairs of dominant weights $\lambda$ and $\mu$.  In order to
construct homomorphisms we will need the following result of Carter
and Payne.

\begin{theorem} \label{morphism}
Suppose that $G=\SL_n(k)$ and that $\lambda, \mu \in X^+$
satisfy the following conditions:    
\begin{enumerate}
\item[(i)]{$\mu < \lambda $.  
}
\item[(ii)]{There exists some integer $e>0$ such that:} 
\begin{enumerate}
\item[(a)]{$\lambda$ and $\mu$ are mirror images in some $p^e$-wall
$L$ and}
\item[(b)]{$L$ is the unique $p^e$-wall between $\lambda $ and $\mu$ (possibly containing
$\lambda$ or $\mu$) parallel to $L$.}
\end{enumerate}
\end{enumerate}
Then $\Hom(\nabla(\lambda), \nabla(\mu)) \neq 0$.
\end{theorem} 
\pf See \cite{cpmaps}.  \qed

When we refer to such a homomorphism from $\nabla(\lambda)$ to
$\nabla(\mu)$ as reflection about a $p^e$-wall, we shall assume that
$e$ is chosen to be minimal with that property. For $\SL_3$ we shall
refer to the three cases as reflection to the left (L), right (R), and
below (B), depending on the relative positions of $\lambda$ and $\mu$.
  
Theorem \ref{morphism} allows us to construct homomorphisms between widely
separated weights, but provides no information as to their structure.
For $\SL_3(k)$ we will
be able to provide this extra information using the explicit
description of $p$-good  morphisms 
given in \cite{par1}.
Thus we conclude this section by reviewing this data.

In the following lemma we describe all non-zero $p$-good homomorphisms
between induced modules. 
We will consider each of the cases illustrated in
Figures \ref{alcoves} and \ref{walls} in turn (the remaining cases are
symmetric).  In each case we will describe homomorphisms from
$\nabla(\lambda)$ \emph{in terms of the labelling of facets given in
the figure corresponding to $\lambda$.}  The image of such a
homomorphism $\nabla(\lambda) \to \nabla(\lambda_i)$ will then be
isomorphic to the smallest quotient of the $p$-filtration of
$\nabla(\lambda)$ containing $\nabla_p(\lambda_i)$, which can be read off the
corresponding diagrams in Figures \ref{alcoves} and \ref{walls}.  

\begin{lemma}\label{lem:pgood} For a dominant weight $\lambda$ there
exists a non-zero $p$-good homomorphism
$\nabla(\lambda)\to\nabla(\lambda_i)$ if and only if $\lambda_i$ is
dominant and the pair $(\lambda,i)$ occurs in the following
table (where we identify $\lambda$ with the type of facet in which it
occurs):

\centerline{\begin{tabular}[pos]{|p{2.7in}|p{2.15in}|}
\hline
 \begin{minipage}{2.7in}\centerline{$\lambda$}\end{minipage}
&\begin{minipage}{2.15in}\centerline{$i$} \end{minipage}\\
\hline
  Fig \ref{alcoves}a: internal down alcove& 1, 2, 4, 6, 8 and 9\\
\hline
  Fig \ref{alcoves}b: up alcove& 1, 2, 3, 4, 6, 8 and 9\\
\hline
  Fig \ref{walls}a: just dominant down alcove& 1, 2 and 3\\
\hline
  Fig \ref{walls}b: horizontal wall& 1, 2, 3 and 4\\
\hline
  Fig \ref{walls}c: interior diagonal wall&1, 3 and 4\\
\hline
  Fig \ref{walls}c: just dominant diagonal wall&1 and 2\\
\hline
\end{tabular}}
\end{lemma}

\pf Argue as in the proof of \cite[Lemma 5.1]{par1}, or \cite[II, 9.14
Remark 4]{jantzen}.
\qed

Note that when all terms in the $p$-filtration of $\nabla(\lambda)$
are simple, the above Lemma gives all weights $\mu$ such that
$\Hom(\nabla(\lambda),\nabla(\mu))\neq 0$. 
All but four of the $p$-good
homomorphisms above can be constructed as composites of Carter-Payne
maps. The four exceptions are the map from $\nabla(\lambda)$ to
$\nabla(\lambda_3)$ in the interior up alcove case, 
the map from $\nabla(\lambda)$ to
$\nabla(\lambda_4)$ in the interior horizontal wall case, 
the map from $\nabla(\lambda)$ to
$\nabla(\lambda_2)$ in the just dominant diagonal wall case, 
and the
map from $\nabla(\lambda)$ to $\nabla(\lambda_3)$ in the just dominant
down alcove case. 
We will call these latter two maps 
the {\it exceptional} $p$-good maps.
(It can happen that these exceptional maps can be constructed as the
composite of Carter-Payne maps, as in second example in
section \ref{examples}. This is non-generic behaviour. Generically, the
exceptional $p$-good maps will not be the composite of Carter-Payne
maps.)
Note that all of these $p$-good maps which are generically composites
of Carter-Payne maps are the maps with image the $G_1$-head of
$\nabla(\lambda)$.




\section{Overall strategy}\label{strategy}

As the remainder of this paper involves some rather intricate
calculations, we devote this section to an overview of the strategy
behind the proof. From this point on we will assume that $G=\SL_3(k)$
and that $p>2$ (so that all our facets are non-empty).

Given a dominant weight $\lambda$, we assume that all homomorphisms
from $\nabla(\tau)$ with $\tau<\lambda$ have been classified, and that
all the corresponding Hom-spaces are one-dimensional.  First we claim
that $\dim \Hom(\nabla(\lambda), \nabla(\mu))$ is at most one
dimensional for all weights $\mu \in X^+$.  Second we claim that 
we can find all weights $\mu < \lambda$ such that
$\Hom(\nabla(\lambda), \nabla(\mu))$ is non-zero.  By induction this
will calculate all possible Hom-spaces for $\SL_3$ with $p>2$.

The basic idea is to use translation functors. For weights on a
vertex, the result is known by induction using (\ref{homreduce}), or
the equivalence of categories discussed just before that
equation. Next suppose that $\lambda$ lies in the closure of a down
alcove whose lowest vertex $\nu$ is dominant. By induction the
Hom-spaces for $\nabla(\nu)$ are known.  Let
$$\Theta_{\nu} = \{ \theta \in X^+ \mid \Hom(\nabla(\nu),
\nabla(\theta)) \ne 0 \}.$$ 
We `translate' all the weights in $\Theta_{\nu}$
to a new set $\Gamma_{\lambda}$ with
\begin{equation}\label{downtrans}
\Gamma_{\lambda} = \{ \gamma \in X^+ \mid ( T_\nu^\lambda
\nabla(\theta):\nabla(\gamma)) \ne 0 \mbox{ for some }\theta \in
\Theta_{\nu} \}.
\end{equation}

Now suppose $\lambda$ lies in an up alcove, and choose $\tau$ on the horizontal
wall below it. If $\tau$ lies in the closure of a down alcove with
lower vertex $\nu$, and $\nu$ is dominant, then we further translate
the weights in $\Gamma_{\tau}$ to form
\begin{equation}\label{uptrans}
\Gamma_{\lambda} = \{ \gamma \in X^+ \mid ( T_{\tau}^\lambda
\nabla(\theta):\nabla(\gamma)) \ne 0 \mbox{ for some }\theta \in
\Gamma_{\tau} \}.
\end{equation}

In both cases (\ref{downtrans}) and (\ref{uptrans}) we now remove all
weights in $\Gamma_{\lambda}$ which are not composition factors of
$\nabla(\lambda)$.  We would like to claim that the resulting set
$$\Upsilon_{\lambda}= \{ \gamma \in \Gamma_{\lambda} \mid [\nabla(\lambda):
L(\gamma)] \ne 0 \}$$ 
is precisely the set of weights $\eta$ for which
$\Hom(\nabla(\lambda), \nabla(\eta))\cong k$, and that if 
$\eta \not
\in \Upsilon_{\lambda}$ then $\Hom(\nabla(\lambda), \nabla(\eta)) \cong
0$. However, this is not quite the case and will need to be modified.

In the course of showing that all Hom-spaces are at most
one-dimensional (Theorem \ref{onedim}) we easily show that
$\Gamma_{\lambda}$ is an upper bound for the set of weights for which
non-zero homomorphisms can exist, and it is obvious that this can then
be refined to give $\Upsilon_{\lambda}$ as an upper bound.

We already have one family of homomorphisms: those constructed by
Carter and Payne. There is another obvious class of maps, which we
will now describe. Suppose that $\nabla_p(\tau)$ is the top term in a
$p$-filtration of some $\nabla(\lambda)$. Given any map from
$\nabla(\tau'')$ to $\nabla(\mu)$ (which we will know by induction) we
obtain a map from $\nabla_p(\tau)$ to $\nabla_p(\tau'+p\mu)$, and
hence from $\nabla(\lambda)$ to $\nabla(\tau'+p\mu)$ (by killing all
other terms in the $p$-filtration of $\nabla(\lambda)$). We will refer
to maps arising in this way as {\it twisted maps}.  (Of course maps
can be both twisted maps and a composite of Carter-Payne maps.)

We would like to claim that all homomorphisms are either composites of
Carter-Payne maps, or twisted maps. However, this too is not quite
correct. To see why, we must consider those dominant $\lambda$ not
covered by the cases above. First suppose that $\lambda$ lies in the
closure of a just dominant down alcove. Although we cannot apply the
translation functor approach outlined above, it is in fact possible to
determine all homomorphisms directly (Theorem \ref{symprop}), using
the explicit description of the structure of symmetric powers given by
Doty \cite{doty}. (For $\lambda$ in an up alcove above a just dominant
down alcove, the determination of homomorphisms is now straightforward
by translation arguments.)

We now return to the sets $\Upsilon_{\lambda}$. If $\nu$ and $\theta$
(as used in the definition of these sets) are related by a composite
of Carter-Payne maps then for each weight in $\Upsilon_{\lambda}$ we
can either (i) explicitly construct all possible homomorphisms using
composites of Carter-Payne maps or twisted maps, or (ii) show no map
exists. The latter case is rare: for \lq generic' weights case (ii)
never occurs for $\lambda$ in a down alcove, and only in very special
configurations for $\lambda$ in an up alcove. For the precise
statement for down alcoves see Theorem
\ref{downhomclass}, and for up alcoves Theorem \ref{uphomclass}.

If $\nu$ and $\theta$ are instead related by a twisted map, then in
most cases the analogues of Theorems \ref{downhomclass} and
\ref{uphomclass} are straightforward, as every weight in
$\Upsilon_{\lambda}$ can be obtained by a twisted map. However, when
(for example) $\nu=(p-1)\rho+p\lambda$ with $\lambda$ close to the
boundary of the dominant region, there are additional weights in
$\Upsilon_{\lambda}$ which can be reached neither by composites of
Carter-Payne maps nor twisted maps. Unfortunately in these cases maps
do exist, and Section \ref{exception} of the paper is devoted to the
construction of these {\it exceptional maps}.

We will show that there are a pair of terms at the top of the
$p$-filtration for $\nabla(\lambda)$ which individually map to a pair
of terms at the bottom of the $p$-filtration for $\nabla(\mu)$. By
considering pullbacks and pushouts we will show that this pair of maps
can be \lq glued together\rq\ to give the required exceptional
map. This will complete our classification of homomorphisms.

\section{Dimensions of Hom-spaces}\label{onedimsec}
%
%
In this section we begin an inductive procedure that will continue for
the rest of the paper. We show that all Hom-spaces between induced
modules are at most one-dimensional, provided that there are no
homomorphisms between certain special pairs of induced modules
labelled by smaller weights. The verification of this hypothesis will
follow from the remainder of the paper.

We say that two weights $\lambda$, $\mu$ are in the {\it same
$\nabla_p$-class} if $\lambda' = \mu'$ and $\lambda''$ and $\mu''$ are
in the same $G$-block. Any weight lies in a unique translate of the
set of $p$-restricted weights; note that $\lambda$ and $\mu$ can only
be in the same $\nabla_p$-class if they both lie in the same position
in the respective translates (and also have the Steinberg weights
immediately above each translate lying in the same $W_{p^2}$-orbit).
Clearly if $\lambda$ and $\mu$ are in distinct $\nabla_p$-classes,
then $\nabla_p(\lambda)$ and $\nabla_p(\mu)$ have no common
composition factors. We will also say that two modules are in the same
$\nabla_p$-class if the highest weights of their composition factors
are all in the same $\nabla_p$-class.

\begin{theorem}\label{onedim}
For all $\lambda, \mu \in X^+ $ we have $$\dim \Hom (\nabla (\lambda),
\nabla (\mu)) \leq 1.  $$ Further, for those $\lambda$ for which the
set $\Upsilon_{\lambda}$ from Section \ref{strategy} has been defined
we have  $$\Hom (\nabla (\lambda),\nabla (\mu))=0$$ for all
$\mu\notin\Upsilon_{\lambda}$.
\end{theorem} 
\pf We proceed by induction on $\lambda $, and may assume that $\mu\in
W.\lambda$ with $\mu\leq \lambda$.  If $\lambda $ is in the
fundamental alcove then $\nabla (\lambda) $ is simple and so we are
done by \cite[II 2.8 Proposition]{jantzen}.  Now suppose that $\lambda
$ is in the closure of a just dominant down alcove (i.e there
\emph{does not} exist a dominant weight $\theta $ on the vertex
immediately below this).  By symmetry, it is enough to consider the
case where $\lambda $ is near to the right hand boundary of the
dominant region.

We first consider the two wall cases.  By Figure \ref{walls}b or c we
have a short exact sequence
\begin{equation}\label{usepfilt}
0 \to \nabla_p (\lambda_1) \to \nabla (\lambda) \to \nabla_p
(\lambda_a) \to 0 
\end{equation}
where $a=3$ in case (b) and $a=2$ in case (c). (Recall that two of the
four possible terms in the $p$-filtration are zero when $\lambda$ is
just dominant.)  Hence there exists a long exact sequence
$$0 \to \Hom (\nabla_p (\lambda_a), \nabla (\mu)) \to \Hom (\nabla
(\lambda), \nabla (\mu)) \to \Hom (\nabla_p (\lambda_1), \nabla (\mu))
\to \cdots. $$ It will be enough to show that the first and third
Hom-spaces $\Hom (\nabla_p (\lambda_i), \nabla (\mu)) $ are each at
most one dimensional, as they can never both be non-zero as
$\lambda_1$ and $\lambda_a$ are in different $\nabla_p$-classes.

We begin by showing that the image in either case of any non-zero map
$\phi$ from $\nabla_p(\lambda_i)$ to $\nabla(\mu)$ must be a submodule
of $\nabla_p (\mu) $.  Clearly $\nabla_p (\lambda_i) $ and $\nabla_p
(\mu) $ must be in the same $\nabla_p $-class, as the socle of $\nabla
(\mu) $ must be a common composition factor.  By Figure \ref{walls}b
or c, the only other term in the $p $-filtration of $\nabla_p (\mu) $
in this $\nabla_p $-class occurs when $\mu $ is as in case (b), in
which case it is the factor $\nabla_p (\mu_4) $. (Recall that to be in
the same $\nabla_p$-class two weights must lie in the same position in
their respective translates of the set of $p$-restricted weights.)  As
all composition factors of $\im \phi $ lie in the same $\nabla_p
$-class, it is enough to show that none of the composition factors in
$\nabla_p (\mu_4) $ occur in this image.  However, the extension of
$\nabla_p (\mu_2) $ by $\nabla_p ( \mu_4) $ has a simple socle, as
this extension is the image under the map $\nabla (\mu_1) \to \nabla
(\mu_2) $.  Therefore if any composition factor of $\nabla_p (\mu_4) $
does occur in $\im \phi $ then so must $\soc \nabla_p (\mu_2) $, which
is impossible.

Thus we have shown that $\im \phi \leq \nabla_p (\mu) $.  Now
(\ref{homreduce}), and the induction hypothesis, immediately implies
that the desired Hom-spaces are at most one dimensional, and cannot
both be non-zero.

Next suppose that $\lambda$ is in the interior of a just dominant down
alcove, (and hence that there exists a weight $\theta$ on the diagonal
wall below $\lambda $).  We have a short exact sequence
\begin{equation}\label{offwall}
0 \to \nabla (\nu) \to T_{ \theta}^{ \lambda} \nabla (\theta)
\to \nabla (\lambda) \to 0 
\end{equation}
for some $\nu $. This gives rise to
an exact sequence 
\begin{equation}\label{tiltstep}
0 \to \Hom (\nabla (\lambda), \nabla (\mu)) \to \Hom (
T_{\theta}^{\lambda} \nabla (\theta), \nabla (\mu))
\end{equation}
and it is enough to show that this final Hom-space is at most
one-dimensional.  However
$$\Hom (T_{\theta}^{\lambda} \nabla (\theta), \nabla (\mu)) \cong \Hom
(\nabla (\theta), T_{\lambda}^{\theta} \nabla (\mu)) \cong \Hom
(\nabla (\theta), \nabla (\tau)) $$ for some $\tau $ (by Proposition
\ref{tiltfacts}(iii)), and we are done by induction.

We next consider the case where $\lambda$ is in the lower closure of
an internal 
down alcove (so there exists a dominant weight $\theta$ on the vertex
immediately below this). If $\lambda = \theta $ then $\nabla (\lambda)
\cong \nabla (\lambda'')\frob\otimes \St $ and block considerations show
that to have a non-zero homomorphism we must have $\nabla (\mu) \cong
\nabla (\mu'')\frob \otimes \St $.  Now by (\ref{homreduce}) we have that
$$\Hom (\nabla (\lambda'')\frob \otimes \St, \nabla (\mu'')\frob \otimes
\St) \cong \Hom (\nabla (\lambda''), \nabla (\mu''))$$ and we are done
by induction.

Now suppose that $\lambda$ is in the interior of an internal down alcove.  By
Proposition \ref{tiltfacts} we have that $T_{\theta}^{\lambda} \nabla
(\theta)$ has six factors of the form $\nabla(\alpha)$, and there
is a short exact sequence
$$0 \to X \to T_{\theta}^{\lambda} \nabla(\theta) \to \nabla(\lambda) 
\to 0$$ 
where X has a good filtration.  Hence we have an
exact sequence as in (\ref{tiltstep}), and the argument follows
exactly as above.  For later use we note that (in this case) if the
final Hom-space in (\ref{tiltstep}) is non-zero then $\mu $ must have
been in one of the six alcoves adjacent to $\tau $.  If $\lambda $
lies on one of the two walls then the argument is similar, but there
are only three weights $\mu $ which could give rise to the weight
$\tau $.

When $\lambda $ is on the wall in the upper closure of an internal
down alcove the proof is a little more complicated, so we shall
first show how the up alcove case follows from it. Thus we suppose
that $\lambda $ lies in the interior of an up alcove and let $\theta $
be a weight on the wall immediately below $\lambda $.  We again have a
short exact sequence as in (\ref{offwall}), and the argument proceeds
as in that case.

Finally we consider the case when $\lambda $ is on the wall in the
upper closure of an internal down alcove.  We have a short exact
sequence
$$0 \to M_1 \to \nabla (\lambda)\to M_2 \to 0$$
where $M_1$ and $M_2$ are defined by the sequences
$$0 \to \nabla_p (\lambda_1) \to M_1 \to \nabla_p (\lambda_2) \to
0\quad\wand\quad 0 \to \nabla_p (\lambda_3) \to M_2 \to \nabla_p
(\lambda_4) \to 0.$$
Arguing as for the sequence (\ref{usepfilt}), we see that 
$\dim\Hom(M_i,\nabla(\mu))\leq 1$ for $i=1,2$, and thus from the
sequence 
$$0 \to \Hom (M_2, \nabla (\mu)) \to \Hom (\nabla (\lambda), \nabla
(\mu))\to \Hom (M_1, \nabla (\mu)) $$
we see that we are done unless both these Hom-spaces are non-zero.

Thus we may assume that $\Hom(M_i,\nabla(\mu))\neq 0$ for
$i=1,2$. When $i=1$ this implies that $\mu'=\lambda_1'$ or
$\lambda_2'$, and for $i=2$ that $\mu'=\lambda_3'$ or
$\lambda_4'$. Therefore we must have $\mu'=\lambda_1'=\lambda_4'$, and
so $\mu$ lies on a horizontal wall. Let $\theta$ be a weight lying in
the interior of the alcove immediately below $\lambda$. We have that 
$$\Hom(\nabla (\lambda), \nabla (\mu))\cong
\Hom (T_{\theta}^{\lambda} \nabla (\theta), \nabla (\mu)) \cong \Hom
(\nabla (\theta), T_{\lambda}^{\theta} \nabla (\mu))$$
and an exact sequence
$$0 \to \Hom (\nabla (\theta), \nabla (\mu_l)) \to \Hom ( \nabla
(\theta), T_{\lambda}^{\theta}\nabla (\mu)) \to \Hom (\nabla
(\theta), \nabla (\mu_u))$$ where $\mu_l$ (respectively $\mu_u$) in
$W.\theta$ lies in the lower (respectively upper) alcove adjacent to
$\mu$. By induction the outer Hom-spaces in this sequence are at most
one-dimensional, so it is enough to prove that they cannot both be
non-zero.

As $\theta'=\mu_l'$ (as $\mu$ is in the same block of $\lambda$
and they are both on horizontal walls), we will be done if we can show

\begin{ass} \label{assume}
Suppose that $\gamma \in X^+$ is in an internal down alcove and
$\tau \in X^+$ lies on the horizontal wall 
above $\gamma$. Let $w \in W_p$ with $w \ne 1$.
If $\nu=w.\gamma$ is in
a dominant down alcove such that $w.\tau$ lies on the horizontal wall above
$\nu$ then $\Hom(\nabla(\gamma),\nabla(\nu))=0$.
\end{ass}
   
Assumption \ref{assume} will follow from the results in Sections
\ref{homstart}--\ref{exception},
which will thus complete the proof of Theorem \ref{onedim}. (Note
that the second part of Theorem \ref{onedim} is clear from the
translation arguments we have used.) Consequently, in what follows we
may only apply Theorem \ref{onedim} to weights that are strictly
smaller than $\lambda$. \qed

%



\section{Composites of Carter-Payne maps}\label{cpcomp}

In this section we give a recursive description of the homomorphisms
constructed by Carter and Payne (for $\SL_3$). Using this we can
determine inductively which composites of such maps are non-zero. 


\begin{lemma} 
All factors in the $p$-filtration of $\nabla(\lambda)$ are in
distinct $\nabla_p$-classes if and only if we are in one of the
following situations:
\begin{enumerate}
\item[(i)]{$\lambda$ lies in a down alcove and either
the alcove is just dominant or both}
\begin{enumerate}
\item[(a)]{$\langle \lambda'' + \rho,
\check{\alpha}_i \rangle \notequiv 1\pmod p$ for $i=1,2$.} 
\item[(b)]{$\langle \lambda'' + \rho, \check{\alpha}_1+ \check{\alpha}_2\rangle
\notequiv 1\pmod p$.}
\end{enumerate} 
\item[(ii)]{$\lambda$ lies in an up alcove and either
the alcove is just dominant or both}
\begin{enumerate}
\item[(a)]{$\langle \lambda'' + \rho,
\check{\alpha}_i \rangle \notequiv 0\pmod p$ for $i=1,2$.}
\item[(b)]{$\langle \lambda'' + \rho, \check{\alpha}_1+ \check{\alpha}_2\rangle
\notequiv 1\pmod p$.}
\end{enumerate} 
\item[(iii)]{$\lambda$ lies on a left-hand diagonal
wall and either is just dominant or $\langle \lambda'' + \rho,
\check{\alpha}_1 \rangle \notequiv 0\pmod p$.} 
\item[(iv)]{$\lambda$ lies
on a right-hand diagonal wall and either is just dominant or $\langle
\lambda'' + \rho, \check{\alpha}_2 \rangle \notequiv 0\pmod p$.} 
\item[(v)]{$\lambda$ lies on a horizontal wall and either is just
dominant or $\langle \lambda'' + \rho,
\check{\alpha}_1+\check{\alpha}_2 \rangle \notequiv 1\pmod p$.}
\item[(vi)]{$\lambda$ lies on a vertex.}
\end{enumerate}
\end{lemma} 
\pf This is an elementary calculation using the description of
$p$-filtrations given in the previous section if $p \ge 5$. 
If $p=3$ a little more care is needed for the up and down alcove case
as there are only two $p$-restricted weights which lie inside alcoves.
But looking at the $\lambda''$ that can occur it is clear that we
get the same result as for $p \ge 5$ as the $\lambda''$ are mostly too
close together to be in the same $G$-block. 
\qed

We will say that $\lambda $ is \emph{generic} if it satisfies the
conditions of this lemma. We will say that $\lambda$ is {\it
sufficiently generic} unless either
\begin{enumerate}
\item[(i)]{$\lambda$ lies on a diagonal wall and is not generic, or}
\item[(ii)]{$\lambda$ lies in a down alcove with
$\langle \lambda'' + \rho, \check{\alpha}_i\rangle
\equiv 1\pmod p$ for $i=1$ or $2$.}
\end{enumerate}
We will say that $\lambda$ is \emph{recursively generic} if $\lambda$
is sufficiently generic and either all terms in the $p$-filtration of
$\nabla(\lambda)$ are simple, or all such terms are of the form
$\nabla_p(\lambda_i)$ with $\lambda_i''$ recursively generic.

When $\lambda$ is non-generic, we will have to consider in detail the
structure of extensions between factors in the $p$-filtration in
the same $\nabla_p$-class. The basic properties of these extensions
are summarised in the following lemma.

\begin{lemma}\label{nongen}
Suppose that $\nabla_p(\lambda_i)$ and $\nabla_p(\lambda_j)$ are two
factors in the $p$-filtration of $\nabla(\lambda)$ in the same
$\nabla_p$-class with $\lambda_i> \lambda_j$ and there is a non-split
extension appearing between them in $\nabla(\lambda)$. Then
$$\Ext^1_G(\nabla_p(\lambda_i),\nabla_p(\lambda_j))\cong k$$
and the non-split extension is isomorphic
to $(T_{\mu}^{\lambda_i''}\nabla(\mu))\frob\otimes L(\lambda_i')$ (where
$\mu$ lies on the wall between $\lambda_i''$ and $\lambda_j''$) and has
simple socle $L(\lambda_j)$. This extension also has simple head
provided $\lambda_i''$ is
not just dominant. 
\end{lemma}
\pf By considering the various cases that can arise, it is easy to
verify that $\lambda_i''$ and $\lambda_j''$ are related by a single
left or right reflection. 
By \cite[Lemma 4.2]{par1} we have
$$\Ext^1_G(\nabla_p(\lambda_i),\nabla_p(\lambda_j))\cong
\Ext^1_G(\nabla(\lambda_i''),\nabla(\lambda_j''))$$ and this latter
Ext-group is $k$ 
using the corresponding $\SL_2(k)$
result \cite[Theorem 3.6 and Corollary 4.3]{erd}.


Now consider the translate $T_{\mu}^{\lambda_i''}\nabla(\mu)$ where
$\mu$ lies on the wall between $\lambda_i''$ and
$\lambda_j''$. Without loss of generality we can pick $\mu$ so that 
$T_{\mu}^{\lambda_i''}= \pr_{\lambda_i''}( - \otimes E)$
or $T_{\mu}^{\lambda_i''}= \pr_{\lambda_i''}( - \otimes E^*)$
where $E=\nabla(1,0)$, the natural representation and 
depending on whether $\lambda_i''$ and $\lambda_j''$ are related by a
right hand or left hand reflection, respectively.
We want to check that this translate
$T_{\mu}^{\lambda_i''}\nabla(\mu)$
has simple head and socle. Firstly this translate has short exact
sequence
$$0 \to \nabla(\lambda_j'') \to T_{\mu}^{\lambda_i''}\nabla(\mu)
\to \nabla(\lambda_i'') \to 0.$$
So its socle is contained in $L(\lambda_i'') \oplus L(\lambda_j'')$
and its head is contained in $\hd(\nabla(\lambda_i'')) \oplus
\hd(\nabla(\lambda_j''))$ where $\hd$ denotes 
the head of a module. Thus it is enough to show that one of each of
these respective simple modules cannot be in the socle (or head).

If $\lambda_i''$ and $\lambda_j''$ are in adjacent alcoves then 
$\mu$ is in the closure of the facet containing $\lambda_i''$ so then 
this translate has
simple socle $L(\lambda_j'')$ and simple head by \cite[II, 7.19 Proposition
(b)]{jantzen}. 

\begin{figure}[ht]
\centerline{\epsffile{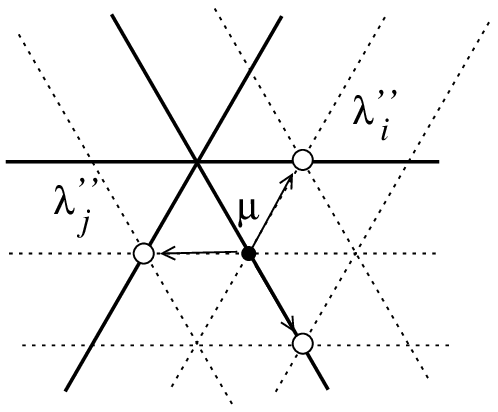}}
\caption{\label{nongenwalls}}
\end{figure}

We now suppose that $\lambda_i''$ and $\lambda_j''$ are related by a
right hand reflection and lie on walls as shown in Figure
\ref{nongenwalls}. 
The argument for the other cases are similar.

Now 
$$\Hom(L(\lambda_i''), T_{\mu}^{\lambda_i''}\nabla(\mu)) 
\cong \Hom(L(\lambda_i''), \nabla(\mu)\otimes E) 
\cong \Hom(L(\lambda_i'')\otimes E^*, \nabla(\mu))$$ and $\lambda_i''=
p\lambda_i'''' + (0,p-2)$. 
So $L(\lambda_i'') \otimes E^*$ has
character $\ch(L(\lambda_i'''')\frob\otimes\nabla(0,p-1)) +
\ch(L(\lambda_i'''')\frob\otimes\nabla(1,p-3))$ and neither of these
factors has $G_1$-type the same as that of $L(\mu)$, which is
$(p-1,p-2)$. Thus the last $\Hom$ space must be zero and the socle of
$T_{\mu}^{\lambda_i''}\nabla(\mu)$ is simple.  This implies that the
translate must be indecomposable and hence by uniqueness and
\cite[Lemma 1.2]{erdhenk2}, our desired extension is isomorphic to
$\bigl(T_{\mu}^{\lambda_i''}\nabla(\mu)\bigr)\frob\otimes
L(\lambda_j'')$ and has simple socle $L(\lambda_j)$.

Similarly for the head (when $\mu$ is an internal weight) we have
$$\Hom(T_{\mu}^{\lambda_i''}\nabla(\mu), \hd(\nabla(\lambda_j'')))
\cong \Hom(\nabla(\mu)\otimes E, \hd(\nabla(\lambda_j'')))
\cong \Hom(\nabla(\mu), \hd(\nabla(\lambda_j''))\otimes E^*)$$
and $\hd(\nabla(\lambda_j''))= L(\eta)=L(p \eta'' + (p-1,0)$. So 
$L(\eta) \otimes E^*$ has character
$\ch(L(\eta'')\frob\otimes\nabla(p-1,1)) 
+ \ch(L(\eta'')\frob\otimes\nabla(p-2,0))$
and neither of these factors has $G_1$-type the same as that of 
$\hd(\nabla((\mu)))$,
which is $(0,p-1)$. Thus the last $\Hom$ space must be zero and the
head of $T_{\mu}^{\lambda_i''}\nabla(\mu)$ is simple.

If $\mu$ is a just dominant weight then the head of $\nabla(\mu)$ has
$G_1$-type $(p-2,0)$ so our argument fails.
In fact in this case the translate has non-simple head so the
condition on $\lambda_i''$ in the lemma is necessary. 
\qed


If $\lambda$ and $\mu$ satisfy the conditions of Theorem
\ref{morphism} for some $e>0$, with $\mu$ to the left (respectively to
the right, below) $\lambda$ we denote the corresponding Carter-Payne
map $\phi$ by $\phi^e_L$ (respectively $\phi^e_R$, $\phi^e_B$). To
each of these maps there is a corresponding local (i.e. $e=1$) map
starting at $\lambda$ (which if $\lambda$ is on a wall may by the
identity map). These local maps are $p$-good homomorphisms, and we
refer to the set of terms in the $p$-filtration of $\nabla(\lambda)$
which survive (necessarily completely) under such a map $\phi$ as {\it
the $p$-factors (of $\nabla(\lambda)$) associated to $\phi$}, unless
both $\lambda$ lies in an up alcove or on a horizontal wall, and
$\phi=\phi_B^{>1}$.  In these remaining cases we refer only to the
{\it top term} in the $p$-filtration of $\nabla(\lambda)$ as such a
$p$-factor.

\begin{proposition}\label{mapsare}
Suppose that $\lambda$ and $\mu$ satisfy the conditions of Theorem
\ref{morphism} for some $e>0$. In all cases the Carter-Payne map
$\phi$ is non-zero only on the $p$-factors associated to $\phi$, and
on such a factor $\nabla_p(\lambda_i)$ it is induced via twisting from
the corresponding Carter-Payne map from $\nabla(\lambda_i'')$ about an
$(e-1)$-wall in all but the following cases:
\begin{enumerate}
\item[(i)]{$\lambda$ in a down alcove with $\lambda_3$ and
$\lambda_8$ in the same $\nabla_p$-class, $\phi=\phi_R$, and
$i=3$ or $8$.}
\item[(ii)]{$\lambda$ in a down alcove with $\lambda_5$ and
$\lambda_6$ in the same $\nabla_p$-class, $\phi=\phi_L$, and
$i=5$ or $6$.}
\item[(iii)]{$\lambda$  on a LH diagonal wall with
$\lambda_2$ and $\lambda_3$ in the same $\nabla_p$-class,
$\phi=\phi_L$, and $i=2$ or $3$.}
\item[(iv)]{$\lambda$ on a RH diagonal wall with
$\lambda_2$ and $\lambda_3$ in the same $\nabla_p$-class,
$\phi=\phi_R$, and $i=2$ or $3$.}
\end{enumerate}
In each of these four cases the map induced on $\nabla_p(\lambda_i)$
is non-zero, but does not come from the twist of the corresponding
Carter-Payne map.
\end{proposition} 

Thus if $\lambda$ is recursively generic then the composition factors
occurring in the image of $\nabla(\lambda)$ under the Carter-Payne
homomorphism $\phi$ can be described inductively using the local data
from the preceding section. In the remaining cases the above
description will be sufficient for our purposes, so we will not
analyse the exceptional cases further. 

\pf We proceed by induction on $e$. When $e=1$ we are done by the
results in the preceding section.  Now suppose that $e>1$. In what
follows we shall assume when presenting lattices in $\nabla(\mu)$ that
$\mu$ is not too close to the edge of the dominant region. The reader
may verify that the modifications necessary in the remaining cases do
not affect the form of the answer.


\noindent{\bf Case (i): Down alcoves.} We will
begin by consider the case when $\lambda$ is as in Figure
\ref{alcoves}a(ii). There are three subcases, corresponding to
reflections about a wall to the left, right or below the alcove
containing $\lambda$. We consider a reflection to a
weight $\mu$ to the right of $\lambda $ (i.e. corresponding when
$e=1$ to the map $\nabla(1)\to \nabla(2)$) and denote the
corresponding morphism (which will be unique up to scalars 
once we know the $\Hom$ space is one-dimensional) by $\phi$. The
case of a reflection to the left is similar, while that for a
reflection below is even simpler; both are left to the
reader.
  
\begin{figure}[ht]
\centerline{\epsffile{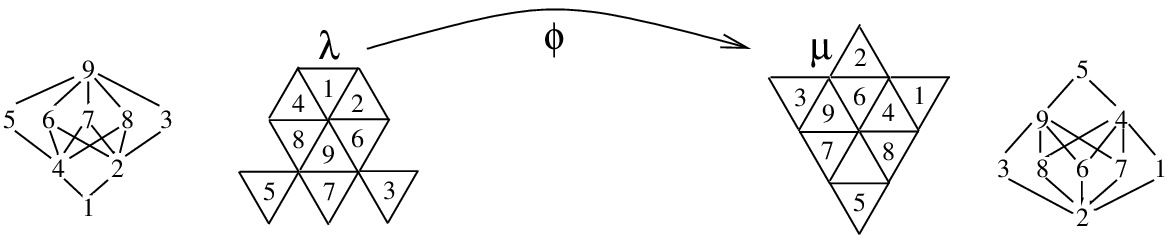}}
\caption{\label{First}}
\end{figure}

In the case under consideration $\mu$ lies in an up alcove,
and corresponds to the situation in Figure \ref{alcoves}b.  We will
renumber the alcoves in this latter diagram to be consistent with the
numbering for $\lambda$.  That is, we will renumber the alcoves so
that alcoves with the same number in each diagram are in the same
$\nabla_p$-class. Note that as $e>1$ this numbering is different from
that in Figure \ref{alcoves}. If $\lambda$ is generic then this is
unambiguous, and we will use the same labelling in the non-generic
case.  The new labelling is illustrated in Figure \ref{First},
together with the corresponding lattices. To
avoid confusion we will refer to a weight in $W.\lambda$ in an alcove
labelled $i$ in the diagram for $\lambda$ as $\lambda_i$ and in the
diagram for $\mu$ as $\mu[\lambda_i]$.

If $\lambda$ is non-generic then the only terms in the 
$p$-filtration which can lie in the same $\nabla_p$-class are 
$\lambda_1$ and $\lambda_7$, or $\lambda_3$ and $\lambda_8$, or $\lambda_5$
and $\lambda_6$.  As we are considering the case in Figure
\ref{alcoves}a(ii), $\lambda_3$ and $\lambda_8$ are in distinct
$\nabla_p$-classes. Also, for any $\lambda$, the terms in the
$p$-filtration labelled by $\lambda_1$, $\lambda_4$ and $\lambda_5$
must lie in the kernel of $\phi$ as none of their composition
factors lie above the socle of $\nabla_p(\lambda_2)$ (by figure
\ref{alcoves}a) and some composition factor of $\nabla_p(\lambda_2)$
maps to the socle of $\nabla(\mu)$ (as this is the only factor in
the same $p$-class as the lowest term in the filtration of $\nabla(\mu)$).

In a similar fashion we see that the image of $\phi$ does not
involve the terms in the $p$-filtration of $\nabla(\mu)$ labelled by
$\mu[\lambda_1]$, $\mu[\lambda_4]$ and $\mu[\lambda_5]$, as the
image of the head of $\nabla(\lambda)$ is a composition factor of
$\nabla_p(\lambda_9)$.  Thus the six terms in the filtration of
$\nabla(\lambda)$ that may survive under $\phi$ are in distinct
$\nabla_p$-classes. 

By the above remarks, $\phi$ induces a (non-zero) map $\phi'$ on the
quotient of
$\nabla(\lambda)$ by the submodule with a filtration by
$\nabla_p(\lambda_1)$, $\nabla_p(\lambda_4)$ and
$\nabla_p(\lambda_5)$. We shall denote this quotient module by
$\nabla(\lambda)/[\lambda_1,\lambda_4,\lambda_5]$, and other quotients
similarly. Now consider the restriction of $\phi'$ to 
$\nabla_p(\lambda_2)$.  As this module is the only part of the filtration with
any composition factors in common with those of 
$\nabla_p(\mu[\lambda_2])$ (which contains the socle of $\nabla(\mu)$) and
$\phi'$ is non-zero, this restriction must also be non-zero.  Hence
we obtain a non-zero map from $\nabla_p(\lambda_2)$ to 
$\nabla_p(\mu[\lambda_2])$.  By (\ref{homreduce}) we have
$$\Hom (\nabla_p (\lambda_2) ,\nabla_p ( \mu [\lambda_2])) \cong\Hom
(\nabla (\lambda_2'') ,\nabla ( \mu [\lambda_2]'')). $$ 
A simple calculation shows that this latter pair of weights satisfy
the conditions of Theorem \ref{morphism} for $e-1$. 
Using the inductive hypothesis
this $\Hom$-space is one-dimensional, and hence any such homomorphism is
unique (up to scalars).  Therefore, by induction, we can describe the
composition factors occurring in the image of such a homomorphism (up
to the exceptional ambiguity in Theorem \ref{morphism}).  As
the isomorphism of $\Hom$-spaces is induced by the map $V \longmapsto V\frob
\otimes L (\lambda_2') $, we can thus describe the composition factors
occurring in the image of our restriction of $\phi'$.

We next consider the map obtained by following $\phi$ by the
quotient map from $\nabla (\mu) $ which kills $\nabla_p (\mu
[\lambda_2]) $.  As the image under $\phi $ of $\nabla_p (\lambda_2)
$ is killed by this quotient map, this induces a map $\phi$ from
$\nabla (\lambda)/[ \lambda_1, \lambda_2, \lambda_4, \lambda_5] $ to
$\nabla (\mu)/[\mu [\lambda_2]] $.  We wish to argue as above and show
that the restriction of $ \phi$ to each of the submodules $\nabla_p
(\lambda_i) $, for $i \in \{3,6,7,8\}$, can be determined by induction
from a reflection about an ($e-1$)-wall.  For this it is enough to
show that the restriction is non-zero, as then the argument given
above also holds in this case.  Once we have shown this, the only
remaining factor in $\nabla (\lambda) $ is $\nabla_p (\lambda_9) $ and
a similar argument to the above shows that $\phi $ induces a
(non-zero) map from this factor to $\nabla_p (\mu [\lambda_9]) $,
which again can be determined by induction.

Thus it only remains to show that the map induced on $\nabla_p
(\lambda_i) $ is non-zero for $i\in\{3,6,7,8\}$.  Suppose that this
map is zero, and hence that this factor is in the kernel of the
original map $\phi $.  As the image under $\phi$ of $\nabla_p
(\lambda_2) $ is non-zero, and this latter module has a simple head,
this head must survive.  Hence the submodule $E_i$ of $\nabla
(\lambda)/[\lambda_1, \lambda_4, \lambda_5] $ with a filtration by
$\nabla_p (\lambda_2) $ and $\nabla_p (\lambda_i) $ must have a non
simple head. We will show that the head is always simple,
contradicting our assumption that the induced map is zero.

When $i=3$ (respectively $i=7$) the head of $E_i$ is simple as there
is a $p $-good homomorphism into $\nabla (\lambda) $ from some $\nabla
(\tau) $ (which has a simple head) such that the image contains $E_i$
as a quotient.  For the remaining cases will imitate an argument in
\cite[pages 360-1]{par1}.  As the head of $E_i$ can contain at most
two simple modules, it will be enough by \cite[II, 3.16 (3)]{jantzen} to
show that $\Hom_{G_1} (E_i,L( \lambda_2')) = 0 $.  From the defining
sequence for $E_i$ (and writing $L=L( \lambda_2') $) we have the exact
sequence
$$0 \to \Hom_{G_1}( \nabla_p (\lambda_i), L) \to\Hom_{G_1} (E_i,L)
\to\Hom_{G_1}( \nabla_p (\lambda_2), L) \to \Ext_{G_1}^1( \nabla_p
(\lambda_i), L) $$ 
which by (\ref{infext}) and the fact that
$\lambda_2$ and $\lambda_i$ are in distinct $\nabla_p$-classes gives
$$0 \to\Hom_{G_1} (E_i,L) \to \Delta (\lambda_2''^* )\frob
\stackrel{\theta}{\to} \Delta (\lambda_i''^*)\frob \otimes
\Ext_{G_1}^1(L( \lambda_i'), L).$$
 
Thus we will be done if we can show that $\theta $ is an embedding.
First suppose that $\theta = 0 $.  Then $\Hom_{G_1} (E_i,L) =\Delta
(\lambda_2''^* )\frob $ and this implies that $E_i$ is semisimple as a
$G_1$-module. By \cite[II, 3.16 (3)]{jantzen} we see that $\soc_G (E_i) $
is not simple, which contradicts the fact that $E_i$ embeds in $\nabla
(\lambda_2) $.  Therefore we must have that $\theta \neq 0 $.  To show
that $\theta $ is an embedding, we will show that
\begin{equation} \label{phiis}
\theta^{*(-1)}: \nabla (\lambda_i'') \otimes (\Ext_{G_1}^1(L(
                     \lambda_i'), L)^{(-1)} )^* \to \nabla
                     (\lambda_2'')
\end{equation}
is onto. 

When $i=6$ the results of Yehia \cite{yehia} summarised in
\cite[Proposition 4.1]{par1} give that the $\Ext^1 $-group in
(\ref{phiis}) is isomorphic to $k $, and $\lambda_6'' = \lambda_2'' $.
As $\Hom (\nabla (\lambda_6''), \nabla (\lambda_2'')) \cong k $ we
deduce that $\theta $ is an embedding.

When $i=8$, \cite[Proposition 4.1]{par1} gives that the $\Ext^1
$-group in (\ref{phiis}) is isomorphic to $\nabla (1,0)\frob$ and
$\lambda_2''= (a,b-1)$, $\lambda_8'' =(a-1,b)$ for some $a,b > 0 $.
Thus $\theta^{*(-1)}$ is in 
\begin{equation}
\label{onto}
\begin{array}{ll}
\Hom (\nabla (a-1,b)\, \otimes& \nabla (1,0)^*,\nabla (a,b-1))\\
&\cong\Hom (\nabla (a-1,b),\nabla (1,0)\otimes\nabla
(a,b-1))\\
&\cong\Hom (\nabla (a-1,b),\pr_{(a-1,b)}\nabla (1,0)\otimes\nabla (a,b-1))\\
&\cong\Hom (\nabla (a-1,b),T_{(a,b-1)}^{(a-1,b)}\nabla (a,b-1))\cong
k
\end{array}
\end{equation}
where the final isomorphism follows because
$T_{(a,b-1)}^{(a-1,b)}\nabla (a,b-1)\cong \nabla (a-1,b)$ (as $b
\notequiv 0\ (\mod p)$). But there is an obvious surjection in the
Hom-space containing $\theta^{*(-1)}$, as $\nabla (a-1,b) \otimes
\nabla (1,0)^*$ has a good filtration with quotient $\nabla (a,b-1)$,
and hence $\theta^{*(-1)}$ is surjective as required.

Next consider the cases when either $a''\equiv 0\ (\mod p)$ or
$b''\equiv 0\ (\mod p)$. It is easy to see that the argument above is
unaffected by the former, so we only need consider the case $b''\equiv
0\ (\mod p)$. We are in the situation where $\lambda$ is as in
Figure \ref{alcoves}a(iii) and $\mu$ is as in Figure
\ref{alcoves}b(iii).  

Just as above, we consider the map induced by $\phi$ on various
subquotients of $\nabla(\lambda)$. Everything goes through unchanged
except for the cases involving $\nabla_p(\lambda_3)$ and
$\nabla_p(\lambda_8)$. 

Let $\nabla_p(\tau,\nu)$ denote a non-split extension of $\nabla_p(\nu)$ by
$\nabla_p(\tau)$, and similarly for $\nabla(\tau,\nu)$. (We will only
apply this notation in situations where the extension is unique.)  By
assumption, we have to consider the map induced by $\phi$ from
$\nabla_p(\lambda_8,\lambda_3)$ to
$\nabla_p(\mu[\lambda_3],\mu[\lambda_8])$. By Lemma \ref{nongen} this
notation is well-defined, and each of the modules has simple
socle. Arguing as in the generic case we see that the restriction of
$\phi$ to $\nabla_p(\lambda_3)$ must be non-zero. (It is clearly
non-zero on $\nabla_p(\lambda_8)$.) 
There is
an obvious homomorphism from $\nabla_p(\lambda_8,\lambda_3)$ to
$\nabla_p(\mu[\lambda_3],\mu[\lambda_8])$ induced by the Carter-Payne
homomorphism from $\nabla(\lambda_8'')$ to $\nabla(\mu[\lambda_8]'')$.
But the restriction of this map to $\nabla_p(\lambda_3)$ is zero, and
hence it cannot be the map we require. Thus the map induced by $\phi$
on $\nabla(\lambda_8'')$ cannot be the twist of the corresponding
Carter-Payne map. As the socle of
$\nabla_p(\mu[\lambda_3],\mu[\lambda_8])$ is simple it is clear that
the same is true of the map induced by $\phi$ on $\nabla(\lambda_3'')$.
%
%
%
This completes the proof in the case of reflection to the right from a
down alcove. Clearly, reflection to the left is entirely
analogous. 

For the reflection $\phi$ below $\lambda$ the argument is much more
straightforward. The only term in the $p$-filtration of
$\nabla(\lambda)$ which can survive under $\phi$ is
$\nabla_p(\lambda_9)$, and this must map into
$\nabla_p(\mu[\lambda_p])$. This map is induced by a Carter-Payne map
from $\nabla(\lambda_9'')$ to $\nabla(\mu[\lambda_9]'')$, and hence is
known by induction.



\begin{figure}[ht]
\centerline{\epsffile{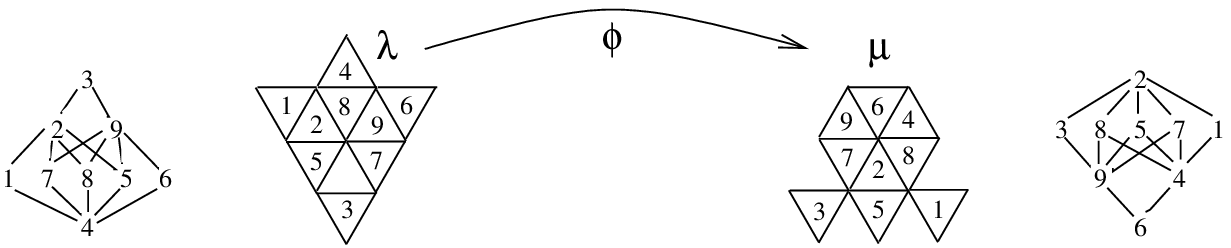}}
\caption{\label{Second}}
\end{figure}

\noindent{\bf Case (ii): Up alcoves.} We begin by considering
reflection to the right of $\lambda$.
As in case (i), we renumber the terms in the $p$-filtration of
$\nabla(\mu$) by their $\nabla_p$-classes with respect to
$\nabla(\lambda)$. (Again we use the generic labelling even in the
non-generic case.) This is illustrated in Figure \ref{Second} (in the
generic case). 

If $\lambda$ is generic then arguing as in case (i) we see that all
but the terms labelled by $\lambda_3$, $\lambda_6$ and $\lambda_9$
must lie in the kernel of $\phi$, and the image is contained in the
submodule filtered by the terms labelled by $\mu[\lambda_3]$,
$\mu[\lambda_6]$, and $\mu[\lambda_9]$. Arguing as in case (i) we
see that $\phi$ is non-zero on each of the three subquotients, and
corresponds to a Carter-Payne map from $\nabla(\lambda_i'')$ to
$\nabla(\mu[\lambda_i]'')$, which are known by induction.

We wish to argue that the same is true in the non-generic case. Here
it is possible that the term labelled by $\lambda_5$ may survive under
$\phi$, and/or the term labelled by $\mu[\lambda_8]$ may contain part
of the image. In the former case we must have that
$\nabla_p(\lambda_5)$ maps into $\nabla_p(\mu[\lambda_6])$ (not
$\nabla_p(\mu[\lambda_5])$ as this would contradict the fact that the
image of $\nabla(\lambda)$ has simple head). This implies the
existence of a non-zero homomorphism from $\nabla(\lambda_5'')$ to
$\nabla(\mu[\lambda_6]'')$. But $\mu[\lambda_6]''\not\leq
\lambda_5''$, so this is impossible. Similarly, if $\mu[\lambda_8]$
intersects $\im\phi$ then this corresponds to a non-zero homomorphism
from $\nabla(\lambda_3'')$ to $\nabla(\mu[\lambda_8]'')$, but
$\mu[\lambda_8]''\not\leq \lambda_3''$. Thus the same argument as for
the generic case holds here.

\begin{figure}[ht]
\centerline{\epsffile{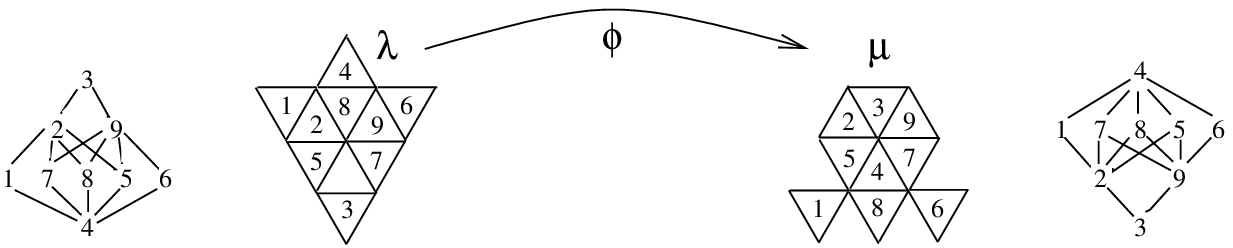}}
\caption{\label{Third}}
\end{figure}

So far each of the cases we have considered has corresponded to the
corresponding result when $e=1$. However we shall see in the next
case, the reflection below an up alcove, that this is not always the case.
Let $\lambda$ be in an up alcove, and $\phi$ be a reflection below
$\lambda$. In Figure \ref{Third} we have (as usual) labelled the terms in the
$p$-filtrations using the generic labels for $\lambda$.

In the generic case we see that $\phi$ must kill all but the top
quotient $\nabla_p(\lambda_3)$, and corresponds to a Carter-Payne map
from $\nabla(\lambda_3'')$ to $\nabla(\mu[\lambda_3]'')$. {\it Note that
this differs from the case $e=1$, when four terms in the filtration survive.}
Clearly this map also exists in the non-generic case, but we also have
the possibility that (as for $e=1$) the terms labelled by $\lambda_2$,
$\lambda_3$, $\lambda_8$ and $\lambda_9$ are not in the kernel of
$\phi$. 

This would contradict the one-dimensionality of Hom-spaces, but we
cannot use this to eliminate the possibility as we have not
established that result for $\nabla(\lambda)$ at this stage. However,
such a map would induce a non-zero homomorphism from
$\nabla_p(\lambda_3)$ to $\nabla_p(\mu[\lambda_8])$, and hence a
non-zero homomorphism from $\nabla(\lambda_3'')$ to
$\nabla(\mu[\lambda_8]'')$. As we are in the non-generic case, it is
easy to verify that $\gamma=\lambda_3''$ and $\nu=\mu[\lambda_8]''$
satisfy the conditions in Assumption \ref{assume}. Hence this
homomorphism cannot exist, as Assumption \ref{assume} follows from the
results in Sections \ref{homstart}--\ref{exception} which we know
holds for weights smaller that $\lambda$ by induction. Therefore there
are no extra homomorphisms in the non-generic case.



\begin{figure}[ht]
\centerline{\epsffile{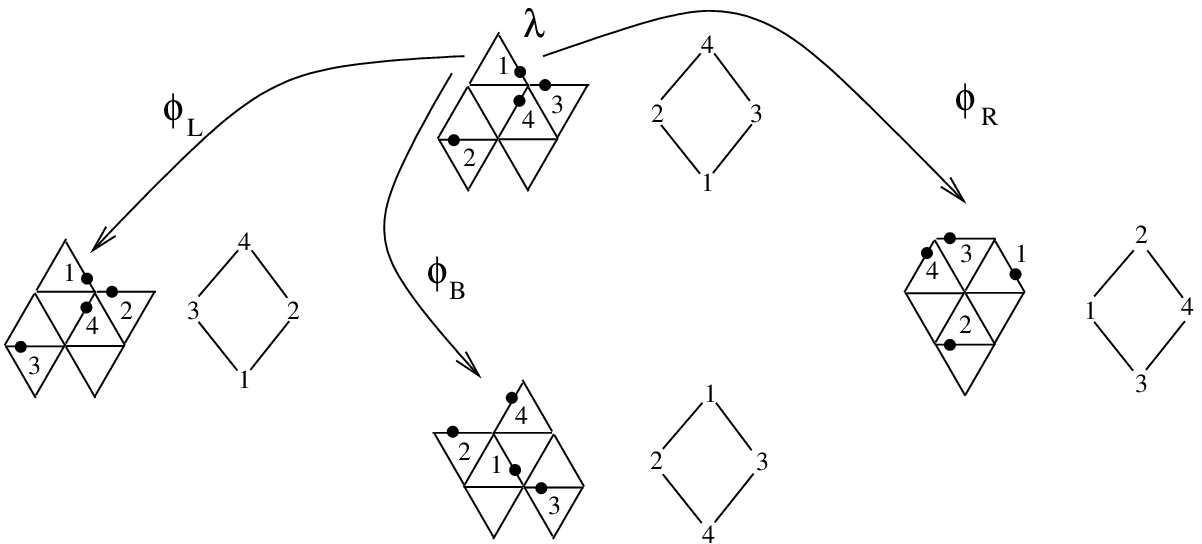}}
\caption{\label{Fourth}}
\end{figure}

\noindent{\bf Case (iii): Diagonal walls.} 
We shall consider the case
where $\lambda$ lies on a wall as in Figure \ref{Fourth}; the other
case is symmetric. As usual there are three possible maps from
$\nabla(\lambda)$: to the left ($\phi_L$), right ($\phi_R$) and
below ($\phi_B$). As usual we label the terms in the
$p$-filtration of the target module $\nabla(\mu)$ by their (generic)
$\nabla_p$-classes with respect to the labelling of
$\nabla(\lambda)$. The three cases are illustrated in Figure
\ref{Fourth}, in the generic case.

When $\lambda$ is generic it is easy to verify as above that
reflections to the right (respectively below) correspond to
Carter-Payne maps (which are known by induction) from the two
(respectively one) term(s) in the $p$-filtration that survive under
the corresponding $e=1$ map. The same is true for reflection to the
left {\it provided one realises that the corresponding $e=1$ map here
is the identity morphism, i.e that all four terms survive.}

For non-generic $\lambda$ the same results hold, but the proofs
require more care. First consider reflection to the right. The only
possible difference here is that $\nabla_p(\lambda_2)$ might survive
under $\phi_R$. However, $\mu[\lambda_3]''\not\leq\lambda_2''$, and
so this is impossible (as in the case of reflection to the right from
an up alcove).

Next consider reflection to the left. The argument here is just as for
the corresponding reflection from a down alcove.

Finally, the argument for reflections below $\lambda$ is unchanged in
the non-generic case, as the only pair of terms in the $p$-filtration
which might be affected are already in the kernel of the map.



\begin{figure}[ht]
\centerline{\epsffile{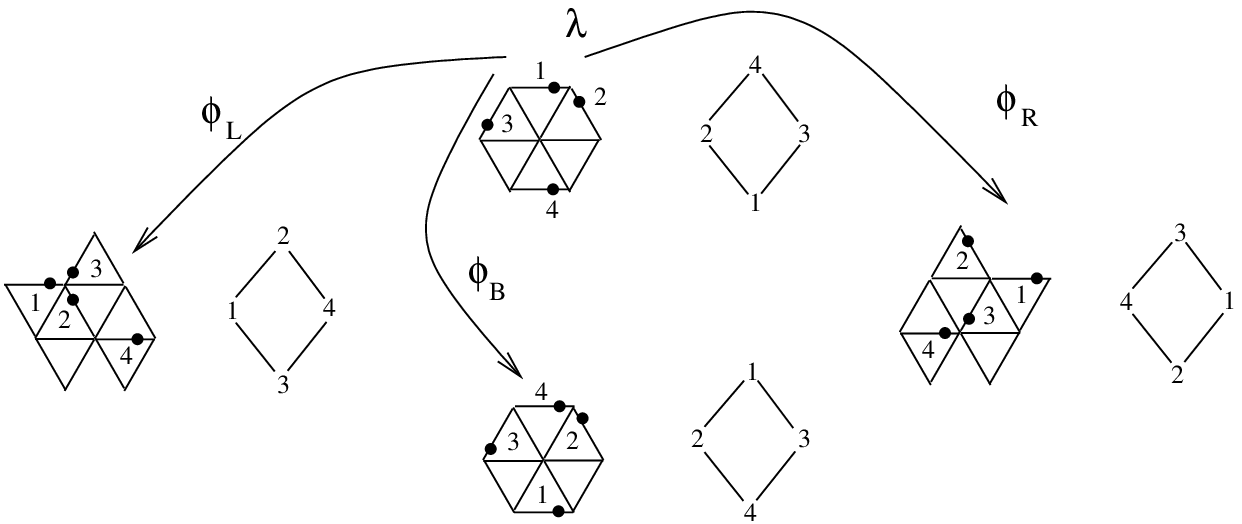}}
\caption{\label{Fifth}}
\end{figure} 

\noindent{\bf Case (iv): Horizontal walls.}  We conclude by
considering the case where $\lambda$ lies on a horizontal wall as in
Figure \ref{Fifth}. As usual there are three possible maps: to the
left ($\phi_L$), right ($\phi_R$) and below ($\phi_B$). As in
the previous cases we number the terms in the $p$-filtration of
$\nabla(\mu)$ by their (generic) $\nabla_p$-classes with respect to
the labelling for $\nabla(\lambda)$. By symmetry we only need consider
reflections to the right and below.

As in the preceding case, when $\lambda$ is generic it is routine to
verify that reflections to the right correspond
to Carter-Payne maps from the two terms in the
$p$-filtration that survive under the corresponding $e=1$ map. For
reflections below the map corresponds to the Carter-Payne map from the
top term of the $p$-filtration of $\nabla(\lambda)$ into $\nabla_p(\mu)$.

When $\lambda$ is non-generic and $\phi$ is a reflection to the
right or down then mimicking the argument for the corresponding
reflection from an up alcove in case (ii) we see that the result is
the same as in the generic case.
\qed



We wish to determine when the composition of two Carter-Payne
homomorphisms is non-zero.  To describe this it will be convenient to
introduce the following notation.  Given a weight $\lambda$, we
denote by $\nabla_p(\tilde{\lambda})$ the term in the $p$-filtration of
$\nabla(\lambda)$ containing the head of $\nabla(\lambda)$.  
(Note
that if $\lambda =p\lambda'' + \lambda'$ where $\lambda'' = (a,b)$ 
with $a>0$ and $b>0$, then 
$\tilde{\lambda} =p(\lambda''-\rho) +w_a \lambda'$ 
where $w_a \lambda' $ is the reflection of $\lambda'$ in the
horizontal wall immediately above the origin.) 

Recall that for a dominant weight $\lambda $, we let $\phi_L^i $
(respectively $ \phi_R^i, \phi_B^i $) denote the Carter-Payne
homomorphism from $\nabla ( \lambda) $ about a $p^i $-wall to the left
of (respectively to the right of, below) $\lambda $.  {\it Recall
further that in using this notation we always chose $i$ minimal for
the given wall.} If we do not wish to specify explicitly the direction
of the reflection, we will replace $L,R,B$ in the above by $\dagger $,
and similarly the degree $i$ of a wall by $* $. If the degree has to
be strictly greater than $1$, we replace the $*$ by $>1$.
Given a composition
$\theta^t\cdots\theta^1$ of homomorphisms, we shall refer to any
composition of the form $\theta^u\cdots\theta^1$ with $u<t$ as a {\it
subcomposition} of the given composition.

We next define for each type of facet a set of {\it maximal
compositions}. (The reason for this terminology will become clear in
the next corollary.) Given a weight $\lambda$, we define the product
$\theta^t\cdots\theta^1$, where each $\theta^i$ is of the form
$\phi_{\dagger}^{*,i}$, to mean the composite of the Carter-Payne maps
$\theta^i:\nabla(\lambda^i)\to \nabla(\lambda^{i+1})$ where
$\lambda^1=\lambda$ and $\lambda^{i+1}$ is the weight $\mu$ such that
the Carter-Payne map $\phi_{\dagger}^{*,i}$ from $\nabla(\lambda^i)$
has image in $\nabla(\mu)$. (Here we interpret any Carter-Payne map
which starts in or leaves the dominant region as the zero map.) Now
the set $\max(\lambda)$ of maximal compositions for $\lambda$ depends
on the type of facet in which $\lambda$ lies:
\medskip\medskip

\centerline{\begin{tabular}[pos]{|p{2.0in}|p{3.4in}|}
\hline
 \begin{minipage}{1.85in}\centerline{$\lambda$$\phantom{^{1^1}}$$\phantom{\underbrace{\!}}$  }
\end{minipage}
&\begin{minipage}{1.85in}\centerline{\
max($\lambda$)$\phantom{^{1^1}}$$\phantom{\underbrace{\!}}$}
\end{minipage}\\ \hline Internal down
alcove$\phantom{^{1^1}}$$\phantom{\underbrace{\!}}$ & $\{\phi_B^*, \
\phi_R^*\phi_L^*\phi_R^*, \ \phi_L^*\phi_B^1\phi_R^*,\
\phi_L^*\phi_R^*\phi_L^*,\ \phi_R^*\phi_B^1\phi_L^*\}$\\ \hline Up
alcove$\phantom{^{1^1}}$$\phantom{\underbrace{\!}}$ & $\{\phi_B^*,\
\phi_L^*\phi_R^*, \ \phi_R^*\phi_L^{>1\!}\phi_R^*,\ \phi_R^*\phi_L^*,\
\phi_L^*\phi_R^{>1\!}\phi_L^*,\ \phi_R^*\phi_B^1,\
\phi_L^*\phi_B^1\}$\\ \hline RH just dominant down alcove&
$\{\phi_R^*\phi_L^{>1\!},\ \phi_L^*\}$$\phantom{^{1^1}}$$\phantom{\underbrace{\!}}$ \\
\hline LH just dominant down alcove&
$\{\phi_L^*\phi_R^{>1\!},\ \phi_R^*\}$$\phantom{^{1^1}}$$\phantom{\underbrace{\!}}$ \\
\hline Horizontal wall$\phantom{^{1^1}}$$\phantom{\underbrace{\!}}$ &
$\{\phi_B^*,\ 
\phi_R^*\phi_L^*\phi_R^*, \ \phi_L^*\phi_R^*\phi_L^*\}$\\ \hline RH
diagonal wall$\phantom{^{1^1}}$$\phantom{\underbrace{\!}}$ &
$\{\phi_L^*\phi_B^*, \ \phi_L^*\phi_R^*\phi_L^*, \ \phi_R^*\phi_L^*\phi_R^*, \
\phi_L^*\phi_B^*\phi_R^*\}$\\ \hline LH diagonal
wall$\phantom{^{1^1}}$$\phantom{\underbrace{\!}}$ & $\{\phi_R^*\phi_B^*, \
\phi_R^*\phi_L^*\phi_R^*, \ \phi_L^*\phi_R^*\phi_L^*, \ \phi_R^*\phi_B^*\phi_L^*\}$\\
\hline
\end{tabular}}
\medskip
\medskip

The following corollary allows us to determine inductively when a
composite of Carter-Payne maps can be non-zero. As all 
$\nabla(\lambda)$ have simple head, condition (i) is clearly
both necessary and sufficient. However, the content of this lemma is
that we can impose the additional condition (ii).

\begin{corollary} \label{composite} Let $\lambda^1,\ldots,\lambda^u$ be a
sequence of distinct dominant weights such that for each $i$ we have a
Carter-Payne homomorphism $\phi_{\dagger}^{*,i}: \nabla (\lambda^i) \to
\nabla (\lambda^{i+1}) $.  Then the composite
$\phi_{\dagger}^{*,u}\cdots \phi_{\dagger}^{*,1}$ is
non-zero if and only if both
\begin{enumerate} 
\item[(i)]{the induced homomorphism
$$\phi_{\dagger}^{*,u}\cdots \phi_{\dagger}^{*,1}: \nabla_p
(\widetilde{\lambda^1}) \to \nabla_p
(\lambda^u[\widetilde{\lambda^1}]) $$
is non-zero, and}
\item[(ii)]{the composite $\phi_{\dagger}^{*,u}\cdots
\phi_{\dagger}^{*,1}$ is a subcomposition of a composition in
$\max(\lambda)$.}
\end{enumerate}
\end{corollary} 
\pf As noted above, the first condition is itself both necessary and
sufficient. That the second condition must also be satisfied for a
composite to be non-zero follows in most cases from the description of
maps via $p$-factors in Proposition \ref{mapsare}. By considering each
case in turn, we see that $\nabla_p (\lambda^u
[\widetilde{\lambda^1}]) $ is never one of the terms in the $p
$-filtration of $\nabla (\lambda^i) $ for which any ambiguity remains
in our description of the Carter-Payne homomorphisms from $\nabla
(\lambda^i) $. 

There are a few additional cases which are not eliminated by the
$p$-factor description in Proposition \ref{mapsare}.  In the cases
where $\lambda^1$ lies on a diagonal wall we can eliminate the
additional possibility of two consecutive reflections in the same
direction by reducing to the $\SL_2$ case using (\ref{homreduce}) and
\cite[(4.3) Corollary]{erd}, where the required Hom-spaces have been
calculated in \cite[Theorem 5.1]{coxerd}. Finally, for left-hand
(respectively right-hand) diagonal walls we must eliminate the
possibility of a composite of the form
$\phi_R^*\phi_L^*\phi_R^*\phi_L^*$ (respectively
$\phi_L^*\phi_R^*\phi_L^*\phi_R^*$).  But for such a map to exist
there must be a corresponding composite on
$\nabla_p(\tilde{\lambda})$, which is impossible by the above
description of max($\mu$) for each possible type of $\mu$ and
induction on the sum of the degrees of the maps in the quartet.  \qed

\begin{remark} The induced morphism in Corollary \ref{composite}(i)
corresponds to a composite of Carter-Payne homomorphisms
$$\nabla(\widetilde{\lambda^1}'') \to
\nabla(\lambda^u[\widetilde{\lambda^1}]'')$$ 
which has already been determined by induction. 
\end{remark}

We will take advantage of the notation used for elements in
max($\lambda$) to talk about maps starting from differing
weights. Given two weights $\lambda$ and $\tau$ and a composite of
Carter-Payne maps $\phi=\phi_{\dagger}^{*,u}\cdots \phi_{\dagger}^{*,1}$
from $\nabla(\lambda)$, we shall call the map
$\theta_{\dagger}^{*,u}\cdots \theta_{\dagger}^{*,1}$ from $\nabla(\tau)$
with $\theta_{\dagger}^{*,i}=\phi_{\dagger}^{*,i}$ for all $1\leq
i\leq u$ the map from $\nabla(\tau)$ {\it corresponding to $\phi$}.

\section{Determining homomorphisms I: interior weights}\label{homstart}


Suppose that $\lambda$ is not in the lower closure of a just dominant
(up or down) alcove. If $\lambda$ is on a vertex then all
homomorphisms are known inductively by (\ref{homreduce}), so we shall
henceforth assume that this is not the case. In Section
\ref{strategy} and Theorem \ref{onedim} we observed that there was an
obvious upper bound $\Upsilon_{\lambda}$ for the set of weights $\mu$
for which there exist homomorphisms from $\nabla(\lambda)$ to
$\nabla(\mu)$. In this section we will refine these sets to give a
precise description of when homomorphisms occur, in most cases. The
remaining, exceptional cases will be dealt with in Section
\ref{exception}. Those $\lambda$ not of the form above will be
considered in Section \ref{symway}.

We will hereafter identify weights with the facets in which they lie
--- and (when we have a particular linkage class in mind) vice versa
--- and simple modules with their labels.  Before continuing further,
we note a convenient property of the vertex weights that simplifies
certain verifications. We will wish to argue inductively from weights
at vertices in our facet diagrams. For any such weight $\theta$ of the
form $\theta=p\theta''+(p-1)\rho$ we will need to consider the set of
weights $\tau$ such that there is a non-zero homomorphism from
$\nabla(\theta)$ to $\nabla(\tau)$. In principle to calculate such
homomorphisms we need to consider $\theta''$ and its images under the
dot action of $W_p$ (and thus reflect our weights in a different
shifted set of hyperplanes). However, as this corresponds to
translation by $\rho$ (and we have rescaled) we can determine the
weights $\tau$ by using the ordinary dot action of $W_{p^2}$ on
$\theta$ itself.

It will be convenient in what follows to have the following easy
consequence of the definition of max($\lambda$) and the associated
Corollary \ref{composite}.

\begin{lemma}\label{lr}
Let $\lambda$, $\mu$, and $\nu$ be dominant weights such that there
exist Carter-Payne maps $\phi_L^a:\nabla(\lambda)\to\nabla(\mu)$ and
$\phi_R^b:\nabla(\mu)\to\nabla(\nu)$ for some $a,b>0$. Then the
composite map
$$\phi_R^b\phi_L^a:\nabla(\lambda)\to\nabla(\nu)$$
is non-zero. A similar result holds with the roles of left and right
reversed. 
\end{lemma}
\pf It is enough to consider the case where $\lambda$ does not lie on
a vertex. Regardless of the type of the facet containing $\lambda$,
the map $\phi_R^b\phi_L^a$ occurs as a subcomposition of an element in
max($\lambda$), by inspection (and the dominance of $\mu$). Hence by
Corollary \ref{composite} it is enough to show that the map from
$\nabla_p(\tilde{\lambda})$ to $\nabla_p(\nu[\tilde{\lambda}])$ is
non-zero. But $\tilde{\lambda}''$ and $\nu[\tilde{\lambda}]''$ are
also related by a map of the form $\phi_R^*\phi_L^*$, and hence the
result follows by induction on $a+b$.\qed

\begin{figure}[ht]
\centerline{\epsffile{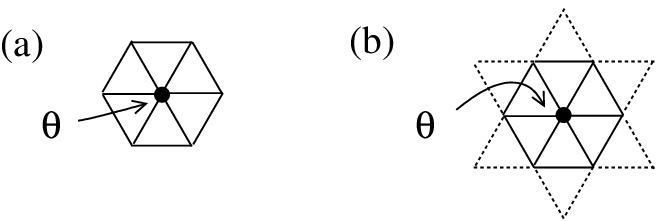}}
\caption{\label{theta}}
\end{figure}

Suppose $\theta \in X^+$ is a vertex weight. We refer to the closure
of the hexagon surrounding $\theta$, as the \emph{neighbourhood} of
$\theta$ (see Figure \ref{theta}(a)). 
We refer to the open star surrounding $\theta$ as the \emph{extended
neighbourhood} of $\theta$ (see Figure \ref{theta}(b)).

The {\it local composition pattern} associated to a weight $\mu$ is
the collection of elements of $W_p.\mu$ labelling terms in a
$p$-filtration of $\nabla(\mu)$. If $\mu \in W_{p^2} . \lambda$  and
$L(\mu)$ is a composition factor of $\nabla(\lambda)$, then
these weights are generically
composition factors of $\nabla(\lambda)$ which may be seen be using
the $p$-filtration of $\nabla(\lambda)$.

For a non-vertex weight $\lambda \in X^+$ and a vertex weight $\theta
\in X^+$, 
we say a dominant weight is a \emph{$\theta$-translate} for
$\nabla(\lambda)$ if it
corresponds to a composition factor of $\nabla(\lambda)$ and it lies   
in the neighbourhood  of $\theta$, if $\lambda$ lies
in the closure of a down alcove, 
and the extended neighbourhood of $\theta$
if $\lambda$ lies in an up alcove. 


Let $\lambda$ be in the closure of an interior down alcove and let
$\eta$ the vertex below $\lambda$. We wish to define a set of eligible
$\theta$-translates for $\nabla(\lambda)$. Unfortunately the
definition is complicated slightly when $\eta$ is on a $p^2$-wall or
$p^2$-vertex, so we first consider the remaining (i.e. $p^2$-regular)
cases.

Let $\theta$ be another vertex dominated by $\eta$ and with $\theta =
w . \eta$ where $w \in W_{p^2}$.  The {\it eligible
$\theta$-translates for $\nabla(\lambda)$} are those highest weights
corresponding to composition factors which are $\theta$-translates for
$\nabla(\lambda)$, and which occur in the local composition pattern
for $\mu = w .\lambda$.  Note that each of these factors lies in the
same $\nabla_p$-class as some weight labelling the $p$-filtration of
$\nabla(\lambda)$.

When $\eta$ is not $p^2$-regular the element $w$ in the definition
above is no longer unique. Let $w$ be the element of $W_{p^2}$
corresponding to the shortest sequence of reflections taking $\eta$ to
$\theta$. Now the {\it eligible $\theta$-translates for
$\nabla(\lambda)$} are those highest weights corresponding to
composition factors which are $\theta$-translates for
$\nabla(\lambda)$, and which occur in the local composition pattern
either for $\mu = w .\lambda$. or for some $\mu^-$ where $\mu^-\leq
\mu$ is obtained from $\mu$ by a sequence of reflections in $W_{p^2}$
all of which fix $\theta$.

\begin{figure}[ht]
\centerline{\epsffile{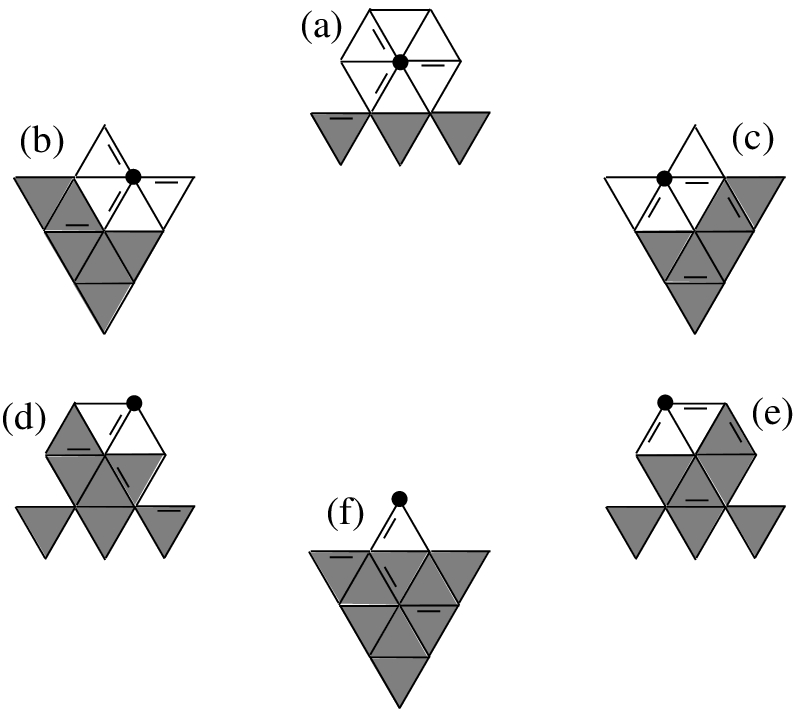}}
\caption{\label{downeligible}}
\end{figure}

We have illustrated the six configurations of eligible
$\theta$-translates that will be of interest to us in Figure
\ref{downeligible} for the case of $\lambda$ in the interior of the
alcove or on a left-hand diagonal wall. In both these cases the case
of $\eta$ non-$p^2$-regular gives the same result as the regular
case. In both cases the vertex $\theta$ is indicated by a dot. For the
alcove case the weight denoted $\mu$ in the preceding paragraph is in
the highest alcove of the cluster of alcoves drawn (which represents
the corresponding local composition pattern), and the unshaded alcoves
in each configuration are those containing eligible
$\theta$-translates.  For the wall case $\mu$ lies on the highest
marked wall (the markings represent the corresponding local
composition pattern), and the marked walls contained in the unshaded
region are those containing eligible $\theta$-translates. The
right-hand diagonal wall case is symmetric, while the eligible
$\theta$-translates in the horizontal wall case are illustrated with
bold walls in Figure \ref{upeligible} for $\eta$ $p^2$-regular. For
$\eta$ non-$p^2$-regular the only modifications occur in cases (a) (b) and
(c), which are illustrated in Figure \ref{nonregular} for $\mu$ a
$p^2$-vertex; the $p^2$-wall cases are similar.

Note that if $\nu$ is the lowest vertex in the closure of the alcove
containing $\lambda$, with a homomorphism from $\nabla(\nu)$ to
$\nabla(\theta)$, then the set of eligible $\theta$-translates in
$\nabla(\lambda)$ is a subset of $\Upsilon_{\lambda}$. The next
Theorem shows that the sets of eligible $\theta$-translates provide
the appropriate refinement of this set.

\begin{theorem}\label{downhomclass}
Suppose $\tau$ is a vertex and we have a composite of Carter-Payne
maps $\nabla(\tau)\to\nabla(\theta)$. If $\lambda$ lies in the closure
of the down alcove with lower vertex $\tau$, and $\mu\in W.\lambda$
lies in the neighbourhood of $\theta$ then
$$\Hom(\nabla(\lambda),\nabla(\mu))\neq 0$$ if and only if $\mu$ is an
eligible $\theta$-translate in $\nabla(\lambda)$.  Further, if there
is a non-zero homomorphism then it is either a composite of
Carter-Payne maps or a twisted map.
\end{theorem}
 
\noindent\pf 
We consider the case
where $\eta$ is $p^2$-regular; the remaining cases are easy
modifications which are left to the reader.
We will also only consider the case where $\lambda$ lies in the interior
of the down alcove; the remaining wall cases are similar, and are also
left as
an exercise for the reader.  Note that $\lambda$ and $\theta$ are such
that we are in one of the 6 configurations illustrated in Figure
\ref{downeligible} (a) to (f) (but it is not necessarily the case that
all alcoves occurring in the local composition pattern of $\mu$
correspond to composition factors of $\nabla(\lambda)$).  We begin by
showing that there exist suitable maps for each eligible
$\theta$-translate.
 
As $\lambda$ is in a down alcove, it is easy to verify that if $\mu$
is the lowest eligible $\theta$-translate (i.e. the one immediately
below $\theta$) then it is a composition factor of
$\nabla_p(\tilde{\lambda})$. Further, the existence of a map from
$\nabla(\tau)$ to $\nabla(\theta)$ immediately implies that there is a
map from $\nabla_p(\tilde{\lambda})$ to $\nabla_p(\mu)$ (as they both
come from taking the map from $\nabla(\tau'')$ to $\nabla(\theta'')$,
twisting by $F$, and tensoring both sides with an appropriate simple
module). Thus if $\mu$ is the lowest eligible $\theta$-translate then
there is a twisted map from $\nabla(\lambda)$ to $\nabla(\mu)$.
 
Next suppose that $\mu$ is the highest eligible $\theta$-translate. We
will show that the desired map can be constructed as a composite of
Carter-Payne maps using Corollary \ref{composite}. In this case $\mu$
is obtained from $\lambda$ by the same sequence of reflections used to
obtain $\theta$ from $\tau$. Write $\tau=(p^i-1)\rho +p^i\tau^{(i)}$
where $\tau^{(i)}$ does not lie on a vertex. Then the sequence of
reflections used must label a subcomposition of a composition in
max($\tau^{(i)}$) by Corollary \ref{composite} (though the maps will
have different superscripts). If the same subcomposition also occurs
as a subcomposition of a composition in max($\lambda$) then the first
condition in Corollary \ref{composite} is satisfied; by the argument
in the previous paragraph, the second condition is also satisfied, and
hence there is a non-zero composite of Carter-Payne maps from
$\nabla(\lambda)$ to $\nabla(\mu)$.

Thus it only remains to consider the cases where a subcomposition
$\phi$ of an element in max($\tau^{(i)}$) does not also occur as a
subcomposition of an element in max($\lambda$). By inspection we see
that the following cases need to be considered:

\begin{enumerate}
\item[(i)]{$\tau^{(i)}$ lies in a down alcove and
$\phi$ is one of $\phi_R^*\phi_B^1\phi_L^*$, $\phi_B^1\phi_L^*$,
 $\phi_L^*\phi_B^1\phi_R^*$, or  $\phi_B^1\phi_R^*$.}
\item[(ii)]{$\tau^{(i)}$ lies in an up alcove and $\phi=\phi_R^*\phi_B^1$ or
$\phi_L^*\phi_B^1$.}
\item[(iii)]{$\tau^{(i)}$ lies on a right-hand diagonal wall and
$\phi=\phi_L^*\phi_B^*$,  $\phi_L^*\phi_B^*\phi_R^*$, or  $\phi_B^*\phi_R^*$.}
\item[(iv)]{$\tau^{(i)}$ lies on a left-hand diagonal wall and
$\phi=\phi_R^*\phi_B^*$, $\phi_R^*\phi_B^*\phi_L^*$, or
$\phi_B^*\phi_L^*$.}
\end{enumerate} 
 
In each case we will show that there is a second non-zero composite of
Carter-Payne maps from $\tau$ to $\theta$ which does occur as a
subcomposition of an element of max($\lambda$). As we know that there is
(up to scalars) a unique map from $\nabla(\tau)$ to $\nabla(\theta)$
this new map must coincide with the original one, and the arguments of
the preceding paragraph will then apply to this new map to give the
desired composite map from $\lambda$ to $\mu$. In what follows all
equalities of maps should be interpreted as being up to some non-zero
scalar multiple.

Consider first case (i), with $\phi=\phi_R^a\phi_B^1\phi_L^b$ for some
$a,b>0$. By Lemma \ref{lr} and the uniqueness of maps we have that
$\phi_B^1\phi_L^b=\phi^b_L\phi^1_R$, and hence
$\phi=\phi_R^a\phi^b_L\phi^1_R$, which does occur as a subcomposition
of an element in max($\lambda$). (Note that as $\phi$ is non-zero it
follows that the conditions of Lemma \ref{lr} are satisfied.) Similar
arguments hold for the other maps in case (i). Case (ii) is similar
as, for example, we have that $\phi_R^a\phi_B^1=\phi_L^1\phi_R^a$.

Cases (iii) and (iv) are symmetric, so we only need consider (iii).
We would like to argue as in the preceding paragraph, but it is no
longer obvious that in applications of Lemma \ref{lr} the intermediate
weight is dominant. Consider the case where $\phi=\phi_L^a\phi_B^b$
(the other cases are similar). Consider the original map $\phi$ from
$\nabla(\tau)$ to $\nabla(\theta)$. By Corollary \ref{composite} there
must be a corresponding composite map from $\nabla_p(\tilde{\tau})$,
and hence from $\nabla(\tau'')$. Repeating this argument for
$\nabla(\tau'')$ (and so on), we deduce that the map must ultimately
have come from some corresponding composite with $a$ or $b$ equal to
$1$. By considering the possible cases (up alcove or RH diagonal wall)
where such a composite can arise, we see that we must have $b\leq a$,
and further that the map from $\lambda$ of the form $\phi_L^a$ must be
to a dominant weight. We have $\phi_R^b\phi_L^a=\phi_L^a\phi_B^b$ from
$\tau$, and hence $\phi=\phi_R^b\phi_L^a$. This is an element of
max$(\lambda)$, and hence the argument above produces the desired map.

The only remaining cases of eligible $\theta$-translates occur in
cases (a), (b) and (c) in Figure \ref{downeligible}. As above, the
sequence of reflections relating $\tau$ and $\theta$ must be a
subcomposition of some element in max($\tau$). By considering the
various possible cases we see that a configuration as in case (b)
(respectively (c)) can only occur when $\theta$ is obtained from
$\tau$ by a single reflection to the left (respectively right), while
for (a) we must have $\theta=\tau$. It is now easy to construct an
explicit composite of Carter-Payne maps for each of the remaining
weights $\mu$.

 

For $\lambda $ not too close to a $p^2 $-wall we have
shown that of the six alcoves in each cluster of possible candidates
for homomorphisms, either we have a morphism of the desired form or
the alcove contains no composition factor of $\nabla (\lambda) $.
Thus in that case we are done.  When $\lambda $ is close to a $p^2$-wall
it is possible that some of the six alcoves in a
cluster are generated as composition factors of $\nabla (\lambda) $ by
a weight on a vertex other than the one at the centre of the cluster.
As long as we can eliminate such alcoves as candidates for a
homomorphism the argument will proceed for the remaining alcoves just
as above. We will consider the case where $\lambda$ is close to a
single $p^2$-wall; the case where it is close to two such walls is
very similar.

\begin{figure}[ht]
\centerline{\epsffile{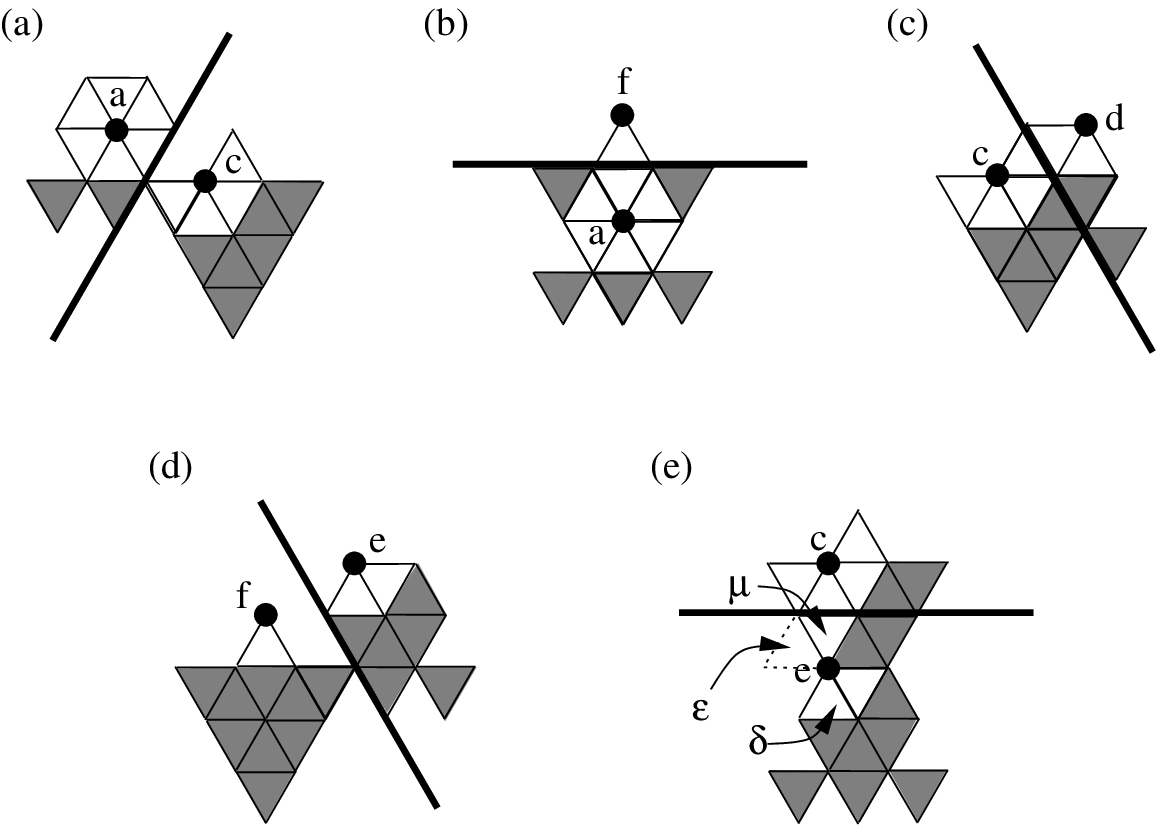}}
\caption{\label{nearcases}}
\end{figure}

\begin{figure}[ht]
\centerline{\epsffile{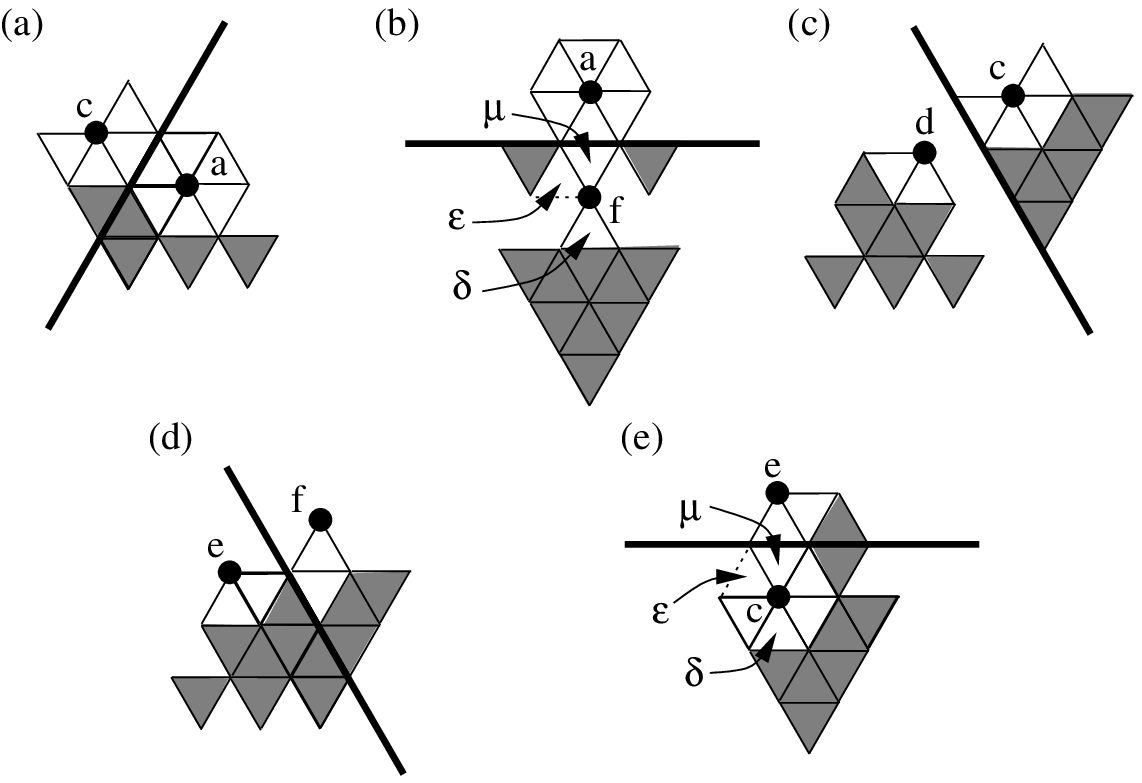}}
\caption{\label{nearcases2}}
\end{figure}

\begin{figure}[ht]
\centerline{\epsffile{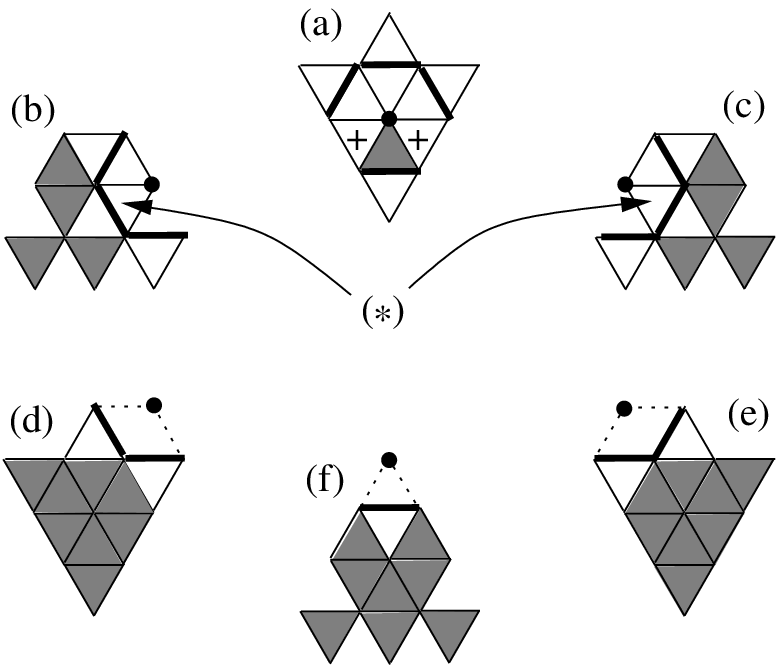}}
\caption{\label{upeligible}}
\end{figure}

\begin{figure}[ht]
\centerline{\epsffile{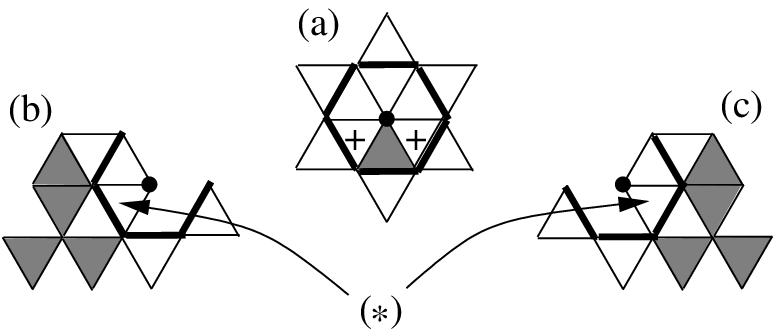}}
\caption{\label{nonregular}}
\end{figure}

If $\lambda $ lies close to a $p^2 $-wall, we see that the
problem case can occur if there exist vertices 
$\theta_1$ and $\theta_2$ separated by a single wall such that
$\theta_1>\theta_2$ and there is a non-zero homomorphism from $
\nabla(\tau)$ to $\nabla (\theta_2) $.  In this case it may be that
there is a local composition pattern clustered around $\theta_1$
formed from composition factors of $\nabla(\lambda)$, which intersects
the neighbourhood of $\theta_2$.  The various configurations that can
occur (after using symmetry considerations to reduce the number of
cases) are shown in Figures \ref{nearcases} and \ref{nearcases2}. Here
we label each of the local composition patterns by the corresponding
label in Figure \ref{downeligible} and shade in all but the 
$\theta_1$- and $\theta_2$-translates. Figure \ref{nearcases}
corresponding to $\tau''$ in a down alcove, and Figure
\ref{nearcases2} to $\tau''$ in an up alcove (the other cases are
similar).  By inspection, we see that the only cases where additional
homomorphisms from $\nabla(\lambda)$ to $\nabla(\mu)$ might arise are
when $\mu$ is as indicated in Figure \ref{nearcases}(e) or Figure
\ref{nearcases2}(b) or (e).

In these cases we need to show that there is no homomorphism
$\nabla(\lambda) \to \nabla(\mu)$: once we have done this we will have
completed the proof for $\lambda$ in the interior of a down alcove. So
suppose there is a homomorphism $\nabla(\lambda) \to \nabla(\mu)$ with
$\mu$ as shown. An easy calculation shows that the factor in the
$p$-filtration of $\nabla(\mu)$ in the same $\nabla_p$-class as
$\nabla_p (\tilde{\lambda})$ is $\nabla_p(\delta)$ where $\delta$ lies
in the alcove shown.  As $\nabla (\lambda) \to \nabla (\mu) $ is
non-zero, the head of $\nabla (\lambda) $ (i.e. the head of $\nabla_p
(\tilde{\lambda}) $) must survive in the image as a composition factor
of $\nabla_p (\delta) $.  But we have local homomorphisms $\nabla
(\mu) \to \nabla (\epsilon) $ (where $\epsilon$ lies in the dotted
alcove shown) where all of $\nabla_p (\delta) $ survives.  Therefore
we obtain a non-zero homomorphism $\nabla (\lambda) \to \nabla (\mu)
\to \nabla (\epsilon) $, which is impossible as $L (\epsilon) $ is not
a composition factor of $\nabla (\lambda) $.
\qed 



We next want to similarly refine $\Upsilon_{\lambda}$ for $\lambda$ in
an up alcove using sets of eligible $\theta$-translates.  

Let $\lambda$ lie in an up alcove and $\eta$ the vertex directly below
it.(That is, $\eta$ will be the centre of the star with $\lambda$ in the
topmost alcove of the star.) We wish to define the set of eligible
$\theta$-translates for $\nabla(\lambda)$; just as in the down alcove
case we begin by supposing that $\eta$ is $p^2$-regular.

Let $\theta$ be another vertex dominated by $\eta$ with $\theta = w
. \eta$ and $w \in W_{p^2}$. The {\it eligible $\theta$-translates for
$\nabla(\lambda)$} are those weights corresponding to composition
factors lying in the extended neighbourhood of $\theta$ which occur in
the local composition pattern for $\mu = w . \lambda$. When $\eta$ is
not $p^2$-regular the modifications are exactly as in the down alcove
case. 
 
As in the down alcove case there are six configurations of eligible
$\theta$-translates that will be of interest to us, and these are
illustrated in Figure \ref{upeligible} for $\eta$
$p^2$-regular. (These six configurations may be superimposed if we are
near a $p^2$-wall.)  As before the vertex $\theta$ is indicated by a
dot, the weight denoted $\mu$ in the preceding paragraph is in the
highest alcove of the cluster of alcoves drawn (which represents the
corresponding local composition pattern), and the unshaded alcoves in
each configuration are those containing eligible
$\theta$-translates. The thick lines indicate the eligible
$\theta$-translates for $\lambda$ in the horizontal wall case, and are
illustrated here for convenience later. For $\eta$ non-$p^2$-regular
the only modifications occur in cases (a) (b) and (c), which are
illustrated in Figure \ref{nonregular} for $\eta$ a $p^2$-vertex (the
$p^2$-wall cases being similar).

Unfortunately the up alcove analogue of Theorem \ref{downhomclass} is
a little more complicated than in the down alcove case. We define an
exceptional alcove for the pair $(\theta,\lambda)$ to be an alcove
($\ast$) as in Figures \ref{upeligible} or \ref{nonregular}, when the
relative positions of $\theta$ and the $\theta$-eligible translates in
$\nabla(\lambda)$ are as shown. We now have
 
\begin{theorem}\label{uphomclass}
Suppose $\tau$ is a vertex and we have a composite of Carter-Payne
maps $\nabla(\tau)\to\nabla(\hat{\theta})$. Let $\lambda$ lie in an up
alcove, with $\nu$ on the horizontal wall below $\lambda$ and $\tau$
the lower vertex in the closure of the down alcove below $\nu$. If we
have a weight $\mu\in W.\lambda$ with $\mu$ in the extended
neighbourhood of $\hat{\theta}$ then
$$\Hom(\nabla(\lambda),\nabla(\mu))\neq 0$$ if and only if there
exists some weight $\theta$ satisfying the following four conditions:
\begin{enumerate}
\item There exists a non-zero  composite of Carter-Payne
maps $\nabla(\tau)\to\nabla(\theta)$.
\item The weight $\mu$ is an eligible $\theta$-translate in
$\nabla(\lambda)$.  
\item \label{cond} There is a non-zero homomorphism from
$\nabla(\nu)$ to $\nabla(\delta)$ with $\delta$ on the wall in the
closure of the alcove containing $\mu$ which lies in the boundary of the
neighbourhood of $\theta$.
\item The weight $\mu$ is not an exceptional weight for the pair
$(\theta,\lambda)$.  
\end{enumerate} 
Further, if there is a non-zero homomorphism then it is either a
composite of Carter-Payne maps or a twisted map.  
\end{theorem}

\begin{remark} (i) As we will see in the following proof, we must have
$\theta=\hat{\theta}$ unless $\lambda$ is too close to a $p^2$-wall
(i.e. unless there are two vertices $\theta_1$ and $\theta_2$ with
overlapping extended neighbourhoods such that $\theta_1$ and
$\theta_2$ are both in $W_{p^2}.\tau$).\\
(ii) Condition (3) eliminates the possibility of maps into the alcoves
marked $(+)$ in Figure \ref{upeligible} in the generic case.
\end{remark}
 
\noindent\pf We proceed as in the proof of the last Theorem,
and as in that proof we will only consider the case
where $\eta$ is $p^2$-regular; the remaining cases are easy
modifications which are left to the reader.
We  begin
by assuming that $\lambda$ is not close to a $p^2$-wall, so that
$\lambda$ and $\theta$ are such that we are in one of the six
configurations illustrated in Figure \ref{upeligible} (a) to (f).  
We now consider the eligible $\theta$-translates.

As in the preceding theorem, it is easy to verify that if $\mu$ is the
lowest eligible $\theta$-translate (i.e. the down alcove below
$\theta$) then it is a composition factor of
$\nabla_p(\tilde{\lambda})$. Just as before, we deduce that there
exists a twisted map from $\nabla(\lambda)$ to $\nabla(\mu)$.

Similarly, if $\mu$ is the highest eligible $\theta$-translate, we
wish to argue exactly as in the preceding proof to show that there is
a composite of Carter-Payne maps from $\nabla(\lambda)$ to
$\nabla(\mu)$. As there we see that the only problem occurs when the
composite of Carter-Payne maps from $\tau$ to $\theta$ does not
correspond to a subcomposition of an element of max$(\lambda)$. With
the notation used in the preceding proof, it remains to consider the
following cases:

\begin{enumerate}
\item[(i)]{$\tau^{(i)}$ lies in a down alcove and
$\phi=\phi_R^*\phi_B^1\phi_L^*$ or $\phi_L^*\phi_B^1\phi_R^*$.}
\item[(ii)]{$\tau^{(i)}$ lies in an up alcove and $\phi=\phi_R^*\phi_B^1$ or
$\phi_L^*\phi_B^1$.}
\item[(iii)]{$\tau^{(i)}$ lies on a right-hand diagonal wall and
$\phi=\phi_L^*\phi_B^*$  or  $\phi^*_L\phi_B^*\phi_R^*$.}
\item[(iv)]{$\tau^{(i)}$ lies on a left-hand diagonal wall and
$\phi=\phi_R^*\phi_B^*$   or  $\phi^*_R\phi_B^*\phi_L^*$.}
\end{enumerate}
 
As before, in each case we will construct a second non-zero composite
of Carter-Payne maps from $\tau$ to $\theta$ which does correspond to
a subcomposite of an element from max($\lambda$). The desired maps
then exist by the arguments in the preceding proof. 

Case (i) is straightforward, as there is only one eligible
$\theta$-translate. Thus the highest eligible $\theta$-translate is
also the lowest, and we have already constructed the required map.
Case (ii) is also easy, as we have that
$\phi_R^a\phi_B^1=\phi_L^1\phi_R^a$ (and similarly with the roles of
left and right reversed), and it is easy to see that this composite
satisfies the conditions of Lemma \ref{lr}.

As the right-hand and left-hand wall cases are symmetric, we only need
consider case (iii). If $\phi=\phi_L^a\phi_B^1$ we argue as in case
(ii). If $\phi=\phi_L^a\phi_B^b\phi_R^c$ we have to be careful (as in
the proof of case (iii) in Theorem \ref{downhomclass}) that in
applying Lemma \ref{lr} we have that the intermediate weight is
dominant. However, as before, we can deduce that the map must
ultimately come from a map from a right-hand diagonal or internal down
alcove, and just as in the corresponding case in Theorem
\ref{downhomclass} we deduce that $\phi_B^b\phi_R^c=\phi_R^c\phi_L^b$,
and hence that $\phi=\phi_L^a\phi_R^c\phi_L^b$ which corresponds to an
element in max($\lambda$) as required.
 
The only remaining cases of eligible $\theta$-translates occur in
cases (a), (b) and (c) in Figure \ref{upeligible}. As in the preceding
proof, we see that a configuration as in case (b) (respectively (c))
can only occur when $\theta$ is obtained from $\tau$ by a single
reflection to the left (respectively right), while for (a) we must
have $\theta=\tau$.  
 
In case (a) we already know that there exist maps of the desired form
into all alcoves except those labelled ($+$). For each of these
latter alcoves no map exists by translation arguments, as there is no
map from $\nabla(\nu)$ to $\nabla(\delta)$, with $\delta$ as in
condition \ref{cond}. (This case is the reason
for this condition in the Theorem, which eliminates these two
possibilities.) 
 
For cases (b) and (c) it is now easy to construct an explicit
composite of Carter-Payne maps for each of the remaining weights
$\mu$, except for the exceptional weights labelled ($\ast$). As these
two cases are symmetric we consider case (c), and must show that there
is no non-zero map from $\nabla(\lambda)$ to $\nabla(\mu)$ when $\mu$
is the weight ($\ast$). Consider the $p$-filtration for $\nabla(\mu)$,
using the labelling of factors given in Figure \ref{alcoves}a(i). If
there were a non-zero map then the head of $\nabla(\lambda)$ must
survive, and by earlier arguments (and our assumption that $\lambda$
is not near a $p^2$-wall) it must be a composition factor in
$\nabla_p(\mu_8)$. As there is a map from $\nabla(\mu)$ to
$\nabla_p(\mu_4)$ which does not kill any of $\nabla_p(\mu_8)$, the
composite of this pair of maps would give a non-zero map from
$\nabla(\lambda)$  to $\nabla(\mu_4)$. But $L(\mu_4)$ is not a
composition factor of $\nabla(\lambda)$, which gives the desired
contradiction. 

By our assumption on $\lambda$, the remaining weights in the extended
neighbourhood of $\theta$ for which we have not yet constructed maps
are not composition factors of $\nabla(\lambda)$. Thus we are done in
the case where $\lambda$ is not too close to a $p^2$-wall.

\begin{figure}[ht]
\centerline{\epsffile{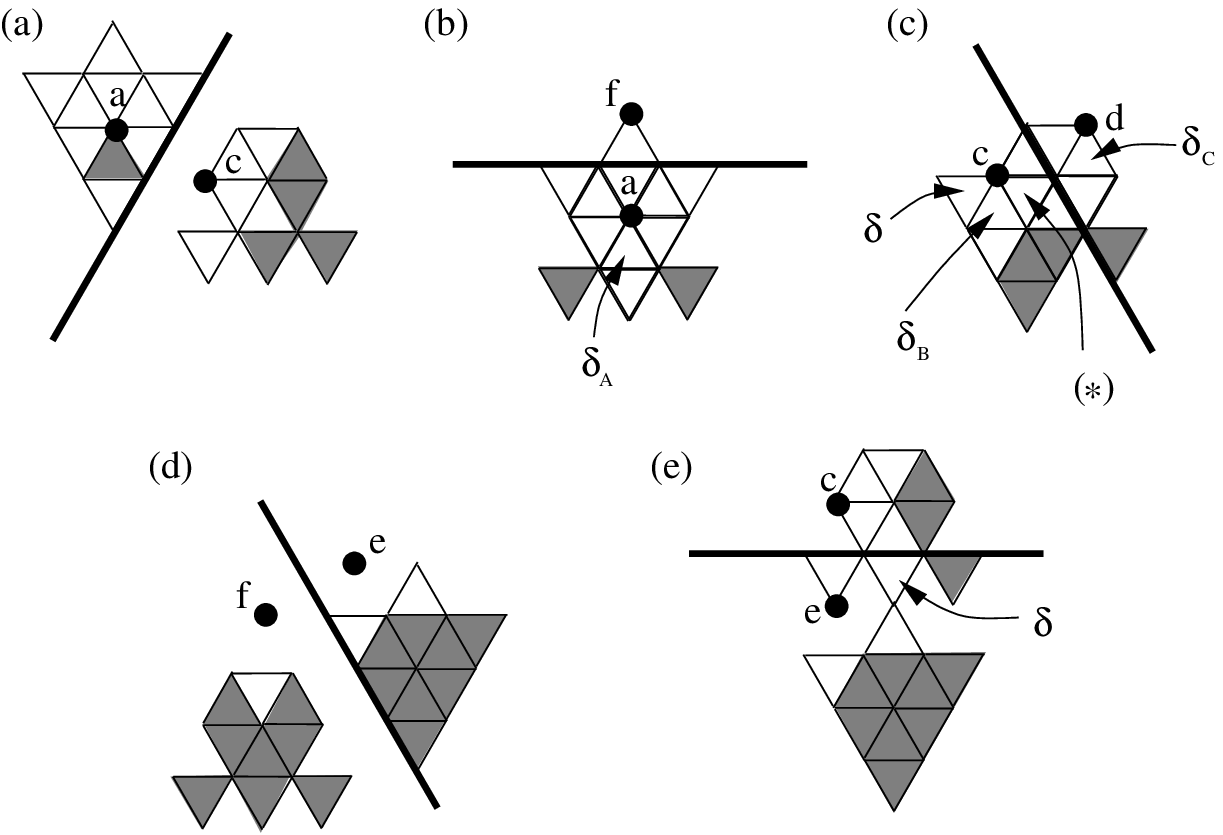}}
\caption{\label{upnear}}
\end{figure}

\begin{figure}[ht]
\centerline{\epsffile{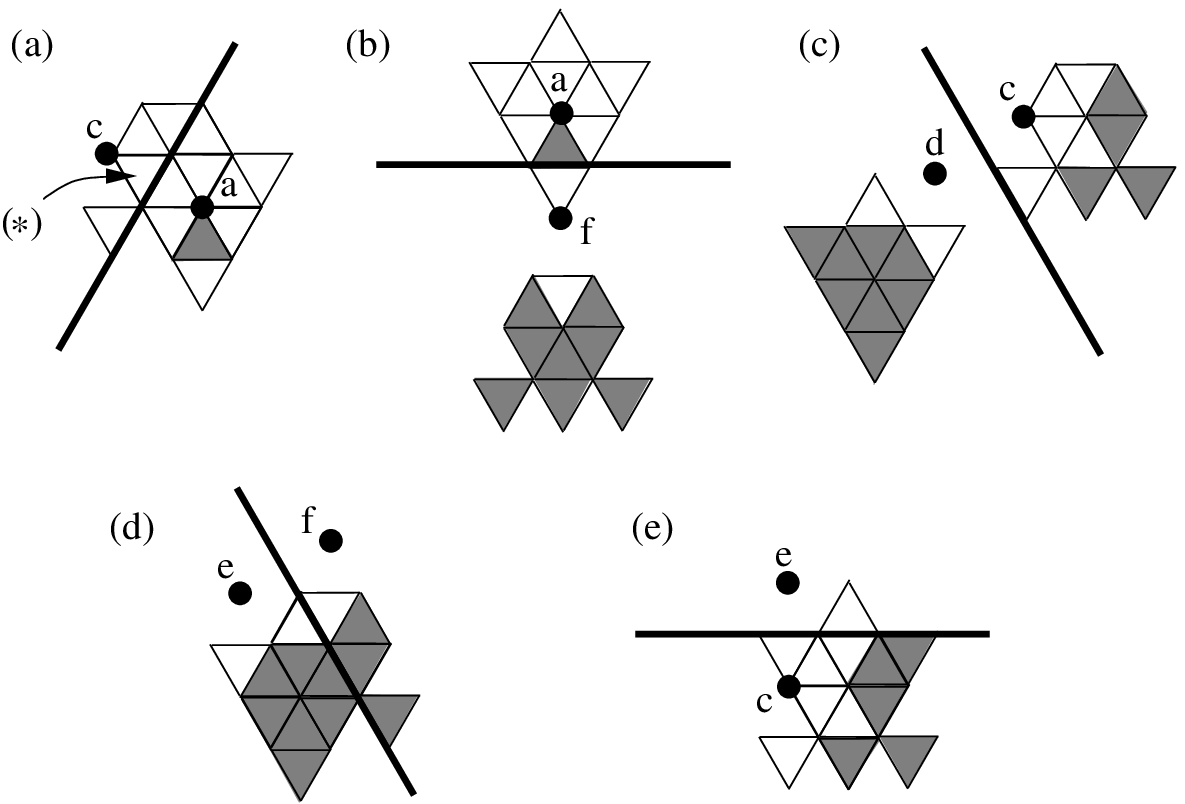}}
\caption{\label{upnear2}}
\end{figure}

It remains to consider the case where $\lambda$ is close to a
$p^2$-wall. If $\lambda$ is close to only one $p^2$-wall then the new
configurations that can occur are illustrated in Figure \ref{upnear}
(for $\tau''$ in a down alcove, the wall cases are similar) and
Figure \ref{upnear2} (for $\tau''$ in an up alcove). Here the
vertices are labelled by the labels of the corresponding
configurations in Figure \ref{upeligible}, and all alcoves not in the
extended neighbourhood of one of the two marked vertices are
shaded. When $\lambda$ is close to two $p^2$-walls it is easy to
verify that the obvious modification of the arguments below
(considering each $p^2$-wall case separately, then superimposing them)
holds.

The construction of maps proceeds just as in the case
considered above, with the understanding that there may be two \lq
highest\rq\ eligible $\theta$-translates: one for each local
composition diagram when a pair of such overlap. However, it remains
to show that there are no additional maps into the remaining unshaded
alcoves in Figures \ref{upnear} and \ref{upnear2}.

If there is a map from $\nabla(\nu)$ to $\nabla(\delta)$ in Figure
\ref{upeligible} with $\delta$ and $\nu$ as in Condition \ref{cond} 
then the alcove marked with a $+$ must coincide with an exceptional
alcove (marked $(*)$) and so this map is eliminated as in cases (b)
and (c). 
For the alcove
marked $(\ast)$ in Figure \ref{upeligible} the argument given above
only breaks down in the case shown in Figures \ref{upnear}(c) and
\ref{upnear2}(a). Consider first the former case. The argument given
above breaks down because $\mu_4$ (here labelled by $\delta_B$) is a
composition factor of $\nabla(\lambda)$. However, it will be enough to
show that there is no map from $\nabla(\lambda)$ to
$\nabla(\delta_B)$, which is proved below.

\begin{figure}[ht]
\centerline{\epsffile{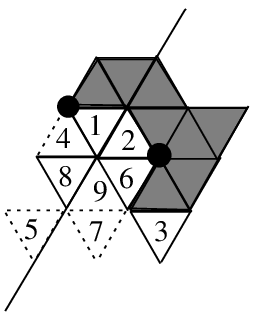}}
\caption{\label{explain}}
\end{figure}

Next consider the case shown in Figure \ref{upnear2}(a). For
convenience we have redrawn this case in Figure \ref{explain} with the
labelling of the $p$-filtration for $\nabla(\mu)$ shown and all other
alcoves shaded. Here we may no longer deduce that the head of
$\nabla(\lambda)$ must be a composition factor of $\nabla_p(\mu_8)$,
as it could also occur in $\nabla_p(\mu_3)$. Suppose that this is the
case. Then our map must induce a morphism from $\nabla_p(\lambda_3)$
to $\nabla_p(\mu_3)$, and hence there exists a map from
$\nabla(\lambda_3'')$ to $\nabla(\mu_3'')$. But these two weights
differ by a multiple of a single root and are not reflections of each
other about a wall. Thus this map cannot exist by \cite[(3.6)
Theorem and (4.3) Corollary]{erd}. Thus we see that the
head of $\nabla(\lambda)$ must be a composition factor of
$\nabla_p(\mu_8)$, and the argument given in the generic case goes
through unchanged.

It remains to show that the unshaded alcoves in Figures \ref{upnear}
and \ref{upnear2} for which we have not yet constructed maps do indeed
have no maps into them. First suppose that  there are maps from
$\nabla(\tau)$ into {\it both} labelled vertices in a given
configuration. Then the remaining cases to consider are the alcoves
labelled by Greek letters in Figures \ref{upnear}(b), (c) and (e).

Those labelled simply by a $\delta$ clearly have no maps, as they fail
condition (3) of the Theorem, which is a necessary condition for a map
to exist by the translation arguments in the proof of Theorem
\ref{onedim}. Next consider $\mu=\delta_A$, $\delta_B$, or $\delta_C$ and
suppose there is a non-zero map from $\nabla(\lambda)$ to
$\nabla(\mu)$. In each case $\mu$ must be a composition factor of
$\nabla_p(\lambda_4)$. Hence the head of $\nabla_p(\lambda_9)$ must
survive under the map, as the submodule of $\nabla(\lambda)$
with $\nabla_p$-factors $\nabla_p(\lambda_4)$, $\nabla_p(\lambda_5)$,
$\nabla_p(\lambda_6)$, $\nabla_p(\lambda_7)$, $\nabla_p(\lambda_8)$
and $\nabla_p(\lambda_9)$ has simple head. But there is no composition
factor of $\nabla_p(\lambda_9)$ in $\nabla(\mu)$ as $\nabla(\mu)$ does
not contain this $\nabla_p$-class.

%

Finally, we have to consider the case where there is a map from $\tau$
to just {\it one} of the two labelled vertices in a given configuration. As
for the cases labelled by $\delta$'s above, an appeal to the necessary
condition (3) in the Theorem is enough to eliminate the remaining
alcoves in most cases.

However, two cases cannot be eliminated in this manner. These are the
configuration shown in Figure \ref{upnear}(b) when there is a map from
$\tau$ only to vertex $f$, and the configuration shown in Figure
\ref{upnear}(c) when there is map from $\tau$ only to vertex $d$. In
each case we must eliminate the possibility of a map to each of the
alcoves adjacent to vertex $f$ (respectively $d$).

In case (b) it is easy to verify that the alcove is in the same
$\nabla_p$-class as $\lambda_4$. Then arguing as above we see that the
head of $\nabla_p(\lambda_9)$ must survive in the image of any map,
and thus the composite of this map followed by a local reflection to
the right would be non-zero. However, no such map exists by our
assumption.

\begin{figure}[ht]
\centerline{\epsffile{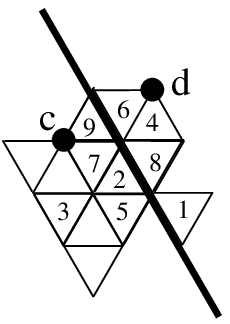}}
\caption{\label{onevert}}
\end{figure}

In case (c) the lower of the two alcoves has already been eliminated
above. The remaining alcove is marked $6$ in Figure \ref{onevert}, where
the numbering indicates the $\nabla_p$-class using the labelling for
$\lambda$. The structure of the $p$-filtrations of
$\nabla(\lambda)$ and $\nabla(\mu)$ are very similar to those shown on
the left (respectively right) sides of Figure \ref{Second} --- except
that in the latter case there is an extension of $3$ by $8$, as these are
in the same $\nabla_p$-class. Any map from $\nabla(\lambda)$ must
preserve the head of $\nabla_p(\lambda_3)$, and its image must lie in
the term in the $p$-filtration of $\nabla(\mu)$ corresponding to the
alcove marked $3$ or $8$.

If this composition factor lies in $\nabla_p(\mu[\lambda_8])$ then
composing the map with the local reflection to the right gives a
non-zero map to $\nabla(\mu[\lambda_4])$, which contradicts the
arguments above. Hence the image of $\nabla(\lambda)$ under our map
must be in the submodule of $\nabla(\mu)$ with $p$-filtration labelled
by $\mu[\lambda_3]$, $\mu[\lambda_6]$, and $\mu[\lambda_9]$. But then
the composite of our map with the obvious quotient map to
$\nabla_p(\mu[\lambda_3])$ must be non-zero, which implies the
existence of a non-zero map from $\nabla(\lambda)$ to
$\nabla(\mu[\lambda_3])$. However, no such map exists by our
assumption.
\qed 

For twisted maps, the analogue of Theorems \ref{downhomclass} and
\ref{uphomclass} is in most cases much more straightforward.

\begin{theorem}\label{twistclass}
Suppose $\tau$ is a vertex and we have a twisted map
$\nabla(\tau)\to\nabla(\theta)$. Let $\lambda$ lie in the closure of
the down alcove with lower vertex $\tau$, or in an up alcove with
$\nu$ on the horizontal wall below $\lambda$ and $\tau$ the lower
vertex in the closure of the down alcove below $\nu$.  Write
$\tau=(p^i-1)\rho+\tau^+$ where $\tau^+$ is not a vertex weight. If
$\widetilde{\tau^+} = p((\tau^+)''-1) \rho +w_a ((\tau^+)')$ then
$$\Hom(\nabla(\lambda),\nabla(\mu))\neq 0$$ if and only if $\mu$ is an
eligible $\theta$-translate in $\nabla(\lambda)$.  Further, if there
is a non-zero homomorphism then it is twisted map.
\end{theorem}
\pf By our assumption on $\tau$, the existence of a twisted map from
$\nabla(\tau)$ to $\nabla(\theta)$ implies that the $\theta$-eligible
factors of $\nabla(\lambda)$ are in one of the configurations (f) in
Figures \ref{downeligible} and \ref{upeligible}. As in both these
cases there is only one eligible factor $\mu$ (which is thus the
lowest such) the proof of the existence of a (necessarily twisted) map
follows just as in the corresponding case in Theorems
\ref{downhomclass} and \ref{uphomclass}.\qed

The assumption in the preceding Theorem that $\tau $ can be written in
the form $\tau=(p^i-1)\rho+\tau^+$ where $\tau^+$ is not a vertex
weight only fails to hold when $\tau^+$ is close to the boundary of
the dominant region. For such weights a twisted map from $\nabla(\tau)$
to $\nabla(\theta)$ will give rise to the first of our exceptional
maps, as well as the twisted maps in Theorem \ref{twistclass}. We
return to this problem in Section \ref{exception}, after we have
considered those weights $\lambda$ not of the form considered in
Theorems \ref{downhomclass}--\ref{twistclass}.




\section{Determining homomorphisms II: just dominant weights}\label{symway}

We next turn our attention to weights $\lambda$ for which the methods
of the previous section do not apply.  We first consider the case
where $\lambda$ is in the closure of a down alcove whose lower vertex
is non-dominant.  By the translation principle, we may assume that
$\lambda = (a, 0) $ or $(0,b) $.  It is easy to show (cf. \cite[II,
2.16]{jantzen}) that $\nabla (a, 0) \cong S^a (V) $, where $V$ is the
natural module for $\SL_3(k)$, and the submodule structure of these
has been determined by Doty \cite{doty}.  We will use (and refine) the
reinterpretation of these results in \cite{cox4}.

Let $a = \sum_{i = 0}^{m} a_ip^i$ where $0 \leq a_i \leq p-1$ for all
$i $ and $a_m \neq 0 $.  Doty gives a procedure for defining a set
$C(a)$ of $m $-tuples of nonnegative integers (called carry patterns)
which is in one-to-one correspondence with the set of composition
factors of $S^a(V)$.  Further, if we impose a poset structure on
$C(a)$ by setting ${\bf c} \leq{\bf d} $ if $c_i \leq d_i$ for all $i
$ then this correspondence induces a lattice isomorphism from the
lattice of order closed subsets of $C(a)$ to the lattice of submodules
of $S^d(V)$. Doty's construction is purely arithmetic; we will use the
alcove theoretic version given in \cite{cox4}.

\begin{figure}[ht]
\centerline{\epsffile{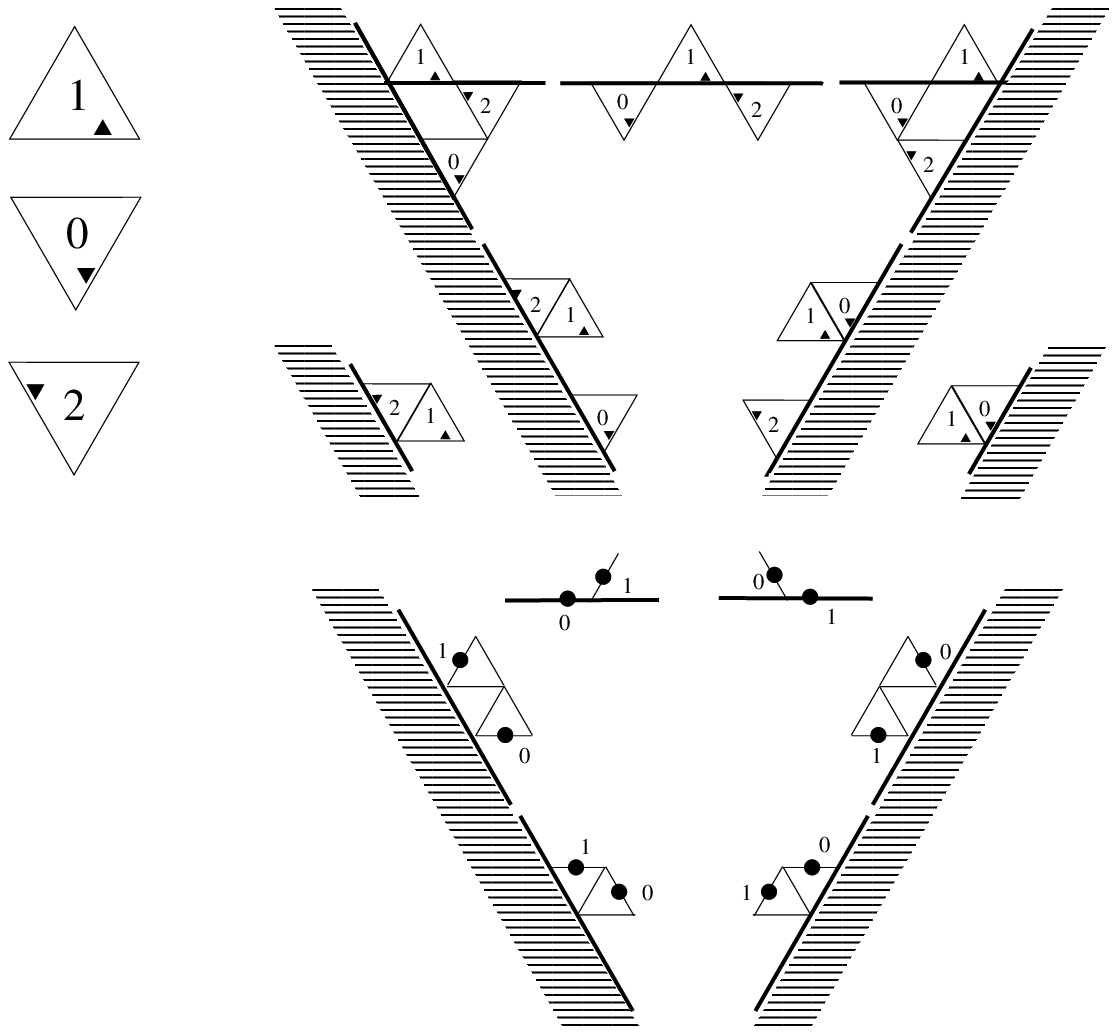}}
\caption{\label{sympower}}
\end{figure} 

Recall that in \cite[Figure 5]{cox4} each type of alcove that can
contain a composition factor of $S^a(V)$ was associated to either $0$,
$1$ or $2$. We will extend this procedure to walls:
we reproduce the
correspondence between alcoves and integers in Figure \ref{sympower},
and extend it to the wall cases as shown. Note that in the new cases
the numbering depends on the type of wall containing $\lambda$.

We now recall the recursive procedure for determining composition
factors of $\nabla(\lambda)$ given in \cite{cox4}, in the special case
where $\nabla(\lambda)$ is a symmetric power. (This allows certain
simplifications to be made.) As this is a little complicated, an
example will be given below. The algorithm starts by finding the
largest $i$ such that $\lambda$ does not lie in the closure of the
lowest $p^i$-alcove. By regarding the facets shown in Figure
\ref{sympower} as $p^i$-facets, we let cf$(\lambda,i)$ be the set of
weights obtained from $\lambda$ by intersecting $W_{p^i}.\lambda$ with
the configuration of $p^i$-facets corresponding to the position of
$\lambda$. (The weight $\lambda$ will always occur in this set.) In
this first iteration of the algorithm the configuration containing
$\lambda$ will be one with a shaded region to the right of it.

For the inductive step, we construct a new set cf$(\lambda,i-1)$ from
cf$(\lambda,i)$. For each weight $\mu$ in cf$(\lambda,i)$ we consider
the unique translate of the set of $p^{i-1}$-restricted weights which
contains $\mu$ and determine the configuration of $p^{i-1}$-facets
given in Figure \ref{sympower} which corresponds to the position of
$\mu$ {\it in this set}. (Thus the shaded regions here represent
weights outside of the translate considered.) Then the intersection of
$W_{p^{i-1}}.\mu$ with these $p^{i-1}$-facets gives the (immediate)
descendants of $\mu$ in cf$(\lambda,i-1)$. All elements of
cf$(\lambda,i-1)$ arise in this manner. In \cite{cox4} it was proved
that cf$(\lambda,1)$ is precisely the set of composition factors of
$\nabla(\lambda)$.

As $S^a(V)$ has no repeated composition factors, each weight $\mu$
that labels such a factor has a unique collection of ancestors; i.e. a
series $\mu^i=(a,0), \mu^{i-1}\ldots,\mu^1=\mu$ such that $\mu^j$
occurs as an immediate decedent of $\mu^{j+1}$ in
cf$(\lambda,j)$. (Note that not all of these weights are necessarily
distinct.) Each such $\mu^j$ comes from a configuration of
facets given in Figure \ref{sympower}, and hence has an associated
integer $a_j(\mu)$ from $\{0,1,2\}$. This is the series of integers
which we assign to $\mu$. We can now prove the following refinement of
\cite[Proposition 5.4]{cox4}.

\begin{proposition}  Given a composition factor $L(\mu)$ of 
$S^a(V)$, the associated carry pattern is the $m$-tuple of integers
whose $j$th entry equals $a_j(\mu)$
\end{proposition}  
\pf Argue as in the proof of \cite[Proposition 5.4]{cox4}, noting that
in the extra cases (where $(a, 0) $ lies on a wall), the $p^i $-facets
containing $(a, 0) $ are all the same type until the first such
$p^j$-facet which is a alcove, after which the remaining $p^i $-facets
(for $i \geq j $) are also all of the same type.\qed

\begin{figure}[ht]
\centerline{\epsffile{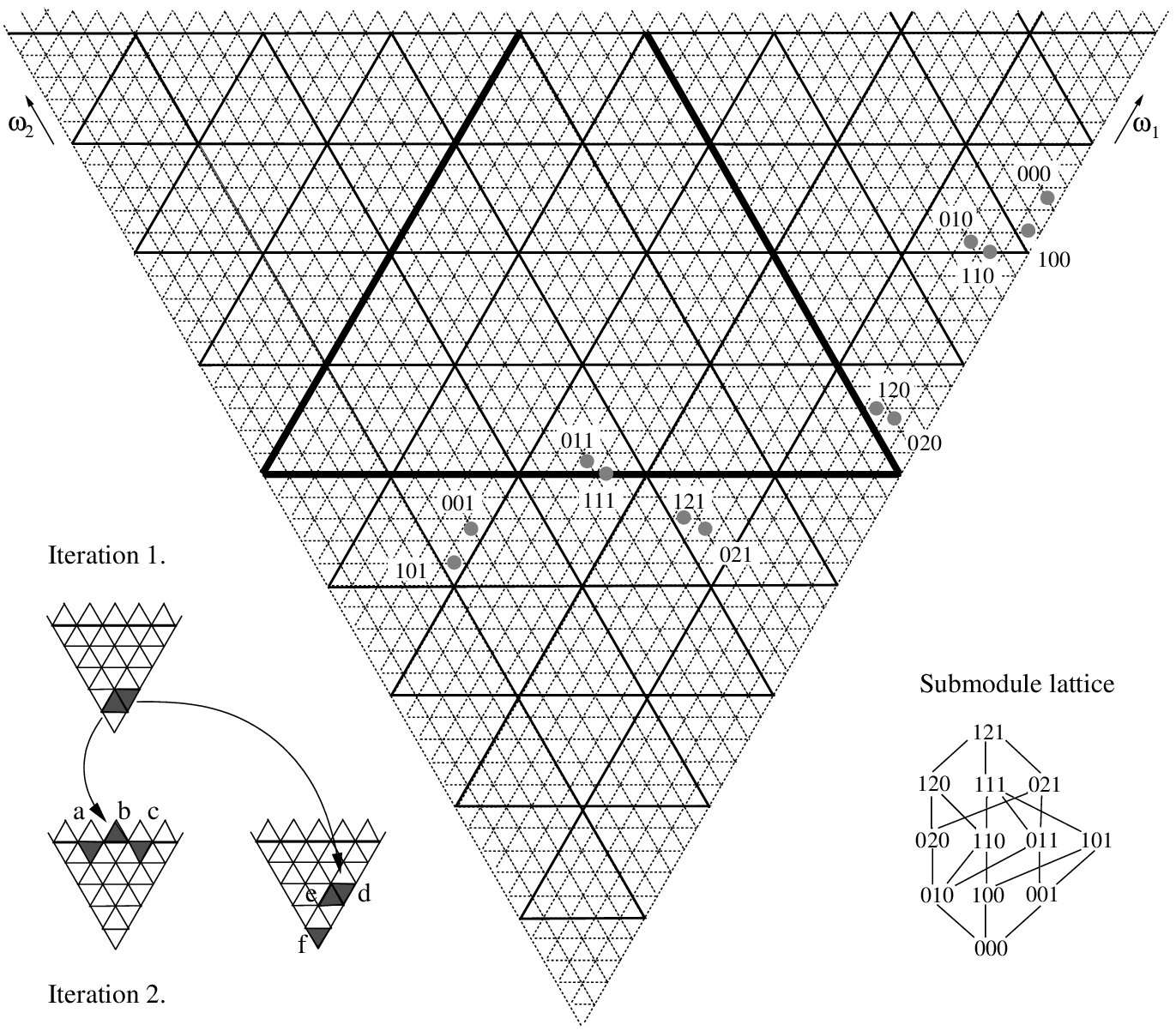}}
\caption{\label{symexample}}
\end{figure}

As an example of this Proposition consider $S^{184}(V)$ wit $p=5$. The
composition factors and their associated carry patterns are
illustrated in Figure \ref{symexample}, together with the first two
iterations of the composition factor algorithm and the complete
submodule lattice. 

At the first iteration of the algorithm, the configuration of alcoves
used is as shown in Iteration 1, and hence the final digit of each
carry pattern is either $0$ or $1$ by Figure \ref{sympower}. Using
Iteration 2 we can determine the penultimate digit of each factor; for
example those arising from Iteration 2(c) will have final pair
$21$. For the final iteration we use the wall cases given in Figure
\ref{sympower} to give the initial digit of each carry pattern.

Another, slightly more  detailed, example of the mechanics of the
composition factor algorithm (in the alcove case, with
$\lambda=(181,0)$, again with $p=5$) is given in \cite[Figure 2]{cox4}. The 
corresponding submodule lattice is given in \cite[Figure
6]{cox4}. (Note that there is an edge missing in the latter diagram
from the node 120 to the node 121.)

\begin{remark} Note that the first digits of the carry patterns
distinguish the terms in the $p $-filtration of $S^a (V) $ in which
the corresponding composition factors lie. Similarly for a fixed
choice of the first $i$ digits, the $(i+1)$st digits distinguish the terms 
in the $p^{i+1}$-filtration of an appropriate $\nabla_{p^i}(\mu)$.
\end{remark} 

We are now in a position to prove

\begin{theorem} \label{symprop} 
Suppose that $\lambda=(a,0)$ or $(0,b)$. All non-zero homomorphisms
$\nabla(\lambda) \to \nabla(\mu)$ are either composites of
Carter-Payne maps (as described in the preceding section) or induced
from non-zero homomorphisms $\nabla_p (\tilde{\lambda}) \to \nabla_p
(\mu) $.  

In particular, for $\lambda$ of the form $(a,0)$ all
non-twisted maps are of the form $\phi^*_L$ or $\phi^*_R\phi^*_L$, and
twisted maps correspond (after translation and twisting) to
homomorphisms from some $\nabla(c,0)$ where $(c,0)$ is any weight
lying in the same facet as $(\tilde{\lambda})''$. The case
$\lambda=(0,b)$ is similar.
\end{theorem}
\pf As the contravariant duals of $\nabla (a ,b )^*$ and $L (a ,b )^*$
are $\nabla (b ,a) $ and $L(b,a) $ respectively, we may assume that
$\lambda = (a, 0) $.  For there to be a non-zero homomorphism it is
clearly necessary that the quotient of $\nabla (\lambda) $ with simple
socle $L( \mu) $ must only contain composition factors $L (\tau) $ with
$\tau \leq \mu $.  (As $\nabla ( \lambda) $ is a symmetric power, and
hence is multiplicity free, this quotient is well defined.)  By the
universal property of $\nabla(\mu)$ this condition is also sufficient;
we will determine precisely the set of such $\mu$ (and show in doing
so that the maps are of the form described in the statement of the
Theorem). 

Suppose that $\lambda $ lies in the interior of an alcove.  For $\mu $
to be a possible candidate for a homomorphism, the series of steps in
the composition factor algorithm which leads to $\mu $ must be
such that, at every stage, the immediate descendants of $\mu^{j+1}$
occurring above $\mu^j$ in the relevant configuration from Figure
\ref{sympower} must be labelled by smaller integers in
Figure \ref{sympower}.  (Otherwise it is easy to verify that there
exists a weight $\tau\geq\mu $ whose carry pattern is above that of
$\mu $.) 
By considering the various
configurations given in Figure \ref{sympower}, we see that the only
candidates for homomorphisms are those weights for which maps are
already known to exist by the results in the preceding section.  That
completes the proof when $\lambda $ lies in the interior of an alcove;
the two wall cases are similar.\qed

The reader may find it helpful to consider the example in Figure
\ref{symexample}. In this case there are homomorphisms into the
induced modules labelled by 000 (identity), 011 ($\phi_L^3$), 121
($\phi_R^2\phi_L^3$), 111 ($\phi_R^1\phi_L^3$), 010 ($\phi_L^2$), 110
($\phi_R^1\phi_L^2$), 100 (twist of identity map from $100''$), 110
(twist of $\phi_L^1$ from $100''$), and 120 (twist of twist of
identity from $(120'')''$ in $100''$). (Here we are abusing notation
by combining carry patterns with our conventions on weights with
primes on.) Note that these are precisely the vertices $v$ in the
submodule lattice such that all vertices above them in the lattice
label weights occurring below the weight labelled by $v$.

The only class of facets which we have not yet considered are just
dominant up alcoves. The argument here is an easier form of Theorem
\ref{uphomclass}. Although there is no longer a vertex below the
weight $\lambda$ under consideration, we will use our knowledge of
homomorphisms from a weight $\nu$ on the horizontal wall below it
instead (which was already required as condition (3) in the statement
of Theorem \ref{uphomclass}). By the translation arguments in the
proof of Theorem \ref{onedim} we know that a necessary condition for
the existence of a map from $\nabla(\lambda)$ to $\nabla(\mu)$ is that
there is a map from $\nabla(\nu)$ to $\nabla(\delta)$ for some $\delta$
on one of the walls in the closure of the alcove containing $\mu$. 

\begin{figure}[ht]
\centerline{\epsffile{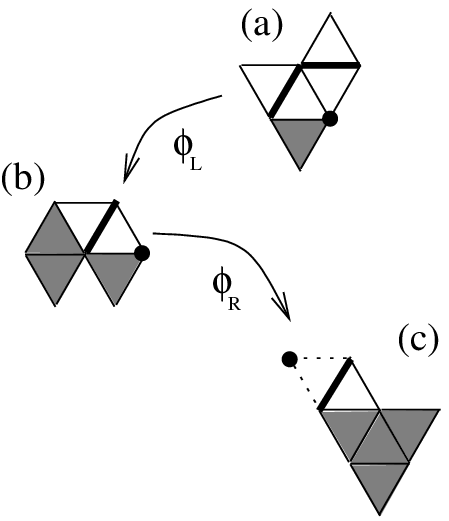}}
\caption{\label{justup}}
\end{figure}

\begin{theorem}\label{justuphoms}
Suppose that $\nu$ lies on a just dominant horizontal wall, and we
have a composite of Carter-Payne maps from $\nabla(\nu)$ to
$\nabla(\delta)$. Let $\lambda$ lie in an up alcove with $\nu$ on the
wall below $\lambda$ and $\tau$ the lower vertex in the closure of the
down alcove below $\nu$. If we have a weight $\mu\in W.\lambda$ with
$\mu$ in either of the alcoves adjacent to the wall containing
$\delta$ then
$$\Hom(\nabla(\lambda),\nabla(\mu))\neq 0.$$ 
When $\lambda$ is a right-hand just dominant weight the set of such
maps is precisely those maps in
$$\{\phi_R^*\phi_L^*, \phi_L^*, \phi_B^1, \phi_L^1\phi_B^1\}$$
which only involve reflections between dominant weights. The case of
$\lambda$ a left-hand just dominant weight is similar.
\end{theorem}
\pf We consider the right-hand just dominant case; the other is
similar. By Theorem \ref{symprop} the map from $\nabla(\nu)$ to
$\nabla(\delta)$ must be of the form $\phi_L^*$ or $\phi_R^*\phi_L^*$,
and hence the various configurations that must be considered are as
shown in Figure \ref{justup}.

By Lemma \ref{lr} and the local data in Section \ref{localdata} we see
that each of the composites given in the statement of the Theorem is
non-zero, and by considering Figure \ref{justup} we see that no
further maps are required.
\qed

In Theorem \ref{twistclass} we constructed a second class of maps for
interior weights, coming from twisted maps. Thus we also need to
consider the corresponding case for $\lambda$ in a just dominant up
alcove. The untwisted form of the top $\nabla_p$-factor in
$\nabla(\lambda)$ will be a symmetric power. For any map from such a
symmetric power which is a composite of Carter-Payne maps, it is easy
to verify using Figure \ref{justup} and the classification of such
maps in Theorem \ref{symprop} that the twisted version is actually one
of the Carter-Payne maps of the original module considered in the
Theorem above. Thus the only remaining problems are the cases arising
from twists of maps which are not Carter-Payne composites. These are
the subject of the following section.



\section{Determining homomorphisms III: exceptional maps}\label{exception}

Let us review the cases which have been covered so far. If $\lambda$
is just dominant then all possible homomorphisms have been
constructed in Section \ref{symway}.

If $\lambda$ is a vertex then $\nabla(\lambda)=\nabla_p(\lambda)$ and
the result is clear from (\ref{homreduce}).  For all remaining
$\lambda$ we have used translation arguments to classify all
homomorphisms induced from vertex homomorphisms which are either
composites of Carter-Payne maps or twisted maps for certain
weights. (Note that this includes all weights in the lowest
$p^2$-alcove, where all translations are of the identity map.)

There remains one extra class of possible maps, which we call {\it
exceptional maps}. These maps are those which are scaled versions of
the exceptional $p$-good maps. That is, they are obtained either by 
translating a twisted
map from a vertex which is a smaller exceptional map, or by
translating an exceptional map from a vertex. In what follows we
consider the case of exceptional maps coming from weights near the
right-hand edge of the dominant region --- the left-hand case is
entirely symmetric.
Since the exceptional $p$-good maps are generically
not composites of Carter-Payne
maps neither are the scaled versions of these maps. So to construct those maps
which cannot be found by factoring through the $G_1$-head of
$\nabla(\lambda)$ we must use some other method.

We first define the notion
of two weights $\lambda$ and $\xi$ being in \emph{exceptional
configuration} inductively.
A map from $\nabla(\lambda)$ to $\nabla(\xi)$ is then an
\emph{exceptional map} if the two weights $\lambda$ and $\xi$ are in
exceptional configuration.

First, the weights labelled by $1$ and $3$ in
Figures \ref{walls}a(i) and $1$ and $2$ in \ref{walls}c(i) are in
an \emph{exceptional configuration}. (These are the weights that
give the exceptional $p$-good maps.)

We say that two Steinberg weights $\lambda$ and $\xi$ are in an
\emph{exceptional configuration} if they are in the same $G$-block and
they come from two smaller weights in exceptional configuration.
That is we have $\lambda= p^d \lambda_1 + (p^d -1)\rho$ and $\xi = p^d \xi_1
+ (p^d -1)\rho$ and $\lambda_1$ and $\xi_1$ are in exceptional configuration
with $d \in \NN$ and $\lambda_1$ and $\xi_1$ both non-Steinberg
weights.  

\begin{figure}[ht]
\centerline{\epsffile{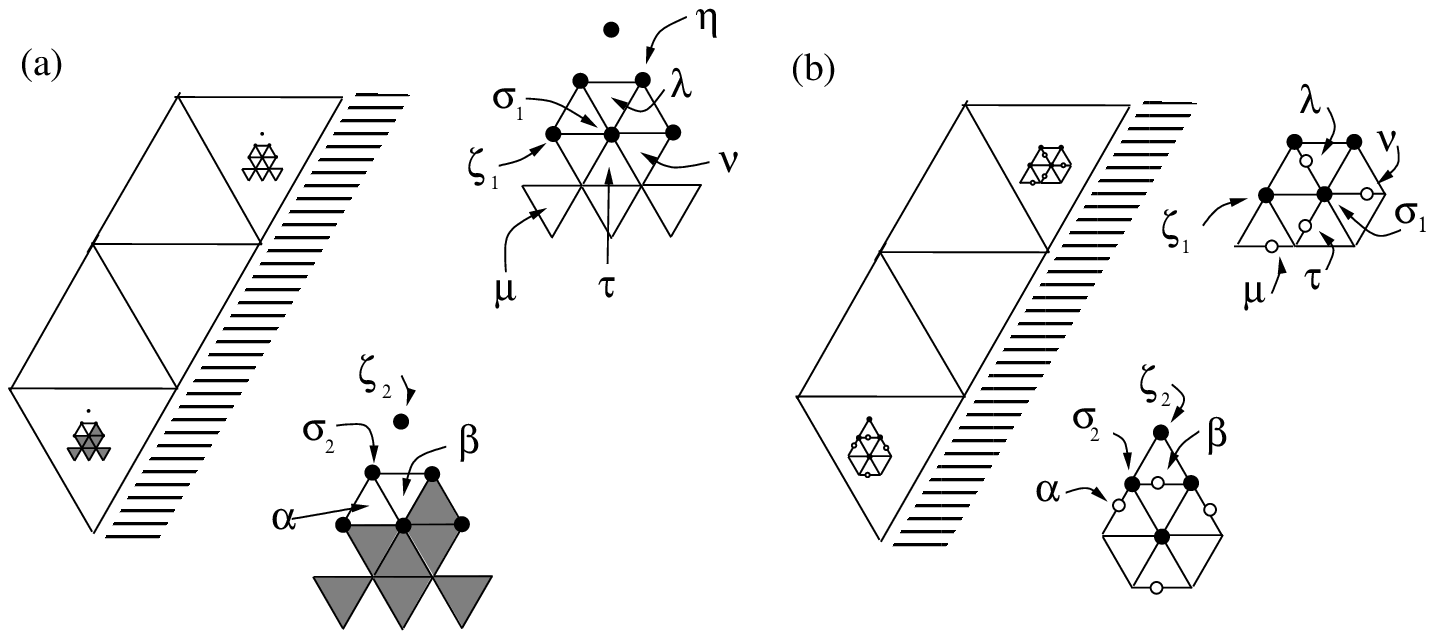}}
\caption{\label{exceptcase1}}
\end{figure}

\begin{figure}[ht]
\centerline{\epsffile{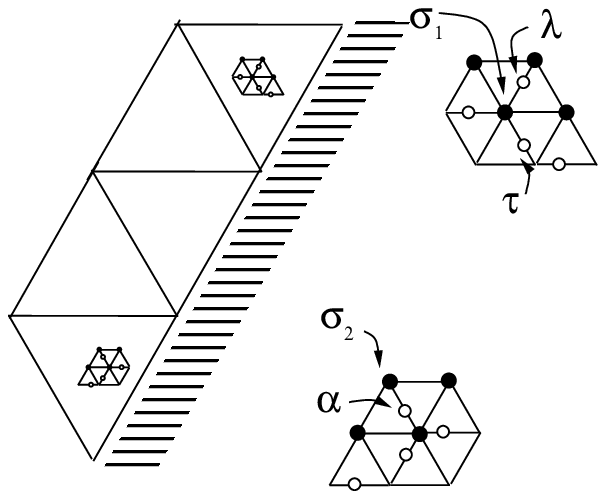}}
\caption{\label{exceptwallr}}
\end{figure}

\begin{figure}[ht]
\centerline{\epsffile{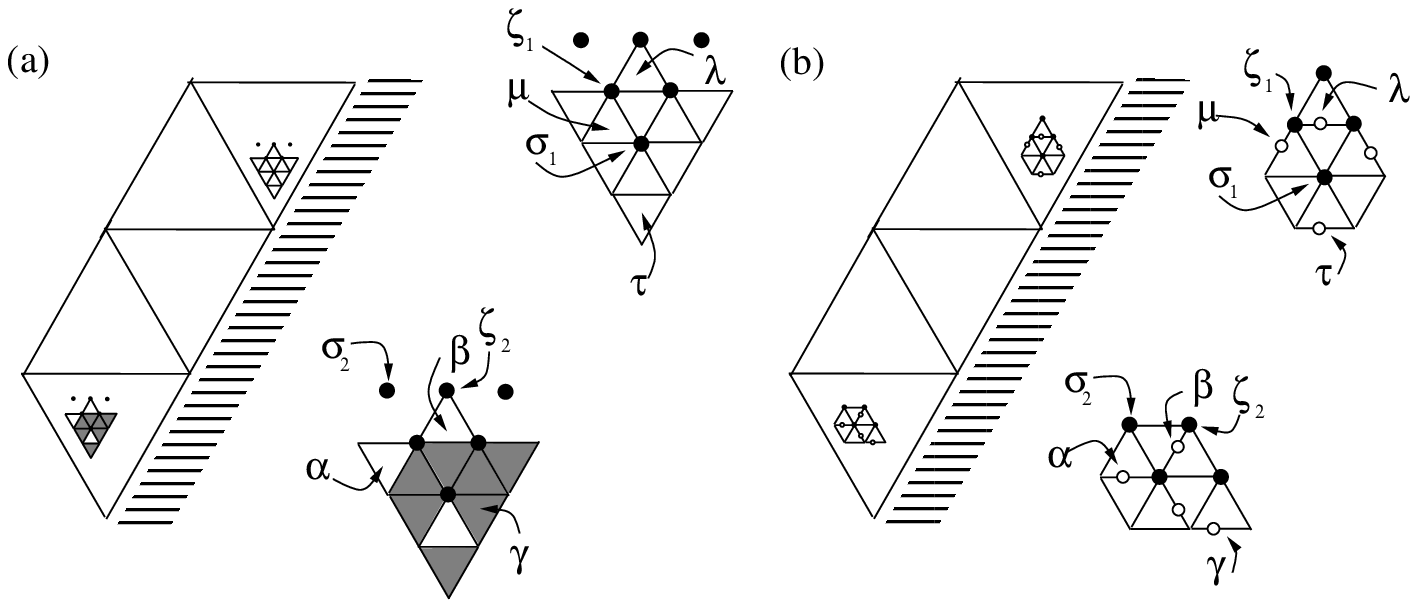}}
\caption{\label{exceptcase3}}
\end{figure}

\begin{figure}[ht]
\centerline{\epsffile{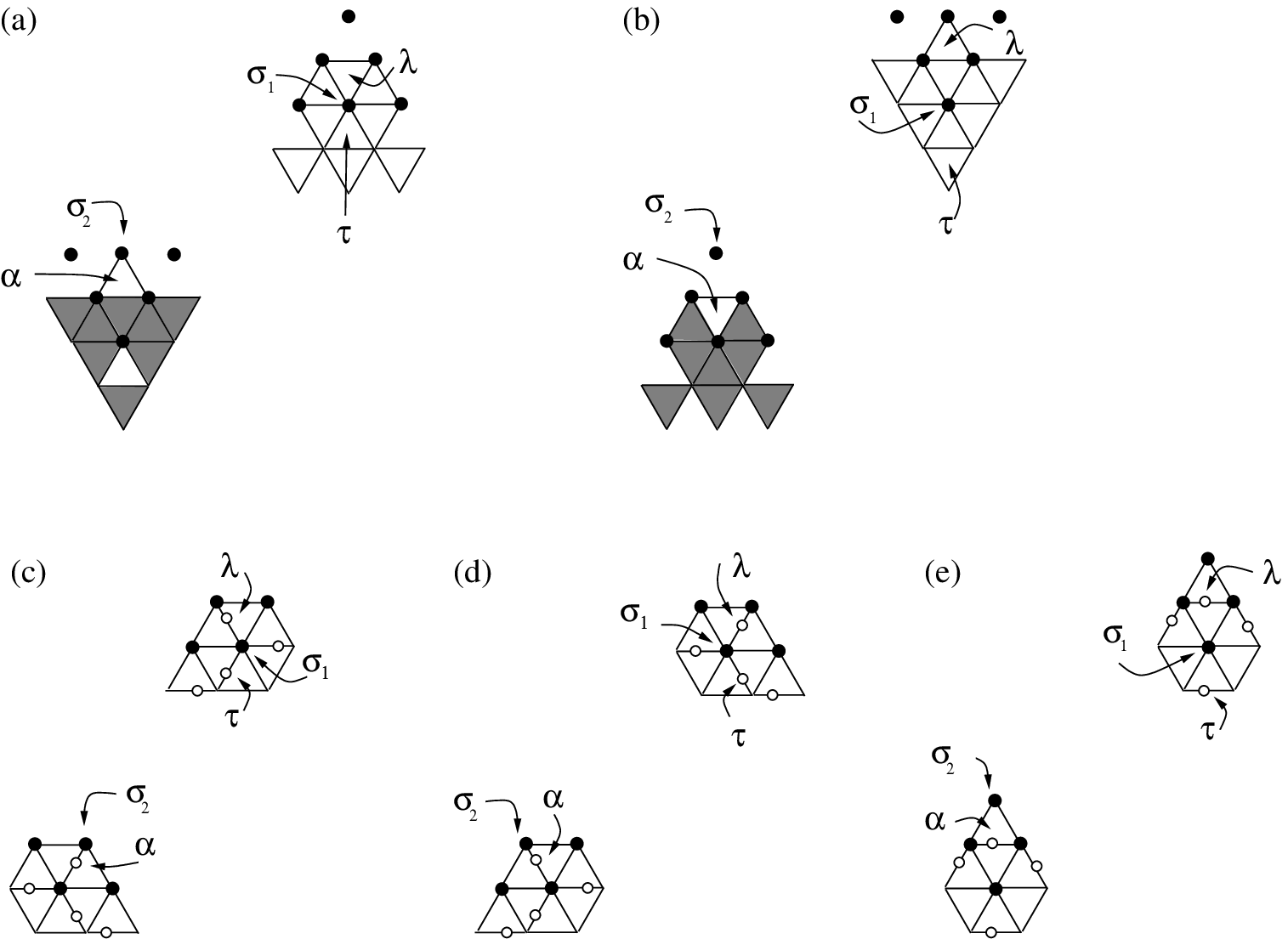}}
\caption{\label{exceptcase4}}
\end{figure}

We say that two non-Steinberg weights $\lambda$ and $\xi$ are in an
\emph{exceptional configuration} if the Steinberg weight $\sigma_1$ directly
below $\lambda$ and $\sigma_2$ with $\sigma_2  < \sigma_1$ are in exceptional
configuration and $\xi$ is an eligible $\sigma_2$ translate.
This is illustrated for various cases in 
Figures \ref{exceptcase1}, \ref{exceptwallr}, \ref{exceptcase3} and
\ref{exceptcase4} where the weights $\lambda$ and $\alpha$ 
and $\lambda$ and $\beta$ are in exceptional configuration.
These figures represent the generic situation and the weights are
related by large $p^d$-reflections. It is possible for
$\sigma_1''$ to lie on a wall or be a Steinberg weight. If $\sigma_1''$ is
a Steinberg weight then it must lie in the interior of larger
$p^d$-alcove and then $\sigma_1$ and $\sigma_2$ are related by
$p^d$-reflections for some $d >1$.
If $\sigma_1''$ is on a wall then it is possible for the $p$-filtration
pattern of $\nabla(\lambda)$ to go into the other $p^d$-alcove. 

If there is no Steinberg weight below $\lambda$ then we say $\lambda$
and $\xi$ are in an \emph{exceptional configuration} if there is
a homomorphism from $\nabla(\lambda)$ to $\nabla(\xi)$ (these were
determined in the previous section) but this is not constructed by a
composite of Carter-Payne maps. That is, it can only be obtained by
factoring through the $G_1$-head of $\nabla(\lambda)$. 

Thus the exceptional configurations are just scaled versions of the
exceptional $p$-good maps and hence always have the same basic
configuration as in Figures 
\ref{exceptcase1}, \ref{exceptwallr}, \ref{exceptcase3} and
\ref{exceptcase4}
(and the dual versions for the other edge of the dominant region).


Before we construct these maps, we note 

\begin{corollary}
Assumption \ref{assume} holds.
\end{corollary}
\pf This follows from Theorems \ref{downhomclass} and
\ref{twistclass}, and the possible configurations coming from
exceptional maps as considered above.\qed

Thus we have completed the proof of Theorem \ref{onedim}, and it only
remains to construct the exceptional maps. 

%
\begin{theorem}\label{makeex}
If $\lambda$ and $\xi$ are in an exceptional configuration then there is
a homomorphism from $\nabla(\lambda)$ to $\nabla(\xi)$.
\end{theorem}

\pf
If $\lambda$ (and hence $\xi$) is $p^2$-restricted then the only
possible exceptional configuration is exactly that which gives the
exceptional $p$-good homomorphisms. Thus the theorem is true in this
case.
We have also proved this theorem (essentially by definition) if
$\nabla(\lambda)$ is as in Section \ref{symway}.

We now take the statement of the theorem as our inductive hypothesis
for the rest of this section.

By left/right symmetry and the results from Section \ref{symway}, 
it remains to consider the case where
$\lambda$ and $\xi$ are in an exceptional configuration 
with $\sigma_1$ below $\lambda$ and we have a homomorphism
from $\nabla(\sigma_1)$ to $\nabla(\sigma_2)$ with eligible
$\sigma_2$-translates $\alpha$ and $\beta$ as shown. (Note that for
the remainder of this section we have adopted a non-standard labelling
of the weights; in particular $\tau$ is the weight usually denoted
$\tilde{\lambda}$.)

By inspection we see that there are several possible type of cases that can 
occur, as illustrated in the figures. 
The lower alcove and lefthand wall case have a similar weight
configuration and the upper alcove and horizontal wall case also has a
similar configuration. 

In all the cases there is an obvious map from $\nabla(\lambda)$ to
$\nabla(\alpha)$. It is the one induced by the map from
$\nabla_p(\tau)$ to $\nabla_p(\alpha)$. 
It remains to construct the map from $\nabla(\lambda)$ to
$\nabla(\beta)$ as in Figures \ref{exceptcase1} and \ref{exceptcase3}.


\begin{lemma}
Assume $\sigma_1,$ $\sigma_2$ , $\zeta_1$ and $\zeta_2$ are as in
Figures \ref{exceptcase1} and \ref{exceptcase3}. Then
$\Hom(\nabla(\sigma_1), \nabla(\sigma_2)) \ne 0$ implies that
$\Hom(\nabla(\zeta_1), \nabla(\zeta_2)) \ne 0$. 
\end{lemma}

\pf Clearly, by (\ref{homreduce}), this is equivalent to proving that
$\Hom(\nabla(\sigma_1''), \nabla(\sigma_2'')) \ne 0$ implies that
$\Hom(\nabla(\zeta_1''),\nabla(\zeta_2'')) \ne 0$. 

%
%
Now $\zeta_1''$ and $\sigma_1''$ are so close together that one must
be in the closure of the facet containing the other. Most of time they
will be in the same facet.  Firstly if $\zeta_1''$ is in the closure
of the facet containing $\sigma_1''$ then the translation principle
gives the desired result.  So we may assume that $\sigma_1''$ is in
the closure of the facet containing $\zeta_1''$.  We now claim that if
$\sigma_1''$ and $\sigma_2''$ are in an exceptional configuration then so
are $\zeta_1''$ and $\zeta_2''$ except for one case (where we will
construct the map directly).  Then our inductive hypothesis gives the
desired result.

To see the claim we observe that if $\zeta_1''$ is not in the same
facet as $\sigma_1''$ then $\sigma_1''$ is either a Steinberg weight
and $\zeta_1''$ is on a wall or $\sigma_1'''$ is on a wall and
$\zeta_1''$ is in an up alcove. We now look at Figures
\ref{exceptcase1} and \ref{exceptcase3} and we see that $\zeta_1''$
and $\zeta_2''$ are in an exceptional configuration.  The
only thing that can go wrong is if $\sigma_1$ is on a left-hand wall
(is $\alpha_1$-singular) and this is the same wall used to get the
exceptional configuration. Although $\zeta_1''$ and $\zeta_2''$ are
not then in an exceptional configuration, we can construct the map
from $\nabla(\zeta_1'')$ to $\nabla(\zeta_2'')$ directly: it is the
map $\phi^d_L \phi^d_R$ which is non-zero using Lemma \ref{lr}.  \qed

Before returning to the proof of Theorem \ref{makeex}, we review the general
theory that we will require. Suppose we have some module category with
non-trivial maps $\phi: A \to B$ and $\psi: C \to D$ and extensions
\begin{gather*}0 \to A \to E \to C \to 0\\
0 \to B \to F \to D \to 0.
\end{gather*}
We want the maps $\phi$ and $\psi$ to glue together to give a map
from $E \to F$.
This happens exactly when the push out of $E$ is the pull back of $F$.
That is we need $H:=\phi E \cong F \psi $ (using the notation of
MacLane\cite[chapter 3]{maclane}). Then
the following diagram commutes
$$\xymatrix@R=15pt@C=15pt{
0 \ar@{->}[r]  &  A \ar@{->}[r]
\ar@{->}[d]_{\textstyle{\phi}}
& E \ar@{->}[r] \ar@{->}[d]_{\textstyle{\bar{\phi}}} &  C \ar@{->}[r]
\ar@{->}[d]_{\textstyle{Id}} & 0 
\\
0 \ar@{->}[r]  &  B \ar@{->}[r] \ar@{->}[d]_{\textstyle{Id}} 
& H \ar@{->}[r] \ar@{->}[d]_{\textstyle{\bar{\psi}}} &  C \ar@{->}[r]
\ar@{->}[d]_{\textstyle{\psi}}& 0 
\\
0 \ar@{->}[r]  &  B \ar@{->}[r] & F \ar@{->}[r] &  D \ar@{->}[r] & 0
}$$
and the composite $\theta: E\to H \to F$ is non-zero. Moreover its image
has the following short exact sequence 
\begin{equation}\label{imseq}
0 \to \im \phi \to \im \theta \to \im \psi \to 0.
\end{equation}

We first consider the down alcove and left hand wall case (depicted in Figure
\ref{exceptcase1})
and start off by assuming that $\mu$ and $\nu$ (in $W.\lambda$) 
in the figure are in
different $\nabla_p$-classes.

In this case we take $A=\nabla_p(\mu)$, $B=\nabla_p(\beta)$,
$C=\nabla_p(\tau)$ and $D=\nabla_p(\alpha)$. We know that there is an
extension of $A$ by $C$ at the top of the $p$-filtration of
$\nabla(\lambda)$ (by our assumption on $\mu$ and $\nu$), and of $B$
by $D$ at the bottom of the $p$-filtration for $\nabla(\beta)$. We
denote these extensions by $E$ and $F$ respectively. Thus it remains
to show that $\phi E\cong F\psi$.

First consider the diagram defining the pushout $E'$ of $E$ along $\phi$:
$$\xymatrix@R=15pt@C=15pt{
0 \ar@{->}[r]  &  \nabla_p(\mu) \ar@{->}[r]
\ar@{->}[d]_{\textstyle{\phi}}
& E \ar@{->}[r] \ar@{->}[d]_{\textstyle{\bar{\phi}}} &  \nabla_p(\tau)
\ar@{->}[r] \ar@{->}[d]_{\textstyle{Id}} & 0 
\\
0 \ar@{->}[r]  &  \nabla_p(\beta) \ar@{->}[r] 
& E' \ar@{->}[r] &
\nabla_p(\tau)  \ar@{->}[r] & 0.}$$
This is exactly the definition of the connecting homomorphism
$\partial$ in the long exact sequence obtained by applying
$\Hom(-,\nabla_p(\beta))$ to the defining short exact sequence for $E$:
$$0\to
\Hom(\nabla_p(\tau),\nabla_p(\beta))\to\Hom(E,\nabla_p(\beta))\to
\Hom(\nabla_p(\mu),\nabla_p(\beta)) \stackrel{\partial}{\to} 
\Ext^1(\nabla_p(\tau),\nabla_p(\beta)).$$
As $\nabla_p(\tau)$ and $\nabla_p(\beta)$ are in different
$\nabla_p$-classes, the first Hom-space is zero. Further, $E$ sits at
the top of a filtration of $\nabla(\lambda)$, and hence has simple
head in the $\nabla_p$-class of $\nabla_p(\tau)$. Therefore the second
Hom-space is also zero, and $\partial$ is an embedding. By our
inductive hypothesis and (\ref{homreduce}) we have
$\Hom(\nabla_p(\mu),\nabla_p(\beta))\cong k$.

Consider $\Ext^1(\nabla_p(\tau),\nabla_p(\beta))$. By the results in
\cite[Lemma 4.2]{par1} and \cite[Proposition 3.3.2]{yehia} (reproduced
in \cite[Proposition 4.1]{par1}) we have 
\begin{equation}\label{exceptext}
\begin{array}{ll}\Ext^1(\nabla_p(\tau),\nabla_p(\beta))&
\cong\Hom(\nabla(\tau''),\nabla(\beta'')\otimes\nabla(0,1))\\
&\cong\Hom(\nabla(\sigma_1''),T_{\beta''}^{\sigma_2''}\nabla(\beta'')).
\end{array}
\end{equation}
This final Hom-space is one-dimensional (by our inductive hypothesis)
either because 
$T_{\beta''}^{\sigma_2''}\nabla(\beta'')\cong \nabla(\sigma_2'')$ (if
$\beta''$ is not on a wall) or because
$$\Hom(\nabla(\sigma_1''),T_{\beta''}^{\sigma_2''}\nabla(\beta''))\cong
\Hom(T^{\beta''}_{\sigma_2''}\nabla(\sigma_1''),\nabla(\beta''))
\cong\Hom(\nabla(\epsilon),\nabla(\beta''))$$
for some $\epsilon$ on the wall containing $\mu''$.

Thus any non-zero homomorphism from
$\nabla_p(\mu)$ to $\nabla_p(\beta)$ pushes out to the unique non-zero
extension
$$0\to\nabla_p(\beta)\to E'\to \nabla_p(\tau)\to 0$$
as $\partial$ is an isomorphism.

The case for the pullback is similar. We have the diagram
$$\xymatrix@R=15pt@C=15pt{
0 \ar@{->}[r]  &  \nabla_p(\beta) \ar@{->}[r]
\ar@{->}[d]_{\textstyle{Id}}
& F' \ar@{->}[r] \ar@{->}[d]_{\textstyle{\bar{\psi}}} &  \nabla_p(\tau)
\ar@{->}[r] \ar@{->}[d]_{\textstyle{\psi}} & 0 
\\
0 \ar@{->}[r]  &  \nabla_p(\beta) \ar@{->}[r] 
& F \ar@{->}[r] &
\nabla_p(\alpha)  \ar@{->}[r] & 0}$$
and corresponding connecting homomorphism
$$0\to
\Hom(\nabla_p(\tau),\nabla_p(\beta))\to\Hom(\nabla_p(\tau),F)\to
\Hom(\nabla_p(\tau),\nabla_p(\alpha)) \stackrel{\partial}{\to}
\Ext^1(\nabla_p(\tau),\nabla_p(\beta)).$$ 

By arguments as above the first pair of Hom-spaces are both zero, and
$\partial$ becomes an isomorphism from $k$ to $k$. Thus any non-zero
homomorphism from $\nabla_p(\tau)$ to $\nabla_p(\alpha)$ pulls back to
the unique non-zero extension $F'\cong E'$ and the composite map
$\bar{\psi}\circ\bar{\phi}$ from $E$ to $F$ is non-zero by
(\ref{imseq}). This gives the required homomorphism 
$$\nabla(\lambda)\twoheadrightarrow E\stackrel{
\bar{\psi}\circ\bar{\phi}}{\to}F\hookrightarrow \nabla(\beta).$$

When $\mu$ and $\nu$ are in the same $\nabla_p$-class we have to
modify the above argument, as
the extension $E$
above does not sit at the top of the $p$-filtration of
$\nabla(\lambda)$. This problem can be rectified by taking
$A$ to be the extension of $\nabla_p(\mu)$ by $\nabla_p(\nu)$ occurring
in the $p$-filtration for $\nabla(\lambda)$, and $E$ to be the
extension of $A$ by $\nabla_p(\tau)$ occurring at the top of this
$p$-filtration. Note that $A$ has been described in Lemma
\ref{nongen};  using this  we deduce that there exists a weight
$\eta$ on the same wall as $\sigma_1''$ such that
$$\begin{array}{ll}\Hom(A,\nabla_p(\beta))
&\cong \Hom(T^{\mu''}_{\eta}\nabla(\eta),\nabla(\beta'')) \\
&\cong \Hom(\nabla(\eta),T_{\mu''}^{\eta}\nabla(\beta''))\\
&\cong \Hom(\nabla(\eta),\nabla(\rho))
\end{array}$$
where $\rho$ is a weight on the same wall as $\sigma_2''$.  By the
translation principle this latter Hom-space in one-dimensional. 
The argument now goes through as
before, with the new choices for $A$ and $E$. 

The argument for Figures \ref{exceptwallr} and \ref{exceptcase3} is
the same in the generic case where $\alpha$ and $\gamma$ are in
different $\nabla_p$ classes. 
When $\alpha$ and $\gamma$ are in the same $\nabla_p$-class then we
need to construct a map from the extension $\nabla_p(\mu)$ by
$\nabla_p(\tau)$ to the submodule of
$\nabla(\beta)$ which has uniserial $p$-filtration by 
$\nabla_p(\beta)$, $\nabla_p(\gamma)$ and $\nabla_p(\alpha)$.
This can be constructed in the similar way as for the down alcove case,
using Lemma \ref{nongen} as before.

As every exceptional configuration is of one of the forms depicted in Figures
\ref{exceptcase1}, \ref{exceptwallr}, \ref{exceptcase3} or \ref{exceptcase4}
(or their left-handed analogues) we have constructed all
remaining maps, which completes our induction for this section and our
classification of homomorphisms.



\section{Examples}\label{examples}

In this section we give some examples for $p=3$ of how we may apply
the various theorems in the paper to give all the possible
homomorphisms starting at a particular induced module.

In the diagrams an alcove is shaded if it contains a weight that is a
composition factor of the initial induced module, $\nabla(\lambda)$.
This shading is dark if there is a homomorphism from $\nabla(\lambda)$
to the induced module corresponding to the weight in the alcove. That
is $\Hom(\nabla(\lambda), \nabla(\mu)) \cong k$ where $\mu$ is in the
$W_p$-orbit of $\lambda$ and inside the shaded alcove.  The
lightly shaded alcoves are those which contain weights $\mu\in
W.\lambda$ which label composition factors but for which there is no
non-zero homomorphism from $\nabla(\lambda) \to \nabla(\mu)$.

\begin{figure}[ht]
\centerline{\epsffile{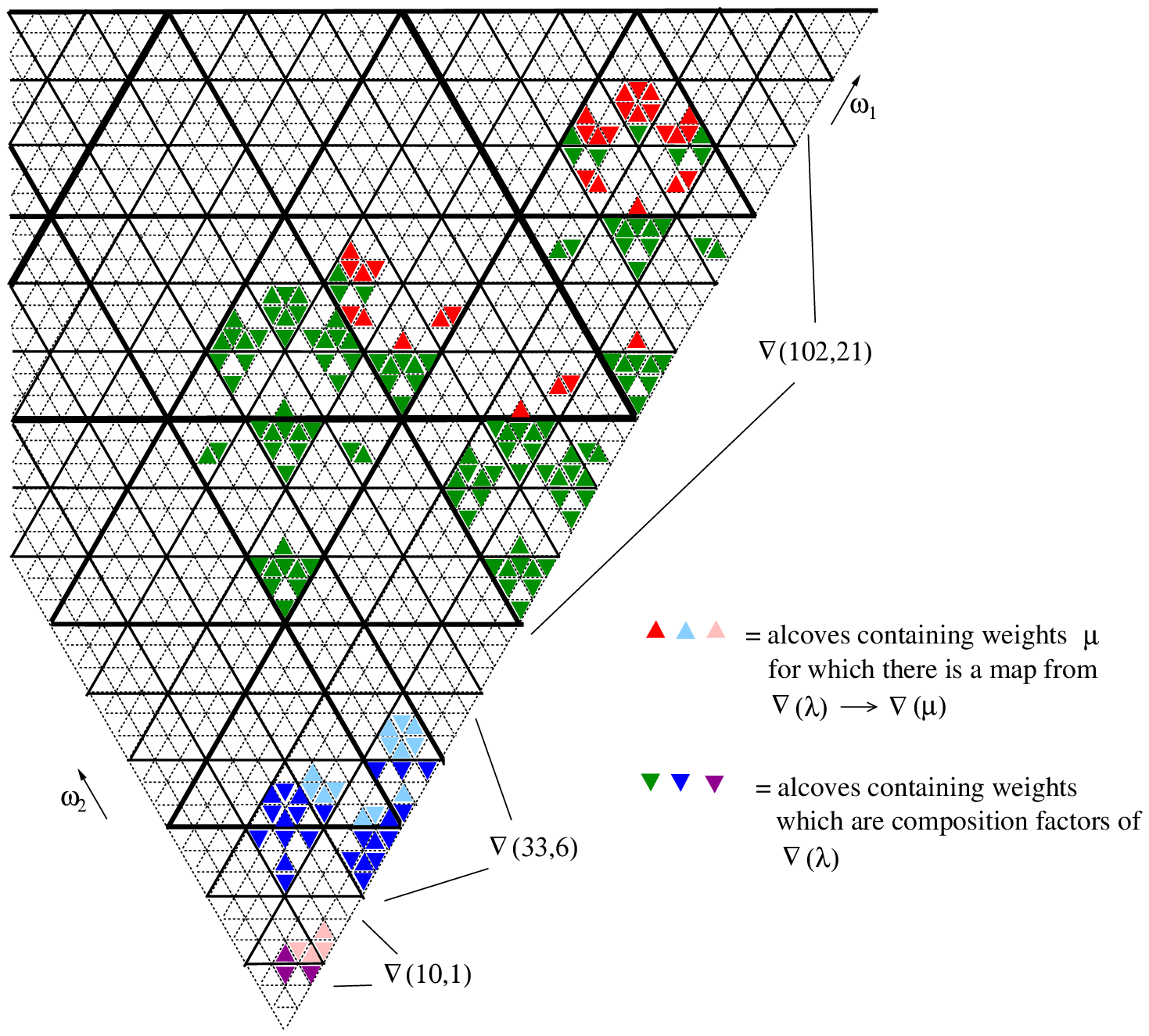}}
\caption{\label{prime3homs}}
\end{figure}

In Figure \ref{prime3homs} we determine all non-zero maps from
$\nabla(10,1)$, $\nabla(33,6)$ and $\nabla(102,21)$.     
We start with
$\nabla(10,1)$. The weight $(10,1)$ is in the interior of a just
dominant up alcove, and $(10,0)$ lies on the horizontal wall just
below it. Thus we will apply Theorem \ref{justuphoms} with this pair
of weights playing the role of $\lambda$ and $\nu$ respectively.
Using Theorem \ref{symprop}, 
we only get two non-zero homomorphisms starting from
$\nabla(10,0)$: the identity map and the map $\nabla(10,0) \to
\nabla(6,1)$. We then apply Theorem \ref{justuphoms} to obtain the
four darkly shaded alcoves shown.
(Alternatively it is not hard to see that the only
homomorphisms we get in this case are the $p$-good maps.)  Thus we have
(after transforming the alcoves back into weights)
$$\Theta_{(10,1)}=\{(10,1), (9,0), (6,3), (7,1)\}.$$
We can now determine the result for $\nabla(33,6)$. The vertex
immediately below $(33,6)$ is $3(10,1) + (2,2)$. So we follow the
procedure in Section \ref{homstart} and use Theorem \ref{downhomclass}
to obtain the dark shaded alcoves
as shown. That is we take the set $p \Theta_{(10,1)} + (p-1,p-1)$ and
then shade in the alcoves according to Figure \ref{downeligible}.
We get
\begin{multline*}
\Theta_{(33,6)}=\{(33,6), (34,4), (33,3), (31,7), (30,6), (31,4), 
(28,1), (19,13), (18,12), (19,10), \\ (21,9), (22,4), (24,4)\}.
\end{multline*}

We then can determine the result for $\nabla(102,21)$ using Theorem
\ref{downhomclass} and Theorem \ref{twistclass}.  Rather than give the
complete list of weights here, we depict them graphically in Figure
\ref{prime3homs}.  In these cases the only types of maps we get are
composites of Carter--Payne maps and twisted maps.

\begin{figure}[ht]
\centerline{\epsffile{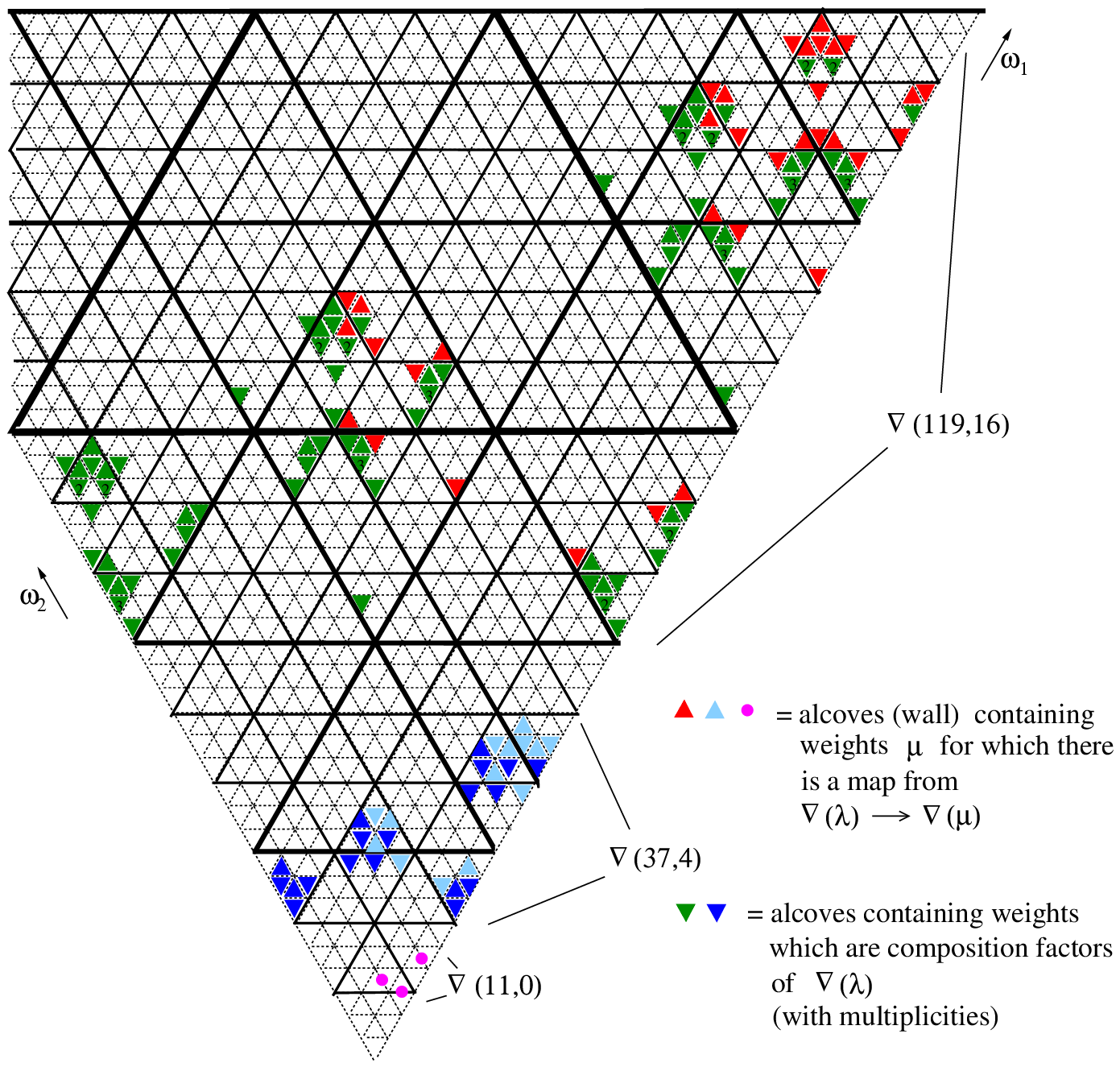}}
\caption{\label{prime3homs2}}
\end{figure}

To illustrate a case involving walls we include Figure \ref{prime3homs2},
where we determine all maps from 
$\nabla(11,0)$, $\nabla(37,4)$ and $\nabla(119,16)$.
The reader may verify that the procedure in the
preceding sections gives rise to the set of maps shown.
The case of $\nabla(37,4)$ is an example where the vertex is not 
$p^2$-regular.
These examples give all the different types of maps. Note that
the exceptional $p$-good map $\nabla(11,0) \to \nabla(6,1)$ gives rise
to the exceptional configurations $((37,4), (18,3))$ and
$((37,4), (22,1))$ which in turn give rise to the
exceptional configurations $((119,16), (54,9))$, $((119,16), (66,3))$ 
and $((119,16), (70,1))$.



\section{Homomorphisms for the symmetric group}\label{rdsect}

In this section we will show how the results for $\SL_3$ can be used to
classify homomorphisms between certain Specht modules for the
symmetric group.

%
Let $\Sigma_d$ be the symmetric group on $d$ letters. For each
partition $\lambda$ of $d$ we can explicitly define a Specht module
$S^{\lambda}$ for $k\Sigma_d$. (See for example \cite[Chapter 3]{jk}.)

Partitions with at most $n$ parts can also be used to label
representations of $\GL_n$. The representation theory of this group is
essentially the same as that for $\SL_n$; in particular given $n$
part partitions $\lambda$ and $\mu$ of $d$ we have that
$$\Hom_{\GL_n}(\nabla(\lambda),\nabla(\mu))\cong
\Hom_{\SL_n}(\nabla(\bar{\lambda}),\nabla(\bar{\mu}))$$ where for a
partition $\tau=(\tau_1,\tau_2,\ldots,\tau_n)$ we set
$\bar{\tau}=(\tau_1-\tau_2,\tau_2-\tau_3,\ldots,
\tau_{n-1}-\tau_n)$. Those representations of $\GL_n$ labelled by
partitions of $d$ are also representations of the corresponding Schur
algebra $S(n,d)$, and again the Hom-spaces are unchanged.

The representation theory of the general linear group and the
symmetric group are related. In particular we  
have the following theorem of Carter and Lusztig
\begin{theorem}[{\cite[Theorem 3.7]{carlus}}]\label{carlus}
If $p \ge 3$ and $\lambda$ and $\mu$ are partitions of $d$ with at
most $n$ parts then
$$
\Hom_{\GL_n}(\nabla(\lambda),\nabla(\mu))\cong
\Hom_{k\Sigma_d}(S^{\lambda},S^{\mu}).
$$
\end{theorem}
%
%
We thus have
$$\Hom_{k\Sigma_d}(S^{\lambda},S^{\mu})\cong
\Hom_{\SL_n}(\nabla(\bar{\lambda}), \nabla(\bar{\mu})).$$
Combining this with  our earlier results gives

\begin{corollary}
Suppose that $p\ge 3$, $d>0$ and $\lambda$ and $\mu$ are three part
partitions of $d$. 
Then $\Hom_{k\Sigma_d} (S^{\lambda},S^{\mu})$ can be
determined from the results in the preceding sections.
\end{corollary}


 
\section{Levi factors and a theorem of Fayers and Lyle}

In this final section we will show how our $\SL_3$ results can be
applied to determine certain Hom-spaces in higher rank cases, and show
how the result of Carter and Lusztig relates work of Donkin to 
a tensor product theorem
for Hom-spaces due to Fayers and Lyle \cite{fayly}. In each case, our
main tool will be a theorem of Donkin concerning certain Hom-spaces
for Levi factors of a reductive group $G$, which we begin by
recalling.

Let $G$ be a  reductive algebraic group as in Section \ref{prelim}. 
Given a subset $I$ of the simple roots $S$, we may also consider
representations of the corresponding Levi factor $G_I$ of $G$. To
distinguish such modules from those for the original group we shall use
the subscript $I$. We have the following result of Donkin (see
\cite[(4.3) Corollary]{erd}).

\begin{theorem}\label{donred} If $\lambda$ and $\mu\in X^+$ are such that
$\lambda-\mu\in\zed I$ then
$$\Ext_G^i(\nabla(\lambda),\nabla(\mu))\cong
\Ext_{G_I}^i(\nabla_I(\lambda),\nabla_I(\mu))$$
for all $i\geq 0$.
\end{theorem}

%

Now let $G$ be $\GL_n$ with the usual choice of simple roots as in
\cite[II, 1.21]{jantzen}. Suppose that $\lambda$ and $\mu$ are two
partitions of $d$ with at most $n$ parts (which we regard as elements
of $X^+$) and that we can write $\lambda=(\lambda(1),\ldots,\lambda(t))$ and
$\mu=(\mu(1),\ldots,\mu(t))$, where for each $i$ the elements
$\lambda(i)$ and $\mu(i)$ are both
partitions of $d_i$ into at most $n_i$ parts. Then it is
easy to verify that $\lambda - \mu$ is an element of $\zed I$, where
$I$ is the subset of $S$ giving rise to the Levi factor $G_I\cong
\GL_{n_1}\times\cdots\times \GL_{n_t}$.

As noted in \cite[(4)]{donnote} we have that
$\nabla_I(\lambda)\cong \nabla_1(\lambda(1))\boxtimes\cdots\boxtimes
\nabla_t(\lambda(t))$ (and similarly for $\nabla_I(\mu)$), where
$\nabla_i(\lambda(i))$ is the induced module of weight $\lambda(i)$ for
$\GL_{n_i}$. Combining this with Theorem \ref{donred} and the
K\"unneth formula we obtain

\begin{corollary} \label{tenhom}
Let $m \in \NN$, $\lambda=(\lambda(1),\ldots,\lambda(t))$ and
$\mu=(\mu(1),\ldots,\mu(t))$ be as above. Then we have
\begin{multline*}
\Ext^m_{G_I}(\nabla(\lambda),\nabla(\mu))\\
\cong
\bigoplus_{m = i_1+\cdots +i_t}
 \Ext^{i_1}_{\GL_{n_1}}(\nabla_1(\lambda(1)),\nabla_1(\mu(1)))
\otimes\cdots\otimes
\Ext^{i_t}_{\GL_{n_t}}(\nabla_t(\lambda(t)),\nabla_t(\mu(t))).
\end{multline*}
\end{corollary}


\begin{corollary} If $\lambda$ and $\mu\in X^+$ differ only by an element of
$\zed I$, $I$ a subset of the simple roots 
and $I$ can be realised as a root system whose largest
connected component is of type $A_2$, then
the results of the preceding sections (together with the known result
for type $A_1$) can be applied to determine
$\Hom_G(\nabla(\lambda),\nabla(\mu))$. 
In particular, the $\Hom$-space is at most one-dimensional.
\end{corollary}

Combining Theorem \ref{carlus} and Corollary \ref{tenhom} we obtain
\begin{corollary} 
Let $\lambda$ and $\mu$ be as above and $p\ge 3$. We have
$$\Hom_{k\Sigma_d}(S^{\lambda},S^{\mu})\cong
\Hom{k\Sigma_{d_1}}(S^{\lambda(1)},S^{\mu(1)})\otimes\cdots\otimes
\Hom{k\Sigma_{d_t}}(S^{\lambda(t)},S^{\mu(t)}).$$
\end{corollary}

This is the theorem of Fayers and Lyle 
\cite[Theorem 2.2]{fayly}.
We also note here that Donkin \cite[Section 10, Proposition 4]{donkpreprint} 
uses Schur functors to prove a generalised version of the above
result for $\Ext$-groups (in small degree) for $p\ge 3$.

\begin{example}
Let $p=5$, $G= \GL_7(k)$, $\lambda =(90,63,40,8,4,4,4)$
and $\mu=(90,54,39,18,4,4,4)$.
We denote the simple roots for $G$ by $\alpha_i = \epsilon_i -
\epsilon_{i+1}$, where $\{\epsilon_i \}$ is the usual basis for $X$.
We have $\lambda -\mu = (0,9,1,-10,0,0,0) = 9 \alpha_2 +11 \alpha_3$.
Thus 
\begin{align*}
\Hom_G(\nabla(\lambda), \nabla(\mu))&\cong 
\Hom_{\GL_3(k)}(\nabla(63,40,8), \nabla(54,39,18))\\
&\cong
\Hom_{\SL_3(k)}(\nabla(23,32), \nabla(15,21)) \cong k.
\end{align*}
We also get $\Hom_{k\Sigma_{213}}(S^{\lambda}, S^{\mu}) \cong k$\\
\end{example}


\section*{Acknowledgements}

The authors would like to thank \O yvind Solberg for some helpful
suggestions used in Section \ref{exception}. We would also like to
thank the organisers of the LMS Durham Symposium: Representations of
Finite Groups and Related Algebras, and of the Conference on
Representation theory for Algebraic Groups at Aarhus University, where
part of this work was carried out.


\providecommand{\bysame}{\leavevmode\hbox to3em{\hrulefill}\thinspace}
\providecommand{\MR}{\relax\ifhmode\unskip\space\fi MR }
\providecommand{\MRhref}[2]{%
  \href{http://www.ams.org/mathscinet-getitem?mr=#1}{#2}
}
\providecommand{\href}[2]{#2}

\end{document}